# The general case of cutting of *GML* surfaces and bodies

Johan Gielis[1] and Ilia Tavkhelidze[2]

[1] University of Antwerp, Bio-Engineering Sciences, Belgium
e-mail johan.gielis@uantwerpen.be

[2] Iv.Javakhishvili Tbilisi,State University, Faculty of Exact and Natural Sciences, Georgia
e-mail ilia.tavkhelidze@tsu.ge

# 1 Propositio

## 1.1 The general problem of cutting regular GML bodies

Given that: a) the cross section of a Generalised Möbius-Listing body $GML_m^n$ (Figure 1) is a disk with regular polygons as boundary with $m$-symmetry, or $m$-vertices and $m$-sides, and b) full cutting along the complete structure is performed with a knife perpendicular to the cross section, dividing the convex polygon or disk in precisely two parts, c) with knives that cut from vertex to vertex ($VV$), vertex to side ($VS$) or side to side ($SS$), then determine 1) in how many ways can a $GML_m^n$ body be cut, and 2) the ways in which resulting shapes are linked.

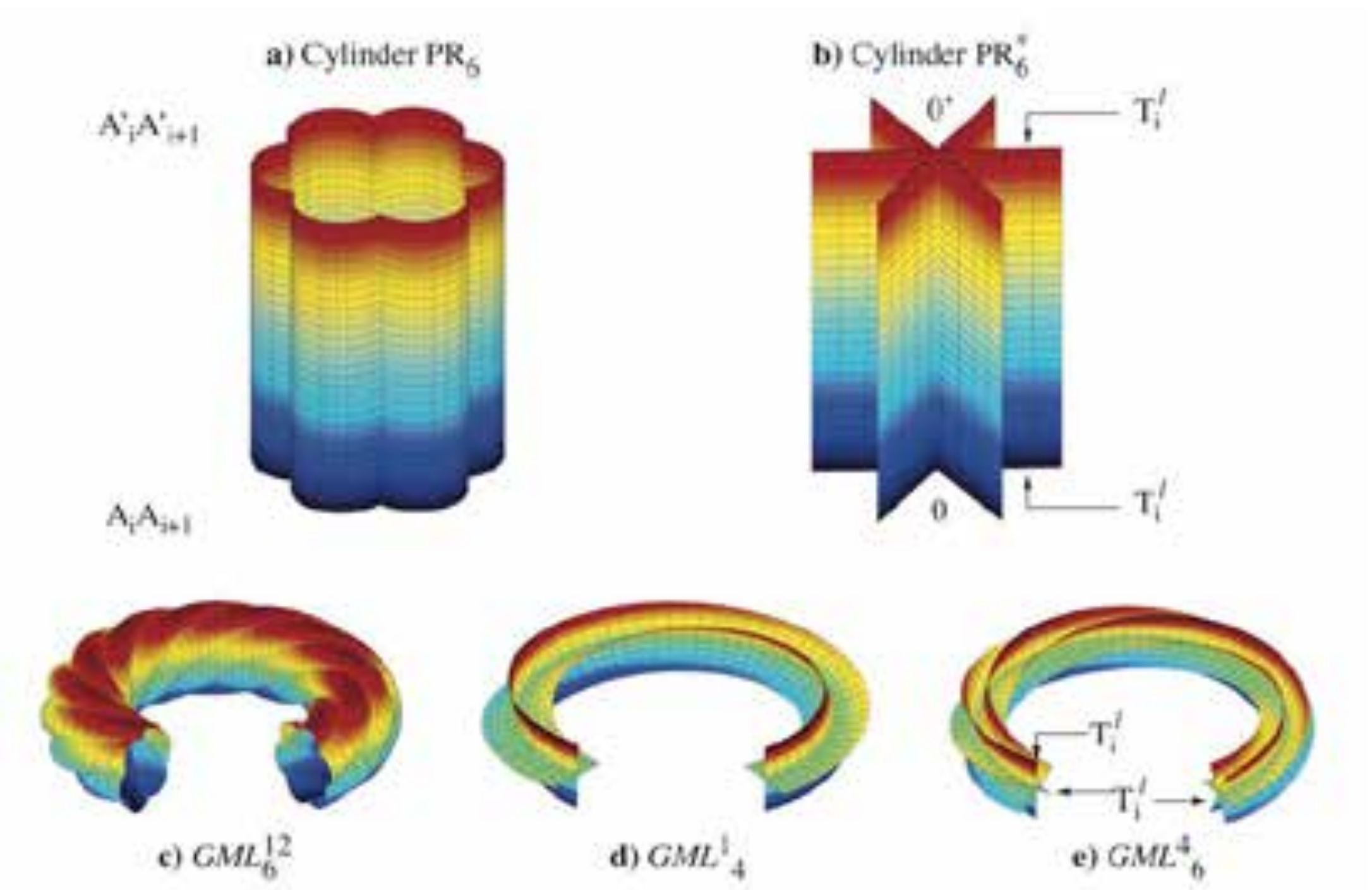

*Figure 1: Identification of vertices of prisms (**a**. $PR_6$ and **b**. $PR_6^*$ with $*$ starlike), leading to GML surfaces with twists (**c,d,e**). The number of twists relative to the symmetry of the cross sections is indicated by the superscript.*

Solutions for this problem have been obtained for $GML$ surfaces and bodies with $m = 2, 3, 4, 5, 6$ [1-6] and revealed connections with the field of knots and links, depending on the number of twists of the original cylinder or prism with $m$-symmetrical cross section, given by the superscript $n$ in $GML_m^n$. The method of study was the full cutting of $GML_m^n$ surfaces and bodies with a moving knife in 3D. More general solutions were obtained when the cross section of the $GML_m^n$ surface or body is a Gielis curve, describing the boundary of the cross section as well as the disk [4-7]. Gielis curves can transform the cross section of the $GML_m^n$ from a circle to a regular polygon and vice versa, and to many other concave or convex curves. Hence using regular polygons allows for generalizing results to other convex figures.

Instead of a moving knife and fixed $GML_m^n$, an equivalent approach is a fixed knife and a moving $GML_m^n$ surface or body. Furthermore, using a fixed knife it is shown that the solution for the general problem of cutting regular $GML_m^n$ surfaces and bodies can be obtained by studying the problem of cutting regular $m$-polygons in the plane, with $d_i$ knives, where $i$ is the ordinal number of the divisor and $d_i$ is the divisor or $m$.

For example, considering $VV$ cuts, a cut with a single knife $d_1$ is made from vertex $V$ to any other vertex $V$ dividing the polygon into two distinct parts. This can be repeated from vertex $V_1$ to $V_j$ , from $V_2$ to $V_{j+1}$ etc...., giving $m$ possible cuts. A $d_i$ knife with $i = m$ on the other hand is a knife with $m$ blades, achieving the same result as for $d_1$ but in one cut.

The focus is on $GML_m^n$ bodies since the cutting problem will lead to separate sectors in the cross sections, and to separate bodies in one or more resulting $GML$'s. This problem has both a geometrical and a topological solution. In the former case, the number and the precise shape of the resulting bodies is counted, while in the topological case, only the number of sides and vertices of the resulting bodies is counted. For example, when the resulting shape is a quadrilateral, in the topological solution all quadrilaterals are treated as one solution, whereas in the geometrical solution the precise shape of the quadrilateral is important, e.g. a square. The topological solution is a subset of the geometrical solution.

## 1.2 The geometrical solution

Defining $N_m^{VV,VS,SS}$ as the number of ways of cutting an $m$-gon for a given divisor of $m$, whereby the superscript details whether this cut is a $VV$, $VS$ or $SS$ cut, giving $N_m^{VV}$, $N_m^{VS}$, $N_m^{SS}$ respectively, then the total number of way of cutting with $d_{i=1,\dots,m}$ knives with $i = 1, \dots, m$ blades is the sum of all possible cuts multiplied by the number of divisors.

The solution can be expressed as a recurrence relation in number of side-to-side cuts. Defining $N_m$ as the number of ways of cutting an $m$-gon for a divisor equal to one, and $N_m^{SS}$ as the number of $SS$ cuts for $m$ and $N_{m-2}^{SS}$ for the number of $SS$ cuts for the polygon with $(m-2)$ symmetry, then for the geometrical solution the number of different ways of cutting an $m$-polygon with a $d_m$ knife is:

- For even $m$ $(= 2k)$: $N_{m=2k} = m + 1 + N_{m-2}^{SS}$ (1a)
- For odd $m$ $(= 2k+1)$: $N_{m=2k+1} = m + 2 + N_{m-2}^{SS}$ (1b)

The total number is then the number of different ways in (1a) and (1b) multiplied by the number of divisors. In regular polygons the cuts $N_m^{VV}$, $N_m^{VS}$, $N_m^{SS}$ have to be made from $V_i \to V_{i+1}$, from $V_i \to S_{i+1}$, or from $S_i \to S_{i+1}$ respectively. $V_i \to V_i$, $V_i \to S_i$, or $S_i \to S_i$ do not divide the regular $m$-polygons into two distinct parts. As a corollary, in a convex polygon or a circle with equally spaced points however, such cuts are possible and the total number of cuts then has to be augmented by $m$.

- For even $m$ $(= 2k)$: $N_{m=2k} = 2m + 1 + N_{m-2}^{SS}$ (2a)
- For odd $m$ $(= 2k+1)$: $N_{m=2k+1} = 2m + 2 + N_{m-2}^{SS}$ (2b)

## 1.3 The topological solution

The topological solution is given by the following general formula:

- If $m = 2k + 1$ and has $N$ nontrivial divisors $d_2, d_3 \ldots d_{N+1}$ and $d_1 \equiv 1, d_{N+2} \equiv d_m \equiv m$, then the number of all possible variants of cutting of $GML_m^n$ bodies is

$$\Xi_m^{top} = 8k + 1 + 3Nk + \sum_{i=2}^{N+1} \left[\frac{k}{d_i}\right] + 2N \qquad (3a)$$

- If $m = 2k$ and has $N$ nontrivial divisors $d_2, d_3 \ldots d_{N+1}$ and $d_1 \equiv 1, d_{N+2} \equiv d_m \equiv m$, then the number of all possible variants of cutting of $GML_m^n$ bodies is

$$\Xi_m^{top} = 8k - 5 + 3Nk + \sum_{i=2}^{N+1} \left[\frac{k-1}{d_i}\right] \qquad (3b)$$

The general formula can be expressed in the variables $k$ and $N$, whereby $m = 2k$ for even or $m = 2k + 1$ for odd numbers. $d_i$ is the $i$-th ($i = 2,3, \ldots N + 1$) non-trivial divisor of the number $m$, so $N$ equals the number of divisors of $m$ *minus* 2 (excluding the trivial divisors $d_{N+2} \equiv d_m \equiv m$, and $d_1 \equiv 1$). $N_{div}^m = N + 2$ is the total number of divisors (including the trivial divisors $m$ and *1*) and square brackets [..] indicate the integer part of the fraction. Also in this case the remark on convex polygons and circles is valid: the number of possible cuts has to be increased by $m$.

## 1.4 The occurrence of Möbius phenomena

The Möbius phenomenon, which led to the discovery of non-orientable surfaces for classical ribbons or strips, also occurs in *GML* surfaces and bodies. Möbius phenomena result in only one body or surface after full process of cutting when a) $m$ = even and b) the knife passes through the centre of the cross section of *GML* bodies. The knives used in making $VV, VS$ or $SS$ cuts are *chordal* knives, dividing the $m$-polygon into two distinct parts and cutting the boundary of the $m$-polygon in exactly two points.

When using *radial* knives, which start from the centre of the polygon and cut the boundary of the $m$-polygon in exactly one point, it turns out that Möbius phenomena can occur both for $m$ even and odd.

## 1.5 Generalizations and interrelations

This paper is organized in the traditional Greek geometrical way, with *Propositio, Expositio, Determinatio, Constructio, Demonstratio* and *Conclusio*. In the Expositio and Determinatio sections the focus is on paths towards generalizations and links to various well-studied problems in mathematics.

Generalizations include the circle or convex polygons, whereby $V_i \to V_{i+1}, V_i \to S_i$, or $S_i \to S_i$ cuts are possible. Furthermore, it is by no means necessary to use straight (chordal or radial) knives as long as the knife-curve is wholly contained within the original domain. The precise shape of the knife can be the solution to some optimization problem. One generalization is cutting of concave polygons, whereby the knife is part of a “more-concave” polygon.

In cutting the polygon or convex shape, the resulting domains can be separated. However, the process of the knives is completely equivalent to connecting vertices and sides via drawing diagonals or connecting equally spaced points on a circle. Hence the problem of using $m$-knives to cut regular $m$-gons has a direct relation to other well studied problems in geometry, such as:

- The problem of dividing $m$-polygons with non-crossing diagonals
- Number of separate parts in a circle when connecting all equally distributed points on a circle with straight lines.

The former problem was studied first by Euler and later by Segner, Lamé and Catalan. The solution of the problem is combinatorial and the results are the Catalan numbers. The combination of geometric and combinatorial problems (linked directly to permutations or identification of corresponding vertices (Figure 1) in the problem of cutting *GML*) is very fruitful. Solutions to the second and related problems are other series of integer numbers, many of which are found in the Online Encyclopedia of Integer Sequences *OEIS*. There are other relationships and many applications in physics and biology, and it opens a way forward to bridge topology and geometry within maths on the one hand, and mathematics and the natural sciences on the other.

# 2 Expositio

## 2.1 *GML* surfaces and bodies

$GML_m^n$ are torus-like surfaces or bodies, whereby the cross sections of the $GML_m^n$ bodies are closed planar curves with symmetry $m$, and which are constructed by identifying opposite sides of a cylinder or prisms (Figure 1) [1-4]. The planar curves with symmetry $m$ can be regular polygons or any closed plane curve, including circles. *GML* surfaces are generated when only the curve itself, as boundary of a region is considered (Figure 2 for $GML_4^n$), or bodies when also the disk enclosed by the curve is considered. For the classical cylinder or Möbius band, the cross section is a line, swept along a path forming a ribbon and twisted an even or odd number of times respectively. Whereas the lower index $m$ determines the symmetry of the cross section, the upper index $n$ describes the number of twists. *GML* surfaces or bodies can either be closed (Figure 1) or not (Figure 2). In the latter case they can be either sections of closed *GML* surfaces or bodies, or complete. They are called Generalized Rotating and Twisting surfaces and bodies *GRT*.

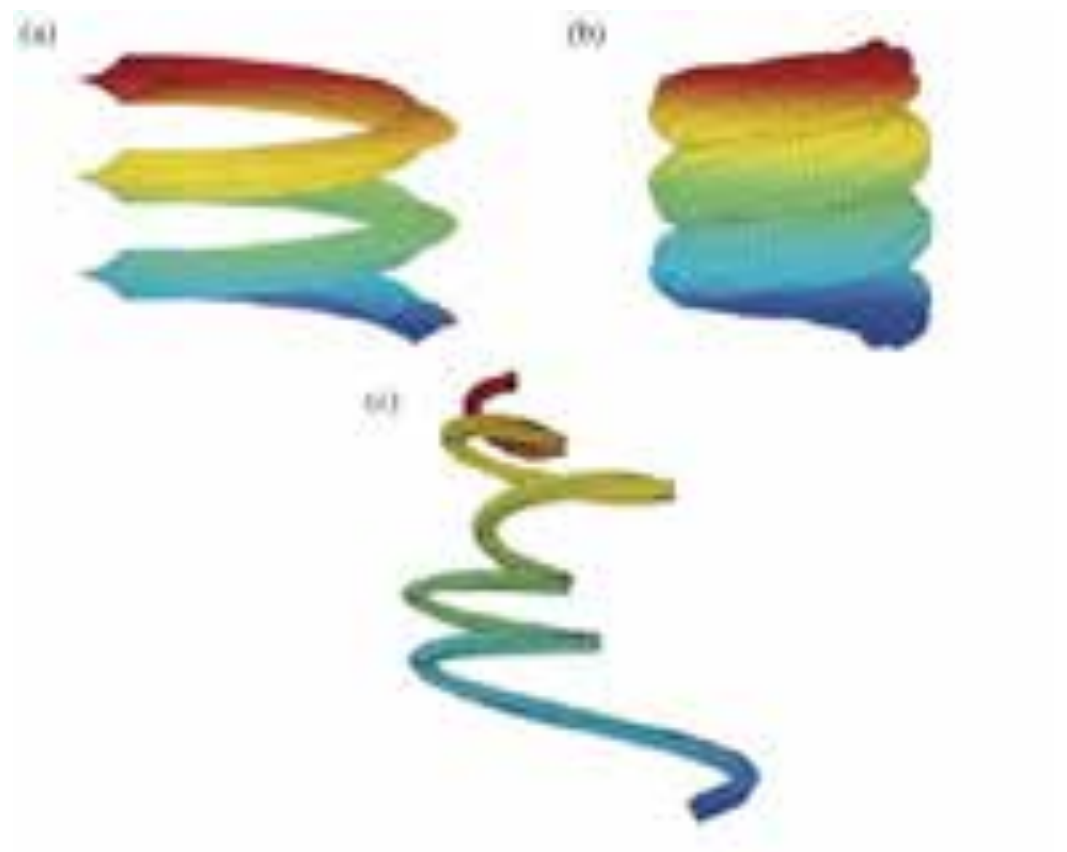

*Figure 2: GRT surfaces and bodies*

### 2.1.1 Analytic definition of *GML* surfaces and bodies

Generalized Möbius-Listing Bodies $GML_m^n$ are defined by the analytic representation:

$$\begin{cases} X(\tau,\psi,\theta) = \left(R(\theta) + p(\tau,\psi)\cos\left(\frac{n\theta}{m}\right) - q(\tau,\psi)\sin\left(\frac{n\theta}{m}\right)\right)\cos(\theta) \\ Y(\tau,\psi,\theta) = \left(R(\theta) + p(\tau,\psi)\cos\left(\frac{n\theta}{m}\right) - q(\tau,\psi)\sin\left(\frac{n\theta}{m}\right)\right)\sin(\theta) \\ Z(\tau,\psi,\theta) = K(\theta) + \left(p(\tau,\psi)\sin\left(\frac{n\theta}{m}\right) + q(\tau,\psi)\cos\left(\frac{n\theta}{m}\right)\right) \end{cases} \quad (4)$$

or, alternatively,

$$\begin{cases} X(\tau,\psi,\theta) = \left(R(\theta) + p(\tau,\psi)\cos\left(\psi + \frac{n\theta}{m}\right)\right)\cos(\theta) \\ Y(\tau,\psi,\theta) = \left(R(\theta) + p(\tau,\psi)\cos\left(\psi + \frac{n\theta}{m}\right)\right)\sin(\theta) \\ Z(\tau,\psi,\theta) = K(\theta) + p(\tau,\psi)\sin\left(\psi + \frac{n\theta}{m}\right) \end{cases} \quad (5)$$

where $p(\tau,\psi)$ and $q(\tau,\psi)$ are given functions which can define surfaces of the radial cross section of bodies (Figure 3), with:

- a planar curve or disc as basic cross section (this can be a very general section, esp. Gielis curves or discs), which
- can be swept around a closed or open path along a basic line (e.g surfaces of revolution, or Monge surfaces), and which
- can be given twists and turns as desired.

### 2.1.2 Permutations and Colours

Each Generalized Möbius Listing's body $GML_m^n$ can be considered as a geometric representation of the element of the permutation group [8], in particular:

$$\begin{pmatrix} A_0A_1 & A_1A_2 & \cdots & A_{i-1}A_i & \cdots & A_{m-2}A_{m-1} & A_{m-1}A_0 \\ A'_nA'_{n+1} & A'_{n+1}A'_{n+2} & \cdots & A'_{n+i-1}A'_{n+i} & \cdots & A'_{n+m-2}A'_{n+m-1} & A'_{n+m-1}A'_n \end{pmatrix} \equiv$$
$$= \begin{pmatrix} A_0A_1 & A_1A_2 & \cdots & A_{i-1}A_i & \cdots & A_{m-2}A_{m-1} & A_{m-1}A_0 \\ A'_{\omega m+k}A'_{\omega m+k+1} & A'_{\omega m+k+1}A'_{\omega m+k+2} & \cdots & A'_{\omega m+k+i-1}A'_{\omega m+k+i} & \cdots & A'_{\omega m+k+m-2}A'_{\omega m+k+m-1} & A'_{\omega m+k+m-1}A'_{\omega m+k} \end{pmatrix} \equiv$$
$$\equiv \begin{pmatrix} A_0A_1 & A_1A_2 & \cdots & A_{i-1}A_i & \cdots & A_{m-2}A_{m-1} & A_{m-1}A_0 \\ A'_kA'_{k+1} & A'_{k+1}A'_{k+2} & \cdots & A'_{k+i-1}A'_{k+i} & \cdots & A'_{k+m-2}A'_{k+m-1} & A'_{k+m-1}A'_k \end{pmatrix} \quad (6)$$

The sub-indices change from 1 to $m-1$ cyclically and the sides are identified by the definition of $GML_m^n$ body. So it may be written as follows:

$$\begin{pmatrix} A_0A_1 & A_1A_2 & \cdots & A_{i-1}A_i & \cdots & A_{m-2}A_{m-1} & A_{m-1}A_0 \\ A'_nA'_{n+1} & A'_{n+1}A'_{n+2} & \cdots & A'_{n+i-1}A'_{n+i} & \cdots & A'_{n+m-2}A'_{n+m-1} & A'_{n+m-1}A'_n \end{pmatrix} \equiv$$
$$\stackrel{def}{\equiv} \begin{pmatrix} A_0A_1 & A_1A_2 & \cdots & A_{i-1}A_i & \cdots & A_{m-2}A_{m-1} & A_{m-1}A_0 \\ A_nA_{n+1} & A_{n+1}A_{n+2} & \cdots & A_{n+i-1}A_{n+i} & \cdots & A_{n+m-2}A_{n+m-1} & A_{n+m-1}A_n \end{pmatrix} \equiv \quad (7)$$
$$\stackrel{not}{\equiv} \begin{pmatrix} 1 & 2 & \cdots & i & \cdots & m-1 & m \\ k & k+1 & \cdots & \mathrm{mod}_m(k+i) & \cdots & \mathrm{mod}_m(k+m-1) & \mathrm{mod}_m(k+m) \end{pmatrix}$$

The first row shows the identifications of the original prism and the second row the identifications after rotation and identification according to definition of $GML_m^n$. These are elements of a special subgroup of the group of permutations, in which the number of elements of this subgroup equals $m$. The first row is always the usual sequence of numbers from 1 to $m$, and the second row may start from any number from **1** to $\boldsymbol{m}$ and then, given the cyclicity of the representation of numbers, their sequence is not changed.

Each element, except the neutral element, has a certain number of cycles, and geometrically, this directly relates to how many different colors can be painted on the surface of a given body, without lifting the brush and not passing the ribs. For example for $m = 4$ (Figure 3):

$\begin{pmatrix} 1 & 2 & 3 & 4 \\ 1 & 2 & 3 & 4 \end{pmatrix}$ - neutral 4-colored ($m$-colored, for any $m$), Figures 3 **a.** & **d.**

$\begin{pmatrix} 1 & 2 & 3 & 4 \\ 2 & 3 & 4 & 1 \end{pmatrix}$, and $\begin{pmatrix} 1 & 2 & 3 & 4 \\ 4 & 1 & 2 & 3 \end{pmatrix}$ - one colored (Figures 3 **c.** & **f.**)

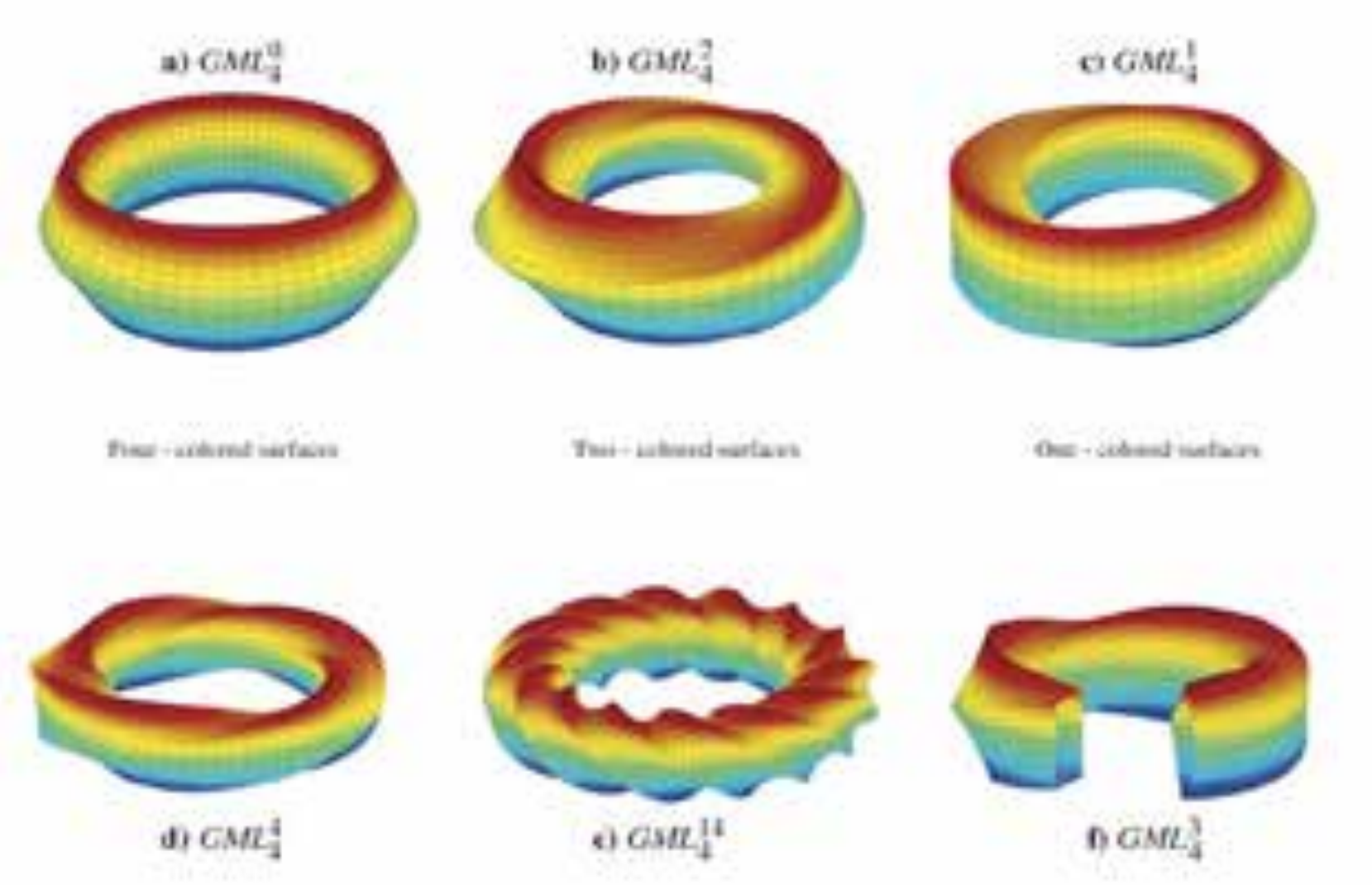

*Figure 3: $GML_4^n$ surfaces with different values of twisting n, leading to one-, two-, or four coloured surfaces.*

**Corollary:** Each $GML_m^n$ body has a $j$-colored surface where $j = \gcd(m, \kappa)$ and $\boldsymbol{n} = \boldsymbol{\omega m} + \boldsymbol{\kappa}$ except when $\kappa = 0$, in this case we have an $m$-colored $GML_m^n$ body [8].

### 2.1.3 From the original motivation to the general problem

The original motivation to study $GML_m^n$ is rooted in the study of boundary value problems and the geometrical description of natural shapes and phenomena [4]. The current study was initiated by the question what happens if the *GML* surfaces and bodies are cut along a certain line or surface, inspired by the cutting of the original Möbius band [2, 3]. For example, the result of cutting **d.**, **e.** and **f.** in Figure 3 along lines containing one vertex, will be very different.

Thus far, the classification of cutting of $GML_m^n$ surfaces and bodies was achieved for classic Möbius bands [3], and for $GML_m^n$ with $m$ = 2, 3, 4, 5 and 6 [4-6]. A full classification was achieved for cutting of classic Möbius bands with any $k$ number of knives [3]. These classifications revealed a close link between the cutting of *GML* bodies and surfaces, the study of knots and links, and with the colouring of surfaces (Figure 3). This can lead to intricate links of various separated *GML* surfaces and bodies (Figure 4). In this case the cutting leads to ribbons, but if the $GML_6^0$ Architon is a body the various structures will have different cross-sections.

The challenge remains to classify the cutting of general *GML* bodies when the cross section of the *GML* body is a regular $m$-polygon, for any value of $m$. This problem can be studied from several points of view, and it will be shown in particular that cutting planar regular $m$-polygons with a certain number of knives is an equivalent problem.

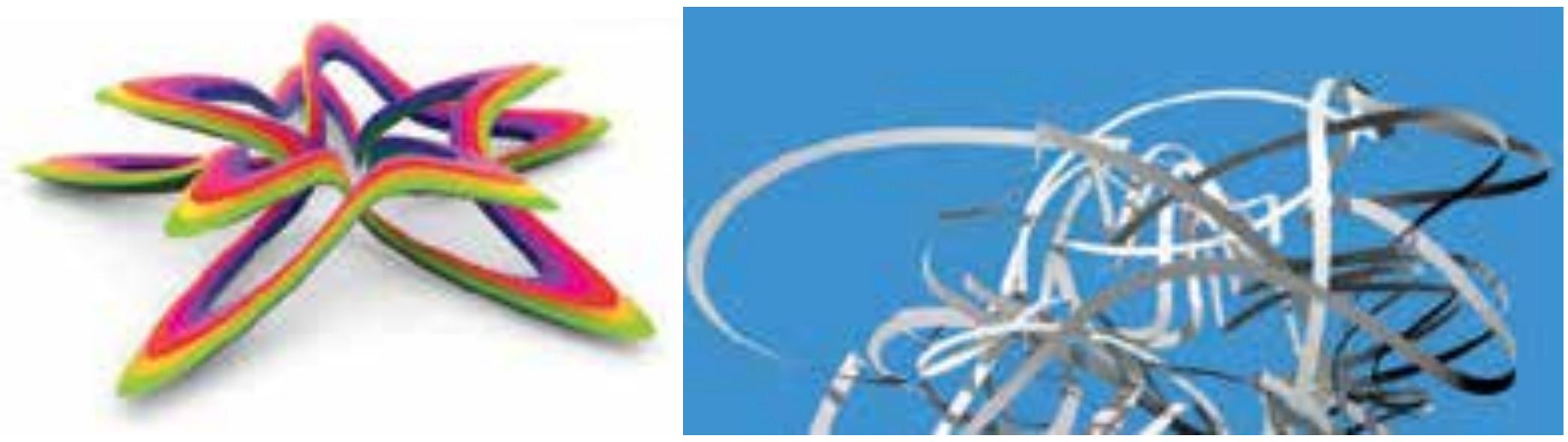

*Figure 4* ***a.*** $GML_6^0$ *surface (Architon) with 6 different sectors in different colours and,* ***b.*** *the result after separating the different sectors.*

Demonstrations will be given for cutting with one knife or with $m$ knives. Other ways of cutting, based on the divisors of $m$, are investigated, and the equivalence of cutting of $GML_m^n$, the cutting of regular $m$-polygons, or in general Gielis curves and rotational symmetries of the circle with inscribed $m$-polygons will be considered, and the relevant connections to other fields of mathematics will be outlined. Using $m$-polygons does not limit the generality of our results. More details are provided in Sections 2.2, 3. *Determinatio.*

## 2.2 Gielis curves, surfaces and transformations

### 2.2.1 Gielis curves and transformations

The theory of $GML$ surfaces and bodies was greatly enriched when it was merged with Gielis curves [7;9-11], a generalization of Lamé curves or superellipses ($A \neq B$) or supercircles ($A = B$) for any symmetry, defined as:

$$\varrho(\vartheta; A, B, n_1, n_2, n_3) = \frac{1}{\sqrt[n_1]{\left|\frac{1}{A}\cos\left(\frac{m}{4}\vartheta\right)\right|^{n_2} \pm \left|\frac{1}{B}\sin\left(\frac{m}{4}\vartheta\right)\right|^{n_3}}} \tag{8}$$

The parameter $m$ defines the symmetry and is thus related to the lower index of $GML_m^n$. The exponents $n_{1,2,3}$ are not related to the upper index $n$ of $GML_m^n$. Equation (8) can also be considered as a generic and geometric transformation on all planar functions, unifying a wide range of natural and abstract shapes [7]. Here we restrict the transformation on a constant function or circle, with $A = B$. Gielis curves (1) can serve both as cross section and as basic line around which the cross section is swept giving rise to $GML$ bodies (Figs 3 and 4). They can define:

- Boundaries and disks, using the zones enclosed by the curves, technically using $\varrho(\vartheta) \leq$ whereby the = sign defines the boundary.
- Regular polygons with symmetry $m$ are defined:

$$\varrho(\vartheta) = \lim_{n_1 \to \infty} \left[ \frac{1}{\sqrt[n_1]{\left|\cos\left(\frac{m}{4}\vartheta\right)\right|^{2(1-n_1 \log_2 \cos\frac{\pi}{m})} + \left|\sin\left(\frac{m}{4}\vartheta\right)\right|^{2(1-n_1 \log_2 \cos\frac{\pi}{m}}}} \right] \tag{9}$$

- Regular Gielis polygons sensu Matsuura [12] are generated for $G_{m,p,q} = G_{m,n_{2,3}=p,n_1=q;} = G_{m,1,(m/4)^2}$. A $G_{5,1,\frac{25}{16}}$ polygon is a very good approximation of a regular pentagon with $m = 5, n_{2,3} = p = 1$ and $q = p.(m/4)^2$. The difference between regular polygons defined in this way is less than 1% for $m \geq 5$ and less than 0.5% for $m \geq 11$.

- Self-intersecting polygons ($m$-polygrams) for $m$ a rational number, are known as Rational Gielis Curves *RGC* [13]; for example for $m = \frac{5}{2}$, a pentagram is generated. To complete the figure with 5 vertices, 2 rotations are needed, corresponding to the numerator and denominator of $m$ respectively. In general, $m$ can be rational ($m = {}^{p}/_{q}$ with *p, q* relative prime), or irrational [7].

Furthermore, it is a continuous transformation, whereby any shape can be transformed into a circle, and then into any other shape [11]. The transformation is equal to 1, yielding the circle, for either $m$=0 or $lim_{m\to\infty}$ The former is a zero-gon (or zero-angle), a figure without corners or vertices, while the latter is the classic notion of a circle as the limiting case of a polygon with an infinite number of sides. But the value is also equal to 1, and thus a circle, for $n_2 = n_3 = 2$, for **any** integer value of $m$ (given $A = B = 1$). So if a pentagon is transformed into a circle, this can be done by changing the values of $n_2 = n_3$ for the pentagon to $n_2 = n_3 = 2$ (in a discrete or continuous way). However, even if this is a circle, it is also a pentagon, with the original equi-spacing of the five points on the circle still imprinted. The polygons can be transformed into circles, where the equidistant points are the roots of unity.

### 2.2.2 Pythagorean compact

Gielis curves are generalizations of Lamé curves using polar coordinates, and for $n_{1,2,3} = n$ and $n_{1,2,3} = 2$ (given $A = B, m = 4$), we have Lamé curves and the classic Euclidean circle, respectively. Hence, these transformations have been named Pythagorean-compact. As all shapes are described in one Pythagorean compact equation, all shapes are equally simple, differing in a few numbers at most. Within the same Pythagorean structure, the arguments of the cosine and sine functions (in the original form $\frac{m}{4}\vartheta$) may be arbitrary functions $f(\vartheta)$, allowing for extremely compact form descriptions [14].

Moreover, all shapes can be continuously (or discontinuously if one chooses discrete steps) transformed into any other Gielis curve. The Lorentz transformation of Special Relativity Theory is one special case [10-11]. There are deep connections to various other parts of mathematics, including Riemann surfaces, approximation theory and number theory. For the natural sciences, the main consequence is that the Gielis Formula allows for a uniform description of a wide variety of natural shapes, and their development [11].

### 2.2.3 Gielis surfaces and bodies

The two dimensional case can be generalized to [11]:

$$r(\vartheta,\varphi) = \frac{1}{\sqrt[n_1]{\left|\frac{\sin\left(\frac{m_1\vartheta}{4}\right)\cos\left(\frac{m_2\varphi}{4}\right)}{a}\right|^{n_2} + \left|\frac{\sin\left(\frac{m_1\vartheta}{4}\right)\sin\left(\frac{m_2\varphi}{4}\right)}{b}\right|^{n_3} + \left|\frac{\cos\left(\frac{m_1\vartheta}{4}\right)}{c}\right|^{n_4}}} . \alpha(\vartheta,\varphi) \quad (10)$$

As transformations of functions, with $\alpha(\vartheta,\varphi)=1$ the unit sphere, Gielis transformations can be generalized immediately to any dimension [9]. Another way of defining 3D shapes is parametrically using two Gielis curves $\rho_1(\vartheta), \rho_2(\varphi)$:

$$\begin{cases} x = \varrho_1(\vartheta)\ \cos\vartheta\,.\,\varrho_2(\varphi)\cos\varphi \\ y = \varrho_1(\vartheta)\ \sin\vartheta\,.\,\varrho_2(\varphi)\cos\varphi \\ z = \varrho_2(\varphi)\ \sin\varphi \end{cases} \qquad (11)$$

with a unit sphere for $\rho_1(\vartheta) = \rho_2(\varphi) = 1$.

In this way a 3D Gielis surface or body can be defined on the basis of two perpendicular Gielis curves $\varrho_1(\vartheta), \varrho_2(\varphi)$. There are a variety of ways of defining surfaces and bodies, with Figure 4 as one example. In Figure 5 other examples are shown, showing both the equivalence with GML and special topological figures.

a. $GML_4^1$ body with one full rotation (circle as basis line)
b. $GML_0^0$ with cross section a circle or zero-angle, and basic line a ½ angle, ($m = p/q = 1/2$).
c. $GML_0^0$ with basic line a one-angle or monogon ($m = 1$), no twists.

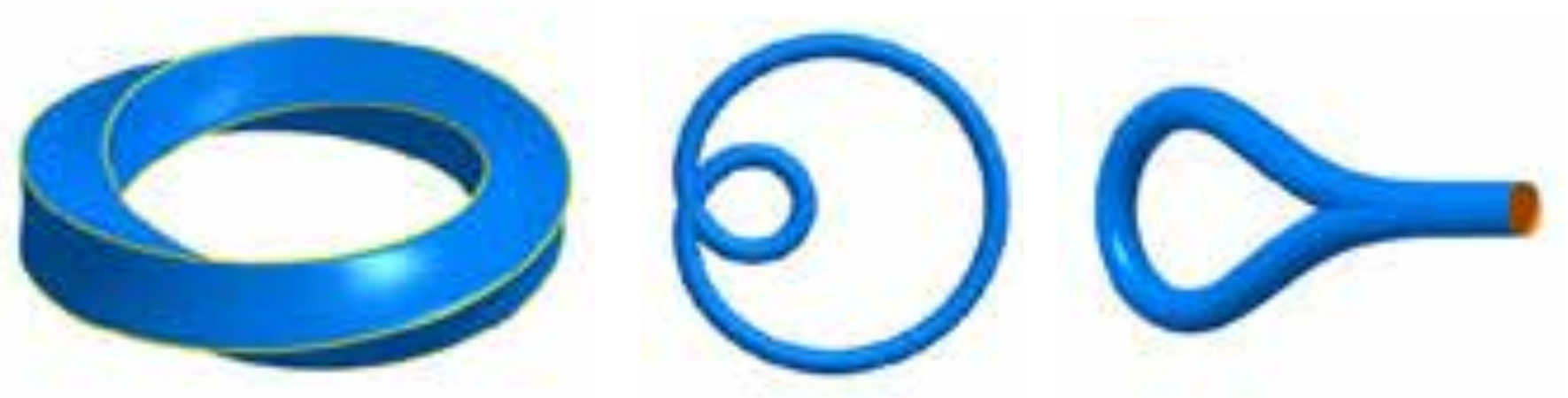

**a** **b** **c**

*Figure 5* ***a.*** *$GML_4$ surface.* ***b****. Torus with a Gielis curve with m=1/2 as path.* ***c.*** *A torus with path m=1, giving rise to a Klein bottle*

Figure 5**c.** is a monogon with parameter $m = 1, A = B$ and all exponents $n$ in Eq.1 equal to 1. When this 2D curve is transformed into a torus, it becomes the Klein bottle. It is remarked that in all cases the cross sections are constant, but also they can be scaled or morphed when swept along the basic line, resulting for example in vortices. In 3D a sphere can be transformed into a torus in a purely geometrical way, since the 3D parametric version is based on two perpendicular sections [11].

When transforming a sphere into a torus one of the cross sectional circles, increases in size. In Figure 5**b.** a torus transforms into the half-angle shape ($m = \frac{1}{2}$), which leads to a torus with two holes. Actually, by having value of $m$ in the Gielis formula as $1/p$ with $p$ an integer, tori with any number of holes equal to $p$ can be generated. The number of holes corresponds to the genus of a surface. One example is the ½ angle torus in Figure 5**b.**, which has two holes, or genus 2. In Figures 1, 3 4 and 7, the basic line of the *GML* surfaces and bodies is a circle, but this can be any knot-like structure or RGC, such as a pentagram.

## 2.3 The equivalence of cutting *GML* bodies and *m*-polygons

### 2.3.1 GML meets Flatland

At the end of the 19th century E. Abbott wrote *Flatland, a story of dimensions* [15]. The goal was to show how people could think about four-dimensional space, by considering Flatland and its inhabitants. When a ball moves through Flatland, the inhabitants will initially see a circle that starts as a point, increases in size, until the maximum size is reached. Then it will decrease until it disappears again, leaving the Flatlanders in shock and awe. Only when one of the Flatlanders travels into 3D space, (s)he realizes what truly happened, but ultimately finds $h_{im}^{er}$self in great difficulty trying to explain this to fellow inhabitants of Flatland.

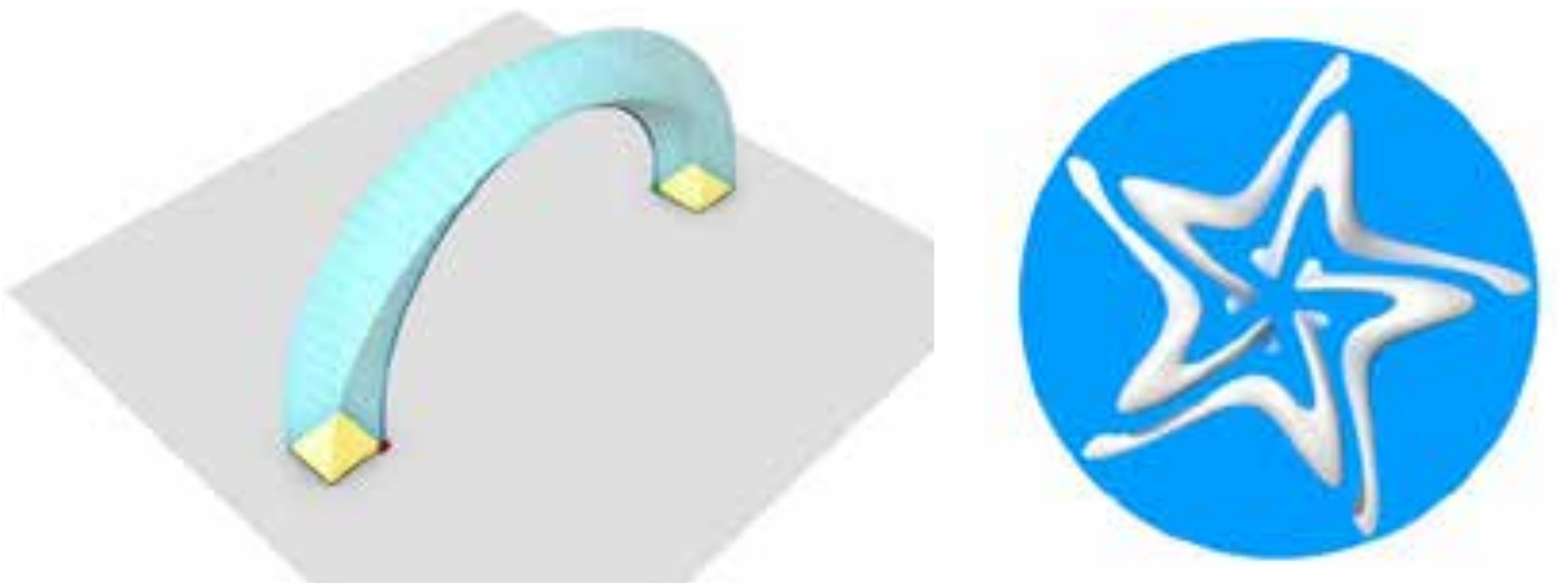

*Figure 6* ***a.*** *Upper half of GML meeting Flatland, with the green and the red dot connected via a rib.* ***b.*** *Architon (Figure 4) moving through a cone.*

We can take the same approach. A *GML* body or surface is a 3D torus-like surface or body, with certain symmetry. Fig. 6**a.** displays part of a *GML* body with square symmetry. The *GML* surface or body crosses Flatland (Figure 6**a.**) and the yellow zones are those observed by inhabitants of Flatland. However, they cannot discriminate between the different orientations of the square, resulting from a twist in 3D and indicated by red and green dots. Figure 6**b.** displays a snapshot of *Architon* (Figure 4) moving through Coneland. To a Conelander various unconnected and mysterious shapes are observed. Undoubtedly to a Flatlander or Conelander, a whole range of mysterious things happens if such structures move through their territories. What the Flatlanders or Conelanders will notice – a boundary or a disk – depends on whether *GML* is a surface or a body. In advanced societies like Flatland and Coneland, its physicists may develop methods to find out whether the shapes are boundaries or disks.

Figure 7**a.** displays part of a complete *GML* body with pentagonal cross-section, cut from side to side, cutting the cross sections into two distinct parts. If the cutting process is continued along the whole twisted *GML* body until the knife arrives at the very same position as the initial one, eventually four different bodies will result, two triangular ones in light blue and brown, and two pentagonal ones, in yellow and grey. Each body or surface will be twisted in a certain way.

Figure 7**b.** displays only one of these shapes, but in fact, after cutting, all four structures are intertwined into one complex shape. An inhabitant of Flatland will only see the overall pentagonal cross section as in Figure 7**a.**, which is cut from side $S_i$ to $S_{i+2}$ with five knives, showing the equivalence of cutting *GML* bodies and regular polygons. Alternatively (s)he may see a number of disconnected shapes, one with a pentagonal shape, ten shapes with one specific triangular shape (either blue as in Figure **7b** or brown) and five quandrangular shapes.

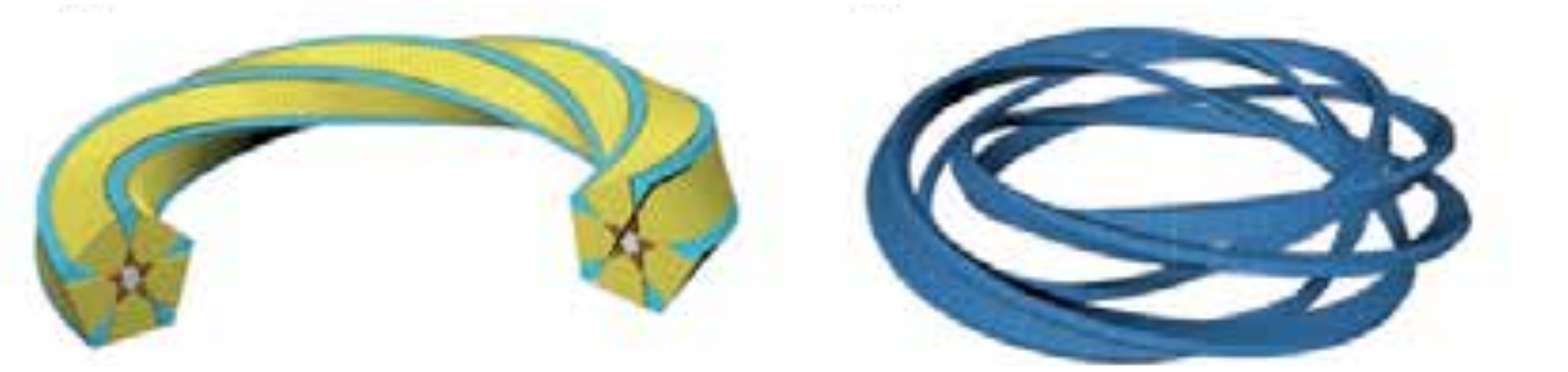

*Figure 7* ***a.*** *Pentagonal GML body.* ***b.*** *One of the resulting structures.*

### 2.3.2 Cutting bamboo poles.

When cutting *GML* with one (or more) knives, one can consider the knife moving along a rib on the *GML* body as in Fig 6**a**. In this case the knife corresponds to one of the diagonals of the cross sections, keeping a constant orientation relative to the cross section, e.g. along one rib of the *GML* body. Alternatively one can fix the knife and move the body through the knife. A *real life* example of cutting with various knives is found in hand tools or in machines used for splitting bamboo culms (Figure 8**a** and **b**). The knives can be adjusted to split the bamboo in any number of pieces. These culms are hollow, but one can easily imagine the same procedure for full prisms. In case the prisms are twisted and closed in a *GML* fashion, the results with one knife or more knives will be the same.

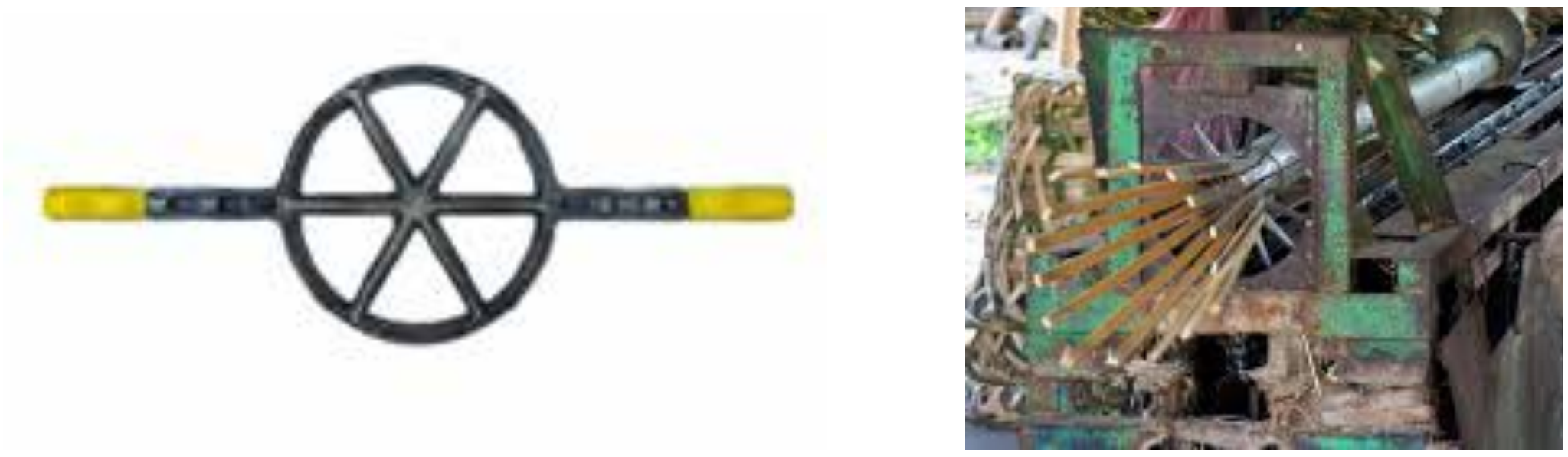

*Figure 8 Knives and machines for splitting bamboo [16].* ***a.*** *hand-operated tool moved to split a fixed culm.* ***b.*** *Bamboo splitting machine with fixed knife and moving culm.*

## 2.4 Cutting regular *m*-polygons with one knife

Having established the equivalence of cutting *GML* and polygons with a moving or fixed knife, now consider the regular *m*-polygon with rotational symmetry $m \in \mathbb{N}$, i.e. a regular polygon with *m* vertices and *m* sides connecting the vertices. There are three ways of cutting a regular *m*-polygon. A cut, which divides the structure into two distinct zones, can be made from:

1. vertex to vertex *VV*, with notation $VV_{i,j}$, a cut from vertex *i* to vertex *j*
2. vertex to side *VS*, with notation $VS_{i,j}$, a cut from vertex *i* to side *j*
3. side to side *SS*, with notation $SS_{i,j}$, a cut from side *i* to side *j*

In the case of *VV* the first *V* is labelled as $V_1$ so that a cut from the first vertex to the third vertex, in clockwise direction is called $VV_{1,3}$, which is for example the diagonal of a square (Figure 9**a**)[1]. The first index *i* of $VV_{i,j}$, $VS_{i,j}$ and, $SS_{i,j}$ refers to the first letter, and the second index *j* refers to the second letter (*V* or *S*). In the case of vertices and sides, the labelling both starts at 1, so at $V_1$ and $S_1$. Then $VS_{1,2}$ defines a cut from $V_1$ to $S_2$, which lies inbetween $V_2$ and $V_3$ (Figure 9**b**). If a cut passes through the centre of the polygon, it is denoted as subscript C (Figure 9**a**,**c**), although this notation may be dropped in obvious cases as in Figure 9**a**. In Annex 1 all possibilities are shown for $m = 6, \ldots, 10$.

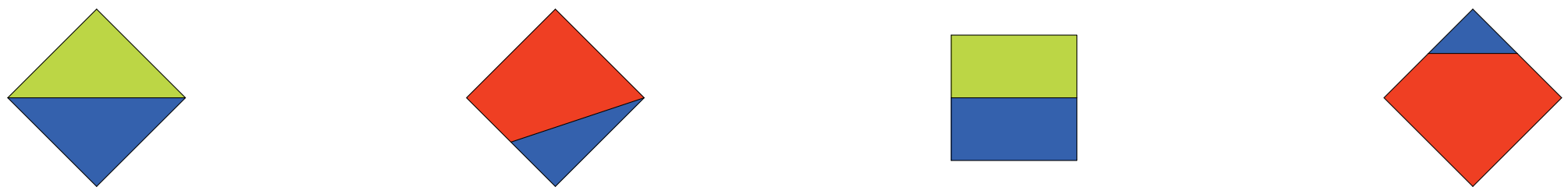

Figure 9. Ways of cutting of a square, left to right: **a**: $VV_{1,3,C}$, **b**: $VS_{1,2}$, **c**: $SS_{1,3,C}$ , **d**: $SS_{1,2}$

*Notes*

1. A cut is similar to drawing a straight line connecting vertices to vertices or sides, and sides to sides. Straight vertex-to-vertex lines are also known as *diagonals*.
2. A single cut divides the regular *m*-polygon into two distinct parts, i.e. the knife is part of a line, which has exactly two points in common with the boundary of the polygon. Such knife is called a *chordal* knife, analogous to the chord cutting a circle, giving rise to definition of sines and cosine. In Figure 8**a** the bamboo knife can be considered as three chordal knives through the centre.
3. Cuts, which are symmetrical with respect to clockwise or counter clockwise rotations are counted as one. If clockwise rotation is indicated as + and counter clockwise as -, then e.g. $VS_{1,-2} = VS_{1,2}$
4. $VV_{i,i+1}$ or $VV_{i,i-1}$ cuts are excluded as they coincide with a side and do no divide the polygon into two distinct parts. When cutting convex polygons with curved sides [4], or circles with equally spaced points $VV_{i,i+1}$ or $VV_{i,i-1}$ are possible.

The number of ways the *m*-polygon can be cut with one knife is equivalent to drawing one line. The total number of cuts is $3k - 2$ for even numbers $m = 2k$, and $3k - 1$ for odd numbers $m = 2k + 1$ (Table 1). For increasing *m* this leads to the following:

- From *m* odd to $m + 1$ (even) the number of possible cuts increases with 2, namely $VV_{1,(m-2)}$ and $SS_{1,(m-2)}$. For example, in $m = 6$ the new possibilities are $VV_{1,4}$ and $SS_{1,4}$, compared to $m = 5$.

[1] Labelling is from 1 to *m* for both vertices and sides. In earlier publications [1-6] numbering starts from 0 in the case of vertices.

- From $m$ even to $m+1$ (odd) the number of possible cuts increases with 1, namely $VS_{1,(m-2)}$. For example, for $m = 5$ the new possible cut is $VS_{1,3}$ compared to $m = 4$.

*Table 1: Possible cuts for m-polygons $m = 3, 4, \dots, 7$ and the general case*

| ***m*** | ***VV* cuts** | | ***VS* cuts** | | | ***SS* cuts** | | | ***Possible cuts*** |
|---|---|---|---|---|---|---|---|---|---|
| ***3*** | | | $VS_{1,2}$ | | | $SS_{1,2}$ | | | 2 |
| ***4*** | $VV_{1,3}$ | | $VS_{1,2}$ | | | $SS_{1,2}$ | $SS_{1,3}$ | | 4 |
| ***5*** | $VV_{1,3}$ | | $VS_{1,2}$ | $VS_{1,3}$ | | $SS_{1,2}$ | $SS_{1,3}$ | | 5 |
| ***6*** | $VV_{1,3}$ | $VV_{1,4}$ | $VS_{1,2}$ | $VS_{1,3}$ | | $SS_{1,2}$ | $SS_{1,3}$ | $SS_{1,4}$ | 7 |
| ***7*** | $VV_{1,3}$ | $VV_{1,4}$ | $VS_{1,2}$ | $VS_{1,3}$ | $VS_{1,4}$ | $SS_{1,2}$ | $SS_{1,3}$ | $SS_{1,4}$ | 8 |
| ***Even*** | $k-1$ | | | $k-1$ | | | $k$ | | $3k-2$ |
| ***Odd*** | $k-1$ | | | $k$ | | | $k$ | | $3k-1$ |

This is equivalent to an increase with 3 possibilities going from $m$ to *m+2. From m odd to m+2 (odd), or from m even to m+2 (even), the number of possible new cuts increases with 3 in both cases.* The number of possible cuts in Table 1, right column, is the sequence 2, 4, 5, 7, 8, 10, 11, 13… and this sequence (from odd to even plus 2 and from even to odd plus 1) will be continued for any increasing value of *m*. In the OEIS database of integer sequences, this sequence of numbers is: *numbers not divisible by 3.*

## 2.5 Divisors and different ways of cutting

The cuts described above (Figure 9 and Table 1) are made with precisely one knife corresponding to the smallest divisor of *m*, namely 1. The number of cuts with 1 knife is equivalent to divisor $d_1$. The results of cutting *GML* with regular polygons as cross-section or cutting regular polygons, aligns with the number of divisors of *m*.

The number of divisors *d* of $m \geq 2$. The number of divisors is 2 if and only if *m* is a prime number and the divisors are then $d_1$ and $d_m$. The smallest divisor $d_1$ is always 1 and the largest divisor $d_m$ is always equal to *m.* For all values of *m*, the divisors are designated $d_1, d_2, \dots, d_m$.

Cutting with *m* knives is equivalent to the largest divisor $d_m$. Using $d_1$ or one knife the cut started either in $V_1$ for *VV* and *VS* cuts, and from $S_1$ for *SS* cuts. For divisor $d_m$ or *m* knives, cuts start in all *m* vertices and on all *m* sides. The cuts are symmetrical to the case $d_1$ of one knife.

The cuts are based on repetition of $d_1$ cuts, by rotation, for example for hexagons and *VV* and *VS* cuts:

- *VV* cuts based on $VV_{1,3}$: $VV_{1,3}, VV_{2,4}, VV_{3,5}, VV_{4,6}, VV_{5,1}$ and $VV_{6,2}$
- *VV* cuts based on $VV_{1,4}$: $VV_{1,4}, VV_{2,5}, VV_{3,6}, VV_{4,1}, VV_{5,2}$ and $VV_{6,3}$
- *VS* cuts based on $VS_{1,2}$: $VS_{1,2}, VS_{2,3}, VS_{3,4}, VS_{4,5}, VS_{5,6}$ and $VS_{6,1}$
- *VS* cuts based on $VS_{1,3}$: $VS_{1,3}, VS_{2,4}, VS_{3,5}, VS_{4,6}, VS_{5,1}$ and $VS_{6,2}$

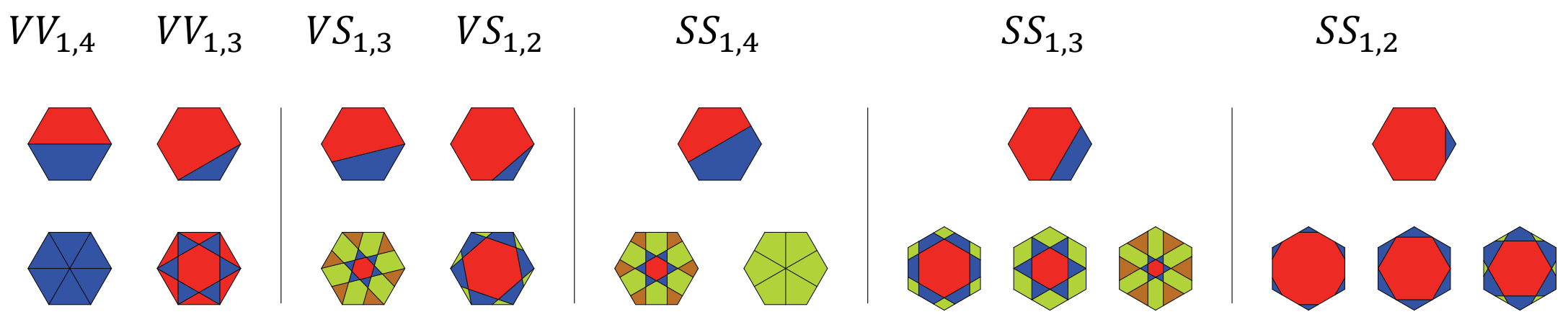

*Figure 10**a**: Cutting a hexagon with 1 ($d_1$upper row) or 6 knives ($d_6$ lower row)*

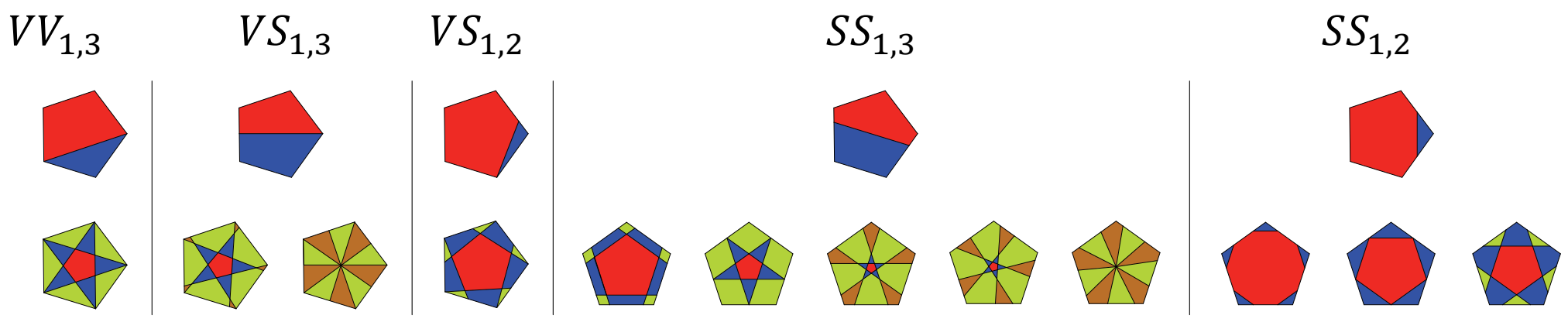

*Figure 10**b**: Cutting a pentagon with 1 ($d_1$upper row) or 5 knives ($d_5$ lower row)*

*VV* and *VS* cuts result in four different figures for both $d_1$ and $d_m = 6$ (Figure 10**a**), but with a one-to-one inheritance from $d_1$ to $d_6$ (or vice versa) i.e. the same total number of figures, namely 4. Additional possibilities will be created however, when *m* side-to-side cuts are considered for $d_m$.

- *SS* cuts based on $SS_{1,4}$: if the cut is made from the middle of $S_1$ to the middle of $S_4$, the cuts cross the centre of the hexagon (Figure 10**a**, column $SS_{1,4}$). In this particular case *m* is even, and the line contains the centre of symmetry.
- For *SS* cuts based on $SS_{1,2}$ or $SS_{1,3}$ the result depends on whether the cuts are made from
  - The middle of $S_1$ to the middle of $S_2$ or to the middle of $S_3$
  - From side to side such that the cutting line is shorter or longer than the cut from middle to middle.

This is shown in Figure 10**a**, columns $SS_{1,2}$ and $SS_{1,3}$. In total, 5 additional different ways of cutting result for $d_m$ compared to $d_1$. For the hexagon, for $d_1$ the total is 7, for $d_m = 6$, it is 12.

The results for cutting a pentagon with 1 knife ($d_1 = 1$) and 5 knives ($d_m = 5$) are shown in Figure 10**b**. For $d_1$ there are 5 possibilities. The additional $VS_{1,3}$ for $d_m = 5$ compared to a square, is generated when the cut is made through the centre of the pentagon. For $d_m = 5$ the total number is 12, as for $d_m = 6$ in a hexagon. For $SS_{1,3}$ there are 4 additional ways of cutting, compared to $d_1$ that fall into two categories.

- First, as in the case of the hexagon, the cut $SS_{1,3}$ can be made parallel to the side of the pentagon, again, from middle of $S_1$ to exact middle of $S_3$, or shorter, or longer, giving three extra possibilities as for the hexagon. The same for $SS_{1,2}$, giving a total of extra 4 possibilities.

- Second, the cut $SS_{1,2}$ gives rise to an additional 2 possibilities, depending on the cut made. From left of middle of $S_1$ to right of middle of $S_3$ or from right of middle of $S_1$ to left of middle of $S_3$

Going from a pentagon to a heptagon (Figure 11), the total number of ways of cutting for $d_m = 7$ is 17. Two extra vertices and two extra sides generate additionally $VS_{1,4}$ and $SS_{1,4}$. The former generates one extra possibility of cutting (Figure 11 column $VS_{1,4}$) like $VS_{1,3}$ in the pentagon, while the latter generates 5 extra possibilities (Figure 11 column $SS_{1,4}$), like $SS_{1,3}$ in the pentagon, compared to $d_1$. $SS_{1,3}$ and $SS_{1,2}$ generate three possibilities each, again dependent on where the cut is made (Figure 11 columns $SS_{1,3}$ and $SS_{1,2}$). This pattern will always continue going from $m$ to $m+2$ for $m$ odd. A similar reasoning leads to same conclusion for $m$ even.

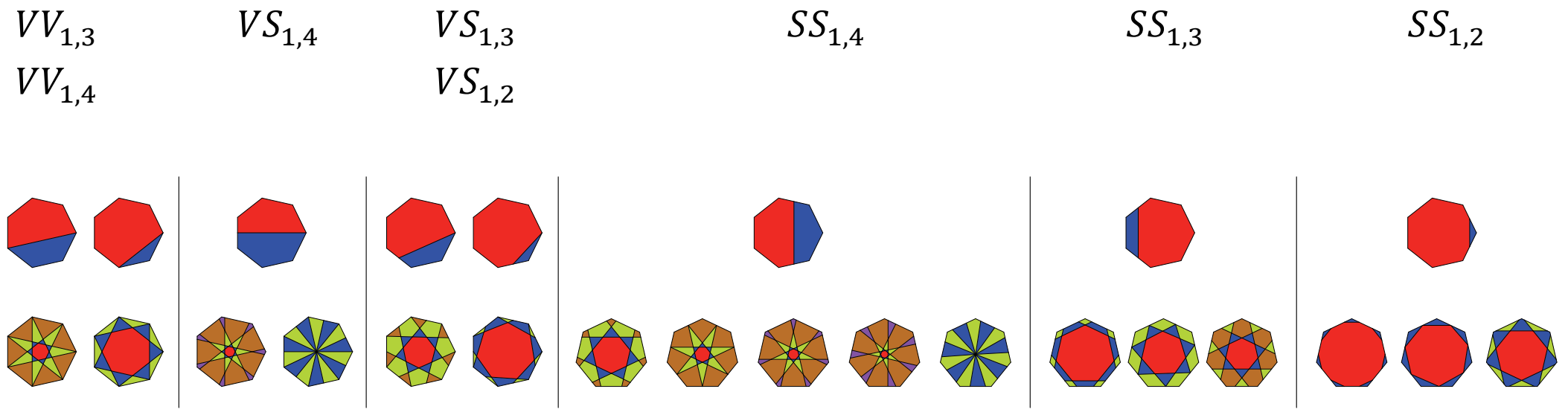

*Figure 11: Cutting a heptagon with 1 or 7 knives*

For $d_m$ the number of ways of cutting increases according the following rule: When $m$ increases from an even number $m = 2k$ to the next odd number $m = 2k + 1$, the number of possible cuttings increases by 5. When $m$ increases to the next even number $m = 2k + 2$ the number of possibilities does not change.

For $m = 3, 4, 5, 6, 7, \ldots.$ the number of possible cuts $d_m$ is the monotonic sequence $7, 7, 12, 12, 17, 17, 22, 22, \ldots$ or taken two by to $14, 24, 34, 44\ldots$ In general, for couples $m = 2k - 1$ and $m + 1 = 2k$ this gives $(k - 1) \cdot 10 + 4$.

## 2.6 Geometry versus topology

The main difference between $d_m$ and $d_1$ is that in $d_1$ certain cuts are not differentiated. For example in $SS_{1,4}$ in Figure 11, the line connecting sides 1 and 4, can lie anywhere on both sides. For $d_m$ however it is important where exactly sides 1 and 4 are cut. In Figure 9 column $SS_{1,4}$ the cut can go either through or not through the centre of the regular hexagon. When it does not go through the centre, various objects are generated, indicated by different colours.

This is the difference between the geometrical and topological identification of the solution, indicated in *Propositio*. In $SS_{1,4}$ in Figure 11, any cut with one knife (upper row) will generate two pieces, one with 5 vertices, one with 6 vertices (the topological condition). If the cuts with one knife need to generate shapes, which are one to one congruent (the geometrical condition) to the $d_m$ case, also the upper row will have 5 different variants.

In Table 2, the number of ways of cutting a regular polygons for smallest and largest divisor is given, and the total ways of cutting in the topological sense, the sum $d_1 + d_m$ .

- From an even number $m_{2n}$ to the next even number $m_{2n+2}$ : $d_1 + d_m$ = +8
- From an odd number $m_{2n+1}$ to the next odd number $m_{2n+3}$ : $d_1 + d_m$ = +8

The general rule is: The number of possibilities of cutting a regular $m$-polygon with one or *(m)* knives ($d_1$+$d_m$) is the number of possibilities of cutting a *(m-2)*-polygon increased by 8**.** As a consequence of the number 8 increase, the sequence of numbers 9,11, 17, 19, 25, 27, 33, 35, ..... is (1,3) *modulo* 8, i.e. the numbers of $d_1 + d_m$ are divisible by 8 with rest 1 or 3. In Table 2 the total number of possible cuttings for smallest and largest divisor of an $m$-polygon is given up to *m = 10.* The value of $d_1 + d_m$ for odd $m$ is 1 *modulo* 8 and the value of $d_1 + d_m$ for even $m$ is 3 *modulo* 8.

*Table 2: Number of cuts for smallest and largest divisors for m-polygons*

| *m* | $d_1$ | $d_m$ | $d_1$+$d_m$ |
|---|---|---|---|
| **3** | 2 | 7 | *9* |
| **4** | 4 | 7 | *11* |
| **5** | 5 | 12 | *17* |
| **6** | 7 | 12 | *19* |
| **7** | 8 | 17 | *25* |
| **8** | 10 | 17 | *27* |
| **9** | 11 | 22 | *33* |
| **10** | 13 | 22 | *35* |
| **....** | .... | .... | .... |
| $\boldsymbol{m = 2k+1}$ | $3k-1$ | $5k+2$ | $8k+1 = 4m-3$ |
| $\boldsymbol{m = 2k}$ | $3k-2$ | $5k-3$ | $8k-5 = 4m-5$ |

The total number of possible cuttings for smallest and largest divisor of an $m$-polygon is always an odd number. It is *4m-3* (or 1 mod 8) for odd values of $m$ and *4m-5* (or 3 mod 8) for even values of *m*. The total number of ways of cutting an $m$-polygon for m = prime is ($d_1$+$d_m$) = *4m-3* or 1 mod 8.

So we have the following:

**Lemma 1:** in the *topological* sense, the cutting with one knife ($d_1$) gives the minimum number of cuts, while cutting with m knives ($d_m$) gives the maximum number of cutting variants. In the *geometrical* sense, the minimum number of cutting variants (with a $d_1$-knife) equals the maximum number of cuts. ($d_m$-knife).

## 2.7 GML versus polygon cutting with 1 or *m* knives

The original study is a moving knife cutting 3D-Generalized Mobius Listing bodies (Figure 6) and continuing until the knife returns to the initial position (i.e. a complete cut). In Table 3 the cutting using a largest and smallest number of knives, related to smallest and largest divisors, is presented for *GML* cutting with moving knife. The condition $n = km$ corresponds to full rotation and $\boldsymbol{d_m = 1}$, and the condition $\boldsymbol{n = km + j\ (gcd(m, j) = 1}$ to $d_m = m$.

*Table 3: GML cutting with moving knives*

**$\boldsymbol{n = km}$ (full rotation)** — **Smallest divisor $\boldsymbol{d_m = 1}$**

| m-odd | All different variants of cutting "SS" (side to side) | All different variants of cutting "VS" (vertex to side) | All different variants of cutting "VV" (vertex to vertex) |
|---|---|---|---|
| **3** | 1 | 1 | 0 |
| **5** | 2 | 2 | 1 |
| **7** | 3 | 3 | 2 |
| **9** | 4 | 4 | 3 |
| **11** | 5 | 5 | 4 |
| | | | |
| 2k+1 | **k** | **k** | **k-1** |

**$\boldsymbol{n = km + j\ (gcd(m, j) = 1}$ - m-knife** — **Largest divisor $\boldsymbol{d_m = m}$**

| m-odd | All different variants of cutting "SS" (side to side) | All different variants of cutting "VS" (vertex to side) | | All different variants of cutting "VV" (vertex to vertex) |
|---|---|---|---|---|
| **3** | 5 | 1 | 1 | 0 |
| **5** | 3 + 5 | 2 | 1 | 1 |
| **7** | 3 + 3 + 5 | 3 | 1 | 2 |
| **9** | 3 + 3 + 3 + 5 | 4 | 1 | 3 |
| **11** | 3 + 3 + 3 + 3 + 5 | 5 | 1 | 4 |
| | | | | |
| 2k+1 | **3·(k-1) + 5** | **k** | **1** | **k-1** |

**$\boldsymbol{n = km}$ (full rotation)** — **Smallest divisor $\boldsymbol{d_m = 1}$**

| m-even | All different variants of cutting "SS" (side to side) | All different variants of cutting "VS" (vertex to side) | All different variants of cutting "VV" (vertex to vertex) |
|---|---|---|---|
| **2** | "1" | | 0 |
| **4** | 2 | 1 | 1 |
| **6** | 3 | 2 | 2 |
| **8** | 4 | 3 | 3 |
| **10** | 5 | 4 | 4 |
| | | | |
| 2k | **k** | **k-1** | **k-1** |

**$\boldsymbol{n = km + j\ (gcd(m, j) = 1}$ (m-knife)** — **Largest divisor $\boldsymbol{d_m = m}$**

| m-even | All different variants of cutting "SS" (side to side) | | All different variants of cutting "VS" (vertex to side) | All different variants of cutting "VV" (vertex to vertex) | |
|---|---|---|---|---|---|
| **2** | 1 | 1 | | 0 | 0 |
| **4** | 3 + 1 | 1 | 1 | 0 | 1 |
| **6** | 3 + 3 + 1 | 1 | 2 | 1 | 1 |
| **8** | 3 + 3 + 3 + 1 | 1 | 3 | 2 | 1 |
| **10** | 3 + 3 + 3 + 3 + 1 | 1 | 4 | 3 | 1 |
| | | | | | |
| 2k | **3·(k-1) + 1** | **1** | **k-1** | **k-2** | **1** |

Table 4 gives a summary. Comparison with Table 2 Totals shows that the correspondence of *GML* cutting and regular polygons is indeed complete.

*Table 4: Comparative table for GML and polygon cutting*

| Odd | GML in 3D | $d_i$ −knives | SS-cuts | VS-cuts | VV-cuts | Total |
|---|---|---|---|---|---|---|
| ***n=2k+1*** | *n=km* | *1* | *k* | *k* | *k-1* | ***3k-1*** |
| | *n=km+j(gcd(m,j)=1* | *m* | *3(k-1)+5* | *k+1* | *k-1* | ***5k+2*** |
| ***n=2k*** | *n=km* | *1* | *k* | *k-1* | *k-1* | ***3k-2*** |
| | *n=km+j(gcd(m,j)=1* | *m* | *3(k-1)+2* | *k-1* | *k-1* | ***5k-3*** |

## 2.8 Defining knives for all divisors

For $d_1$ cuts are made starting from $V_1$ or $S_1$. For $d_m$ this operation was repeated every $2\pi/_m$. For divisors other than $d_1$ or $d_m$, the operation has to be repeated less than every $2\pi/_m$. In Table 5 the rotations which have to be performed is indicated for $m$-polygons with more than two divisors for $m$ = 4, 6, 8, 9 and 10. The denomination $d_3 = 4$ indicates that for $m = 8$ having four divisors 1, 2, 4 and 8, the third divisor is 4. The denomination $d_m$ always refers to the highest divisor, in this case equal to $d_m = 8$ for the octagon. The rotations are $k \cdot 2\pi/_m$. The value of $k$ depends on the other divisors. For example for $m = 10$, $k = 5$ for $d_2 = 2$, ; $k = 2$ for $d_3 = 5$, $k = 1$ for $d_m = 10$ and $k = 10$ for $d_1$ . The results are shown in Figure 12 for $m = 4$.

*Table 5: Traces of d-knives*

| | $d_1$ | $d_2$ | $d_3$ | $d_m$ |
|---|---|---|---|---|
| ***m=4*** | $d_1 = 1$<br>$4 \cdot 2\pi/_4$ or 360° | $d_2 = 2$<br>$2 \cdot 2\pi/_4$ or 180° | | $d_m = 4$<br>$2\pi/_4$ or 90° |
| ***m=6*** | $d_1 = 1$<br>$6 \cdot 2\pi/_6$ or 360° | $d_2 = 2$<br>$3 \cdot 2\pi/_6$ or 180° | $d_3 = 3$<br>$2 \cdot 2\pi/_6$ or 120° | $d_m = 6$<br>$1 \cdot 2\pi/_6$ or 60° |
| ***m=8*** | $d_1 = 1$<br>$8 \cdot 2\pi/_8$ or 360° | $d_2 = 2$<br>$4 \cdot 2\pi/_8$ or 180° | $d_3 = 4$<br>$2 \cdot 2\pi/_8$ or 90° | $d_m = 8$<br>$1 \cdot 2\pi/_8$ or 45° |
| ***m=9*** | $d_1 = 1$<br>$9 \cdot 2\pi/_9$ or 360° | $d_2 = 3$<br>$3 \cdot 2\pi/_9$ or 270° | | $d_m = 9$<br>$1 \cdot 2\pi/_9$ or 40° |
| ***m=10*** | $d_1 = 1$<br>$10 \cdot 2\pi/_{10}$ or 360° | $d_2 = 2$<br>$5 \cdot 2\pi/_{10}$ or 180° | $d_3 = 5$<br>$2 \cdot 2\pi/_{10}$ or 72° | $d_m = 10$<br>$1 \cdot 2\pi/_{10}$ or 36° |

This then leads to the following analytic definition of $d_{i=1,..,m}$-knife: the $d_i$-knife is a construction, with $m$ straight lines [8]:

$$\sin\left(\alpha + \frac{2\pi}{m} i\right) x_i + \cos\left(\alpha + \frac{2\pi}{m} i\right) y_i + \delta = 0\,, \quad i = 0, 1, .., m-1; \quad -\frac{\pi}{m} \le \alpha \le \frac{\pi}{m}, \quad (12)$$

$$y_i = \tan\left(\alpha + \frac{2\pi}{m} i\right) x_i + \delta \left[\cos\left(\alpha + \frac{2\pi}{m} i\right)\right]^{-1}, \quad i = 0, 1, .., m-1; \quad -\frac{\pi}{m} \le \alpha \le \frac{\pi}{m}, \quad (13)$$

The parameter $\alpha$ is a *rotation* parameter and $\delta$ is a *dilation* or zooming parameter. With these parameters all possible ways of cutting can be described. In Figure 12 middle row one can observe how the red square is dilated and rotated. Dilation can shrink the square to a point as in the $VV_{1,3}$ and $SS_{1,3}$ cut through the centre. This will be used in *Demonstratio* Section 4.3 since this is equivalent to projection. In evaluating other divisors, the value of $i$ should be amended, according to Table 5. For example, in the octagon $d_3 = 4$ the cut is repeated every 90°, or $i = 2$: $^{2 \cdot 2\pi}/_{8}$ or 90°.

For $d_2 = 2$ (equivalent to $m$ = even) the number of possible ways of cutting in the topological sense is the number of possibilities for $d_1$ +1**.** In evaluating the number of divisors in relation to *GML* or regular polygon cuts, it is found that for *VV* and *VS* all divisors are equal for even *m*. For odd m the same goes for *VV* cuts, but for *VS* cuts one has to take into account whether or not the cut goes through the centre, at least in the *topological* sense. The number of ways of cutting *VV* or *VS* is in *geometrical* sense inherited one-to-one from *VS* or *VV* cutting for smallest and largest divisors, basically this number is constant for all *VS, VV* cuts for all divisors.

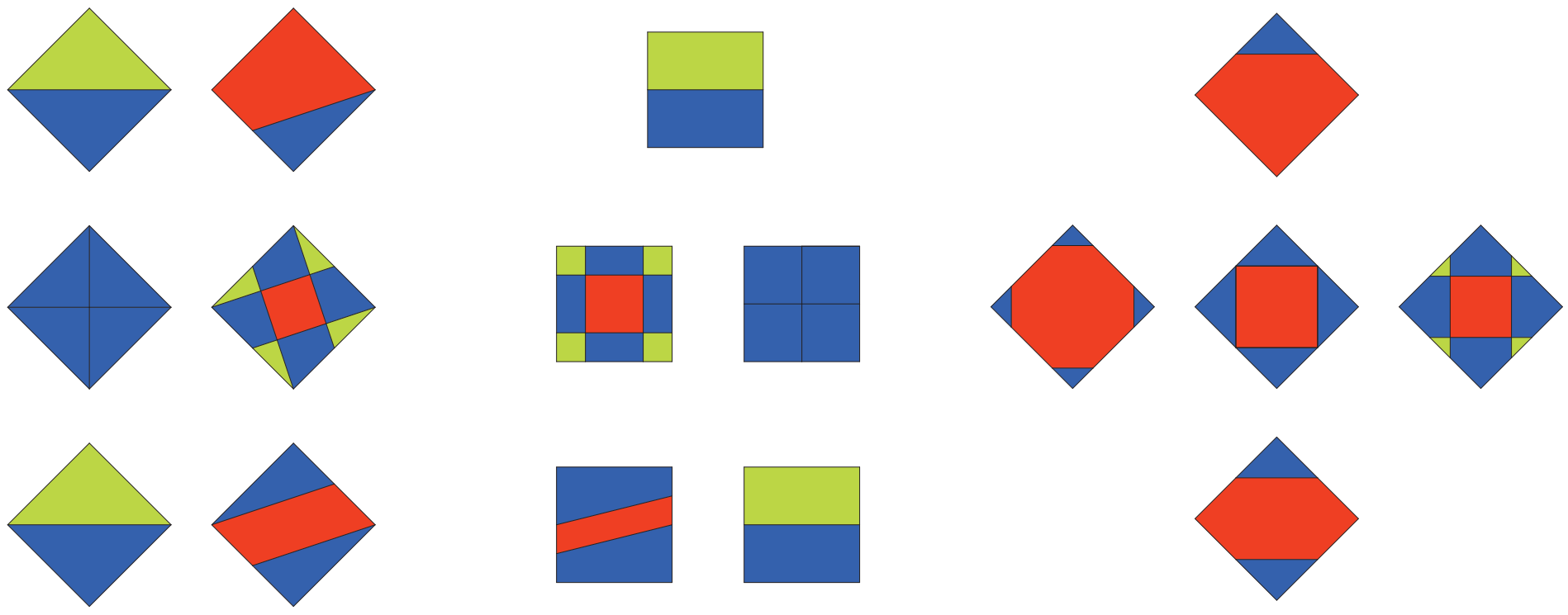

*Figure 12: Ways of cutting for a square for all divisors* $d_1$ *(upper row)* , $d_2 = 2$ *(lower row) and* $d_m = 4$ *(central row).*

The challenge will be determining the regularity for *SS* cuts. In enumerating the total number of ways of cutting of *m*-polygons, we need to take into account the fact that we do not know a priori how many divisors *m* has ( https://oeis.org/A000203).

The divisors of a number *m* may be numbered from minimum to maximum (min and max) leading to the following sequence:

$$d_{min}, d_{min+1}, \dots\dots, d_{max-2}, d_{max-1}, d_{max}$$

where:

$d_{min}$ is the smallest divisor (always equal to 1) and,
$d_{max}$ is the largest divisor (always equal to $m$)

A number has either even or odd number of divisors e.g. $(1, 2, 3, 4, 6, 12)$ for 12 and $(1, 3, 9)$ for 9. If the number of divisors is odd, the sequence has a central value (e.g. 3 in sequence for 9), but this may be used twice, so the sequence is $(1, 3, 3, 9)$ as in $(1, 2, 4, 4, 8, 16)$, so that all numbers have an even number of divisors.

The divisors can be grouped pairwise as:

$$(d_{min}, d_{max}), (d_{min+1}, d_{max-1}), (d_{min+2}, d_{max-2}), \ldots..$$

- For $m = 9$ the pairs are $(1, 9)$; $(3, 3)$.
- For $m = 12$ the pairs are $(1, 12)$; $(2, 6)$; $(3, 4)$
- For $m = 16$ the pairs are $(1, 16)$; $(2, 8)$; $(4, 4)$

The products of these pairs are always equal to *m*. When relating this to rotational symmetries $C$ for $m = 12$, the pairs are $(C_1, C_{12})$; $(C_2, C_6)$; $(C_3, C_4)$ and these relate to the spacing of equidistant points on the circle, and the angles between them. For example, $C_2, C_6$ points are spaced 180° and 60° apart, respectively, corresponding to the rotation parameter for divisor 2 and 6.

## 2.9 Back to Flatland

The cutting of regular polygons is an equivalent method for demonstrating the general cutting of $GML_m^n$, with a one-to-one correspondence of the cutting of regular polygons and the cross sections of $GML_m^n$, bodies. As a consequence the problem may be reversed as follows:

> *Given the possibility of cutting planar regular m-polygons (as discussed in Section 2) obtained with a number of knives equal to* $d_1$, $d_2$ *….,* $d_m$*, can the same result be obtained using only one knife?*
>
> *The answer is affirmative, if the operation is carried out in 3D, via GML bodies with a specific number of twists.*

But what does an inhabitant of Flatland observe when operations with 3D $GML_m^n$ bodies or surfaces are carried out, in other words, when $GML_m^n$ bodies rotate through Flatland? Or, when $GML_m^n$ bodies do not rotate but the knives are moving? Or when we position the fixed knife in Flatland with rotating $GML_m^n$ bodies?

An inhabitant of Flatland may observe at some instance in time what looks like a pentagonal piece of paper lying on his/her desk, and may put with one knife across the paper (*VV, VS* or *SS* direction). Next morning the Flatlander may wake up to see how the whole situation is still the same. However, it might also be that the "pentagonal piece of paper", which is actually a planar section of a $GML_m^n$, rotating through Flatland, is sliced into different zones (Figures 7**a** and 10**b**); how many zones depends on the number of twists of the $GML_m^n$ and the number of rotations through Flatland. The Flatlander may also use more knives. He or she does not have to do anything at all; one dimension higher a rotation of the $GML_m^n$ may have occurred.

Performing this with a square, somewhere else in Flatland a second, identical operation has happened. If both Flatlanders come into contact with each other, they will find that the square is somehow turned 90° (Figure 6**a**). In classic Euclidean geometry, one cannot distinguish between two squares rotated 90° or a multiple thereof. If the $GML_m^n$ bodies are twisted, experimental identification of differences may be performed on the ground. In one case the diagonal runs "North to South", in the other one "East to West".

Without knowing the *GML* structure and number of knives, a Flatlander is bamboozled by this spooky action at a distance. If the Flatlander tries to untangle the mystery, by repeating or experimentally using other knives and positions, this will result in fully predictable structures, and reproducible phenomena elsewhere in Flatland, for example diagonals in the two squares in the same or different directions. The diagonal also has direction, like a spin, and the one influences the other. This is reminiscent of quantum entanglement, now via *GML* bodies and cutting.

Figure 6**a** shows a *GML* touching Flatland as a plane. Only half of the *GML* is shown, and the second half is beneath Flatland. But when the plane is folded, the complete *GML* body or surface can become a generalized cylinder (Figure 13). This cylinder can be a complete *GML* or it can be part of the *GML*. In the latter case the *GML* can be closed either in 3D by a completed *GML*, not shown in the figure and not touching the Brane, or by direct connections on one side of plane or brane (think streamlines, an electronic circuit or graph). In addition, the width of the *GML* need not be constant, so it can be a generalised cone, a generalized hyperboloid, twisted beams or wormholes instead of a generalized cylinder. The brane in Figure 13**b** is a fold of the planar surface in Figure 6**a**, and the grey surface of Figure 6**a** or the grey surfaces of Figures 13**a** and **b** themselves can be folded topologically into a torus. The *GML* body then connects the surface of the torus like a fistula within the torus (Figure 13), or a surface or body of genus 2 results, whereby the upper half of the *GML* in Figure 6**a** is a handle.

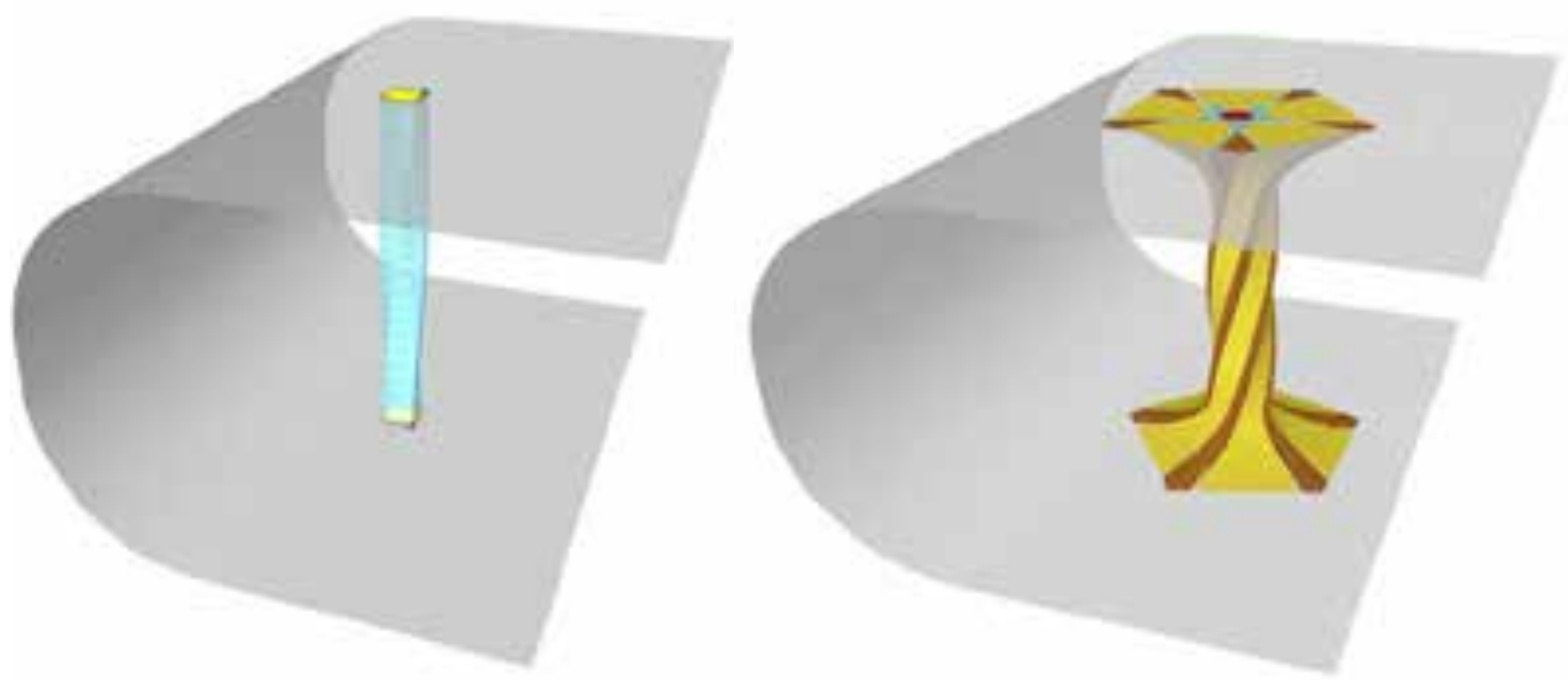

**a** **b**

Figure 13: Partial or complete *GML* bodies, with cross sections in Braneland.

Actually, the cutting via divisors or via knives can also be considered as a process of projection, whereby a certain shape is projected onto a basic shape. The rotation and zoom parameters define the projection. In Figure 13 the upper yellow zone is projected onto the lower part of the brane, along the *GML*. In Figure 12 middle row for example a red square is projected on to a base square of fixed size (blue). This projection can be performed with a square for $d_{max} = 4$, with a line for $d_{min} = 1$ or with a rectangle for $d_2 = 2$. One can observe from Figure 12 that this gives similar results as cutting with $d$ knives and rotations. Projection may occur via a classic cone, a pyramid, a prism or via any projective device, for example with varying width in the *GML*, even when the width becomes small as in a vortex, or goes to zero.

## 2.10 Chordal and Radial knives

In previous publications it was shown how in some cases the process of cutting leads to the Möbius phenomenon. Based on the current method of cutting, the Möbius phenomenon only occurs for $m$-polygons under two conditions: 1) $m$ is even, and 2) the cut is made through the centre. It does not occur in any case of $m$ odd, since the first and foremost condition is a geometrical one: the shapes have to be congruent after rotation to be able to form a single body! If $m$ is odd there is no case in which cutting leads to a single body for the pentagon.

Thus far only cuts were from perimeter to perimeter (including no, one or two vertices), cutting the planar polygons into two different parts. These knives are chordal knives, since it relates to a chord in the classical sense. However, one can consider also shorter knives. A knife starting in the centre and passing through one vertex or side is called a radial knife. The notion of $m$-knives, chordal or radial knives can be viewed in a broader perspective. In this sense chords and radii are intervals of a line and half-line respectively, falling wholly within the $m$-polygon or circle and having at least two points in common or more for self-intersecting shapes (Figure 14). It is noted that curves are defined in this case as regular Gielis polygons *sensu* Matsuura [12]. The pentagon is a $G_{m,p,\, p\cdot\left(\frac{5}{4}\right)^2} = G_{5,1,\left(\frac{5}{4}\right)^2}$; the pentagram is $G_{\frac{5}{2},1,\frac{25}{16}}$ and the circumscribing circle is $G_{5,1,2}$.

In a circle, the length of the radius is half the length of the longest chordal knife (Figure 14). The inspiration can be found in botany once more: when sawing wood, one can use radial or tangential cuts. The latter is a chordal cut. When this cut is from vertex to vertex, the chordal knife may also be called a diagonal knife.

A more refined definition is then:

- $d_{cc}$ chordal knife, through the centre $C$ (e.g. $VS_C$ cut in odd $m$ polygon)
- $d_{c\bar{c}}$ chordal knife, not through the centre
- $d_{rc}$ radial knife originating at the center
- $d_{r\bar{c}}$ radial knife not through the centre

The $d_{r\bar{c}}$ -knife is defined for completeness, but is a special case of a chordal knife with length shorter than the chordal one). A radial cut $d_{rc}$ is also a half-line or ray.

Interestingly our knives, as part of lines and half-lines, can be considered in a very classical way: When cutting a polygon, e.g. a pentagon (Figure 14), the circle circumscribing the $m$-polygon, has the vertices in common (as roots of unity). A chordal cut (the yellow line in Figure 14**a** is a $VV_{i,i+2}$ cut; the dotted lines are $VV_{i,i+1}$ cuts) then connects two points on the circle, and defines the associated trigonometric function *sine* (chord = 2.*sine*), which is maximal when the chordal cut is through the centre and minimal when cutting only one vertex. In this case the chord, when prolonged, is the tangent to a point on the circle. This chordal cut can be taken from any point to any other point on the circle, and the perpendicular can be drawn, which defines the *cosine*. Indeed, the cosine is defined by drawing the perpendicular to the chord through the centre (Fig 14**a.** red). This divides the chord into two equal pieces.

Prolonging the red perpendicular (the cosine) to both sides gives the diameter. The max value of the cosine in a unit circle with radius $R = 1$ is also the magnitude or length of the normal to the circle (radius of curvature). Hence tangent and curvature are defined in the same setting. Figure 14**b** shows the relation with other classical trigonometric functions. A radial knife has length 1, and is the sum of cosine $OC$ and versine[2] $CD$, or the sum of sine $OG$ and coversine $GH$.

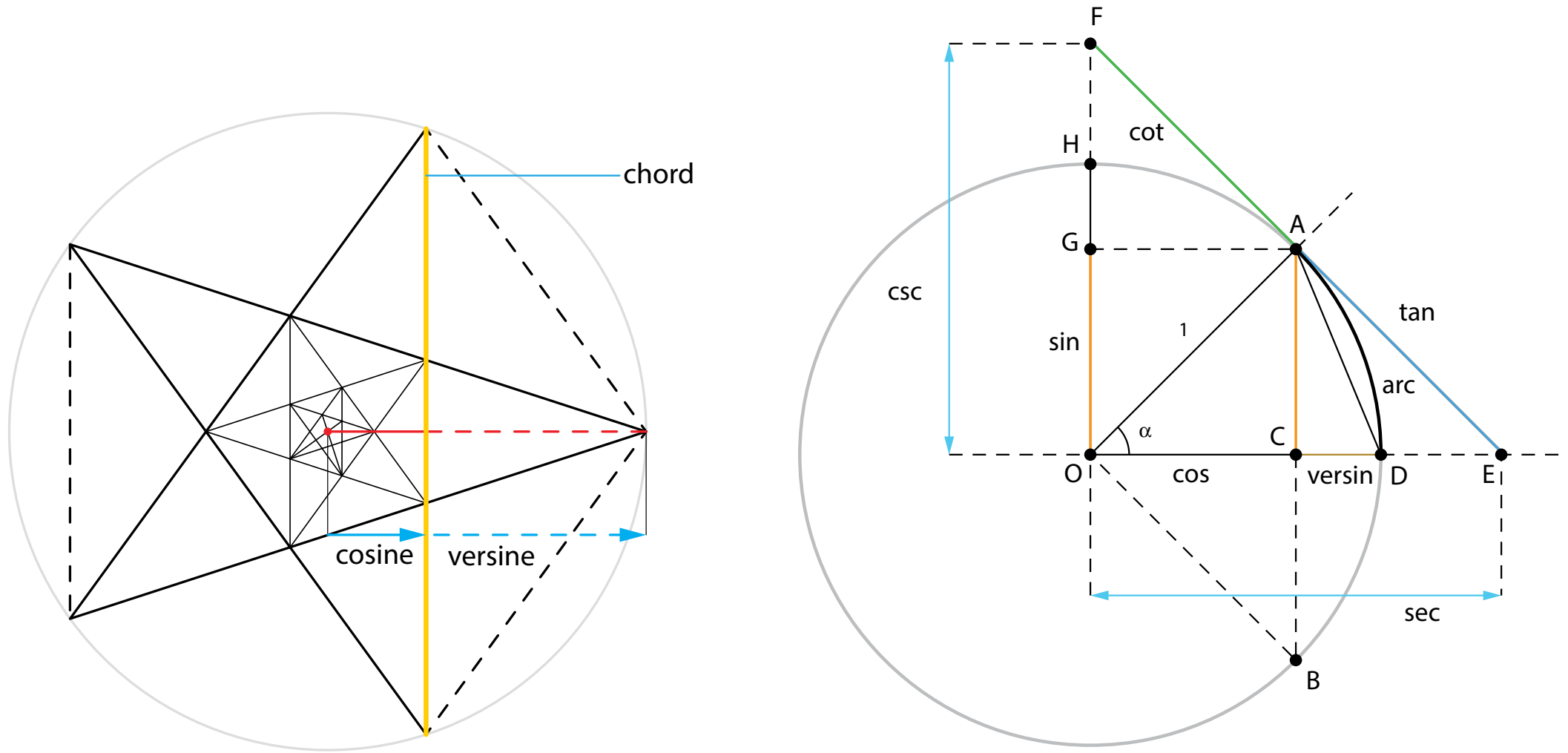

Figure 14: **a**. Radial (red, dotted+solid) and chordal (yellow) knives.
**b.** Classical trigonometric functions

Using chordal or radial knives the results of cutting will be the same for $GML$ surfaces and bodies, but the distinction will become very important when discussing the occurrence of the Möbius phenomenon. Chordal knives cut the boundary, but it is also possible to define knives defined by a Jordan curve completely contained inside the original boundary. In this way two zones are created. If the inner zone is then cut with chordal or radial knives, similar results are obtained. If in Figure 14 a knife is used corresponding to the thick pentagon inside, then the inner pentagon is separated from the rest of the original $GML_m^n$ body. When radial or chordal knives are used, the same way of cutting gives the smaller pentagram. This can be continued *ad infinitum*[3].

---

[2] The versine was used to measure curvedness in railroad tracks, by spanning a chord $AB$ and measuring the distance $CD$ from the chord to the track.

[3] This argument can be used to proof the incommensurability between side and diagonal.

# 3 Determinatio

In *Expositio* it was shown that the processes of cutting *GML* bodies or $m$-polygons are equivalent and related to the divisors and the number of knives. This equivalence will be used later for the demonstration and proof. In this *Determinatio* section, it is investigated how the process of cutting is related to the theory of knots and links, the distribution of points on the circle (roots of unity), the ways of cutting a polygon, graph theory, the theory of numbers and primes, combinatorial aspects etc ….. In this way many viewpoints are combined from various areas (geometry and algebra), with a very broad applicability for the general results of cutting of *GML* surfaces and bodies.

One crucial point of this solution with $d$-knives (radial or chordal), is that this solves a problem concerning 3D shapes and bodies using 2D geometry. In Riemannian geometry for example, one can determine the sectional curvature of an $n$-dimensional manifold, by determining the curvatures of planar 2D-cross sections in the tangent spaces [17]. In the same way, the current method would allow to solve problems in $n$-D *GML* bodies and surfaces, by studying 2D sections and solutions.

## 3.1 *GML* bodies and relation to knots and links after cutting

The use of planar geometry and cutting aims at providing a proof/demonstration for the cutting of *GML* bodies and surfaces. At first it makes abstraction of the intrinsic complexity of the cutting process, namely the number of bodies that are generated, how they are twisted and how they are interconnected. Indeed for 3D *GML* surfaces and bodies, the shapes can be very complex, and a full classification has been achieved for *m= 2, 3, 4, 5, 6.* The case 2 refers to ribbons, but can also be 2-angular figures [1-6].

We recall some results from [5]. In Figure 15**a** the results are shown for cutting a $GML_4^n$ body, which is twisted $4\omega$ times (e.g. a full rotation over 360°; compare Figure 3, but $\omega$ can be 0 for untwisted prism). For example the $SS_{1,2}$ (A) generates both a triangular and a pentagonal shape and depending on the initial value of $\omega$, the shapes are twisted 3 and 5 times respectively. The complete structure consists of the two shapes interlinked (Link-2). This situation corresponds to cutting with 1 knife or $d_1$ (note that the $SS_{1,3}$ and $SS_{1,3,C}$ cut give the same result so that total number is four).

The $SS_{1,3,C}$ cut of Figure 15**b** (BII) generates a single shape whereby the original shape is twisted 8 times ($8\omega + 8$ for $\omega = 0,1,2,..$). These situations correspond exactly to cutting squares with two knives or $d_2 = 2$ in section 2.7, with 5 possibilities, one more compared to one knife. In Figure 12 the same shapes are generated as cutting with four knives or $d_m = 4$. The essential aspect of cutting *GML* surfaces and bodies, and one of the original motivations is that the result of cutting leads to linked bodies with a certain link number (right columns). The surfaces or bodies may be knot like or display the Möbius phenomenon, with one shape and Link number 1, and this occurs in BII and D in both Figures 15**b** and 15**c** (for $m$ even).

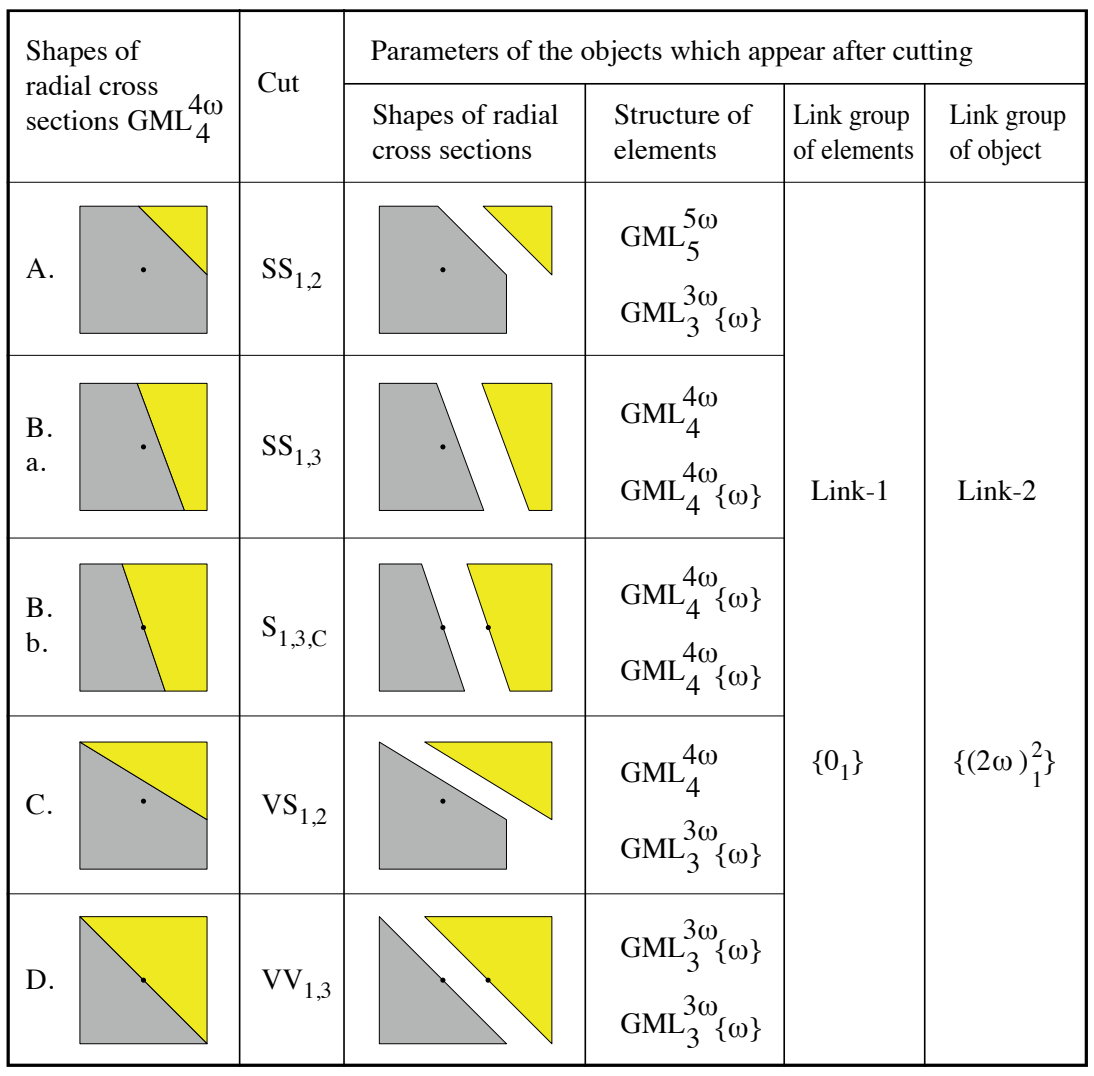

| Shapes of radial cross sections $GML_4^{4\omega}$ | Cut | Parameters of the objects which appear after cutting | | | |
|---|---|---|---|---|---|
| | | Shapes of radial cross sections | Structure of elements | Link group of elements | Link group of object |
| A. | $SS_{1,2}$ | | $GML_5^{5\omega}$ $GML_3^{3\omega}\{\omega\}$ | Link-1 $\{0_1\}$ | Link-2 $\{(2\omega)_1^2\}$ |
| B. a. | $SS_{1,3}$ | | $GML_4^{4\omega}$ $GML_4^{4\omega}\{\omega\}$ | | |
| B. b. | $S_{1,3,C}$ | | $GML_4^{4\omega}\{\omega\}$ $GML_4^{4\omega}\{\omega\}$ | | |
| C. | $VS_{1,2}$ | | $GML_4^{4\omega}$ $GML_3^{3\omega}\{\omega\}$ | | |
| D. | $VV_{1,3}$ | | $GML_3^{3\omega}\{\omega\}$ $GML_3^{3\omega}\{\omega\}$ | | |

***a***

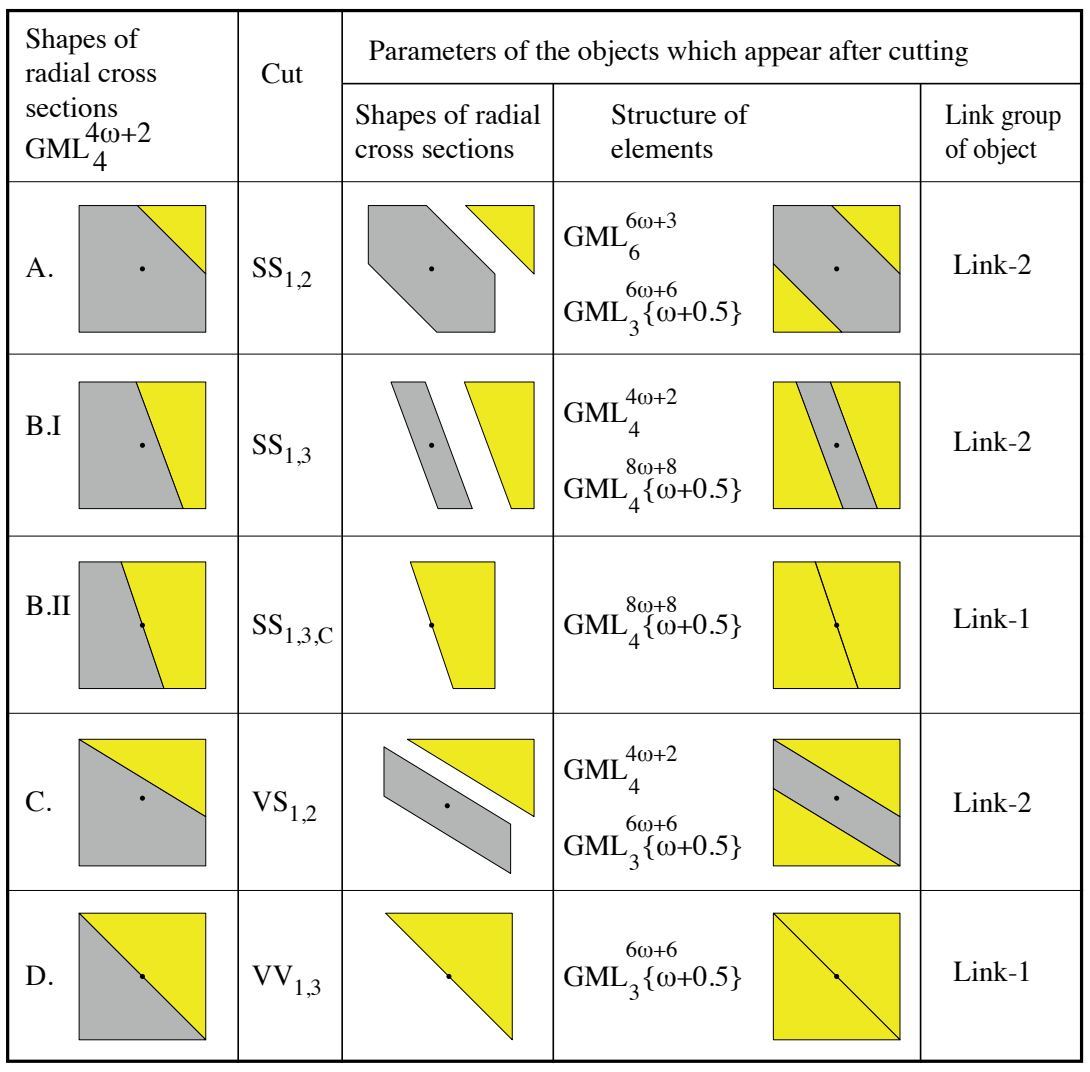

| Shapes of radial cross sections $GML_4^{4\omega+2}$ | Cut | Parameters of the objects which appear after cutting | | |
|---|---|---|---|---|
| | | Shapes of radial cross sections | Structure of elements | Link group of object |
| A. | $SS_{1,2}$ | | $GML_6^{6\omega+3}$ $GML_3^{6\omega+6}\{\omega+0.5\}$ | Link-2 |
| B.I | $SS_{1,3}$ | | $GML_4^{4\omega+2}$ $GML_4^{8\omega+8}\{\omega+0.5\}$ | Link-2 |
| B.II | $SS_{1,3,C}$ | | $GML_4^{8\omega+8}\{\omega+0.5\}$ | Link-1 |
| C. | $VS_{1,2}$ | | $GML_4^{4\omega+2}$ $GML_3^{6\omega+6}\{\omega+0.5\}$ | Link-2 |
| D. | $VV_{1,3}$ | | $GML_3^{6\omega+6}\{\omega+0.5\}$ | Link-1 |

***b***

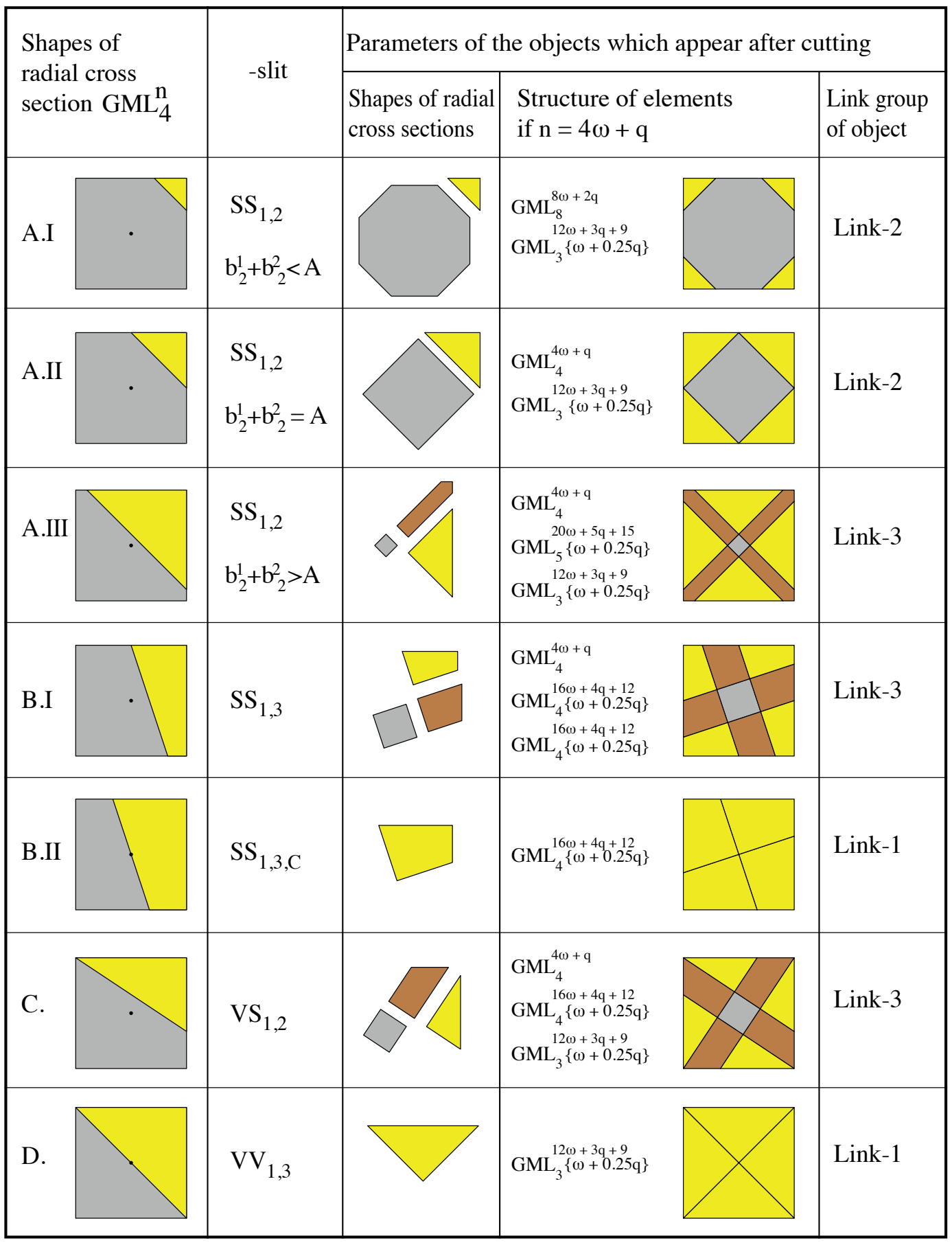

| Shapes of radial cross section $GML_4^n$ | -slit | Parameters of the objects which appear after cutting | | |
|---|---|---|---|---|
| | | Shapes of radial cross sections | Structure of elements if $n = 4\omega + q$ | Link group of object |
| A.I | $SS_{1,2}$ $b_2^1+b_2^2<A$ | | $GML_8^{8\omega+2q}$ $GML_3^{12\omega+3q+9}\{\omega+0.25q\}$ | Link-2 |
| A.II | $SS_{1,2}$ $b_2^1+b_2^2=A$ | | $GML_4^{4\omega+q}$ $GML_3^{12\omega+3q+9}\{\omega+0.25q\}$ | Link-2 |
| A.III | $SS_{1,2}$ $b_2^1+b_2^2>A$ | | $GML_4^{4\omega+q}$ $GML_5^{20\omega+5q+15}\{\omega+0.25q\}$ $GML_3^{12\omega+3q+9}\{\omega+0.25q\}$ | Link-3 |
| B.I | $SS_{1,3}$ | | $GML_4^{4\omega+q}$ $GML_4^{16\omega+4q+12}\{\omega+0.25q\}$ $GML_4^{16\omega+4q+12}\{\omega+0.25q\}$ | Link-3 |
| B.II | $SS_{1,3,C}$ | | $GML_4^{16\omega+4q+12}\{\omega+0.25q\}$ | Link-1 |
| C. | $VS_{1,2}$ | | $GML_4^{4\omega+q}$ $GML_4^{16\omega+4q+12}\{\omega+0.25q\}$ $GML_3^{12\omega+3q+9}\{\omega+0.25q\}$ | Link-3 |
| D. | $VV_{1,3}$ | | $GML_3^{12\omega+3q+9}\{\omega+0.25q\}$ | Link-1 |

***c***

*Figure 15. Cutting* **a.** $GML_4^{4\omega}$, **b.** $GML_4^{4\omega+2}$ $(\omega = 0,1,2,..)$, **c.** $GML_4^n$, with resulting **n**umber, shape and twisting of resulting bodies and their link number [5]. The independent objects and their shapes are given in different colours.

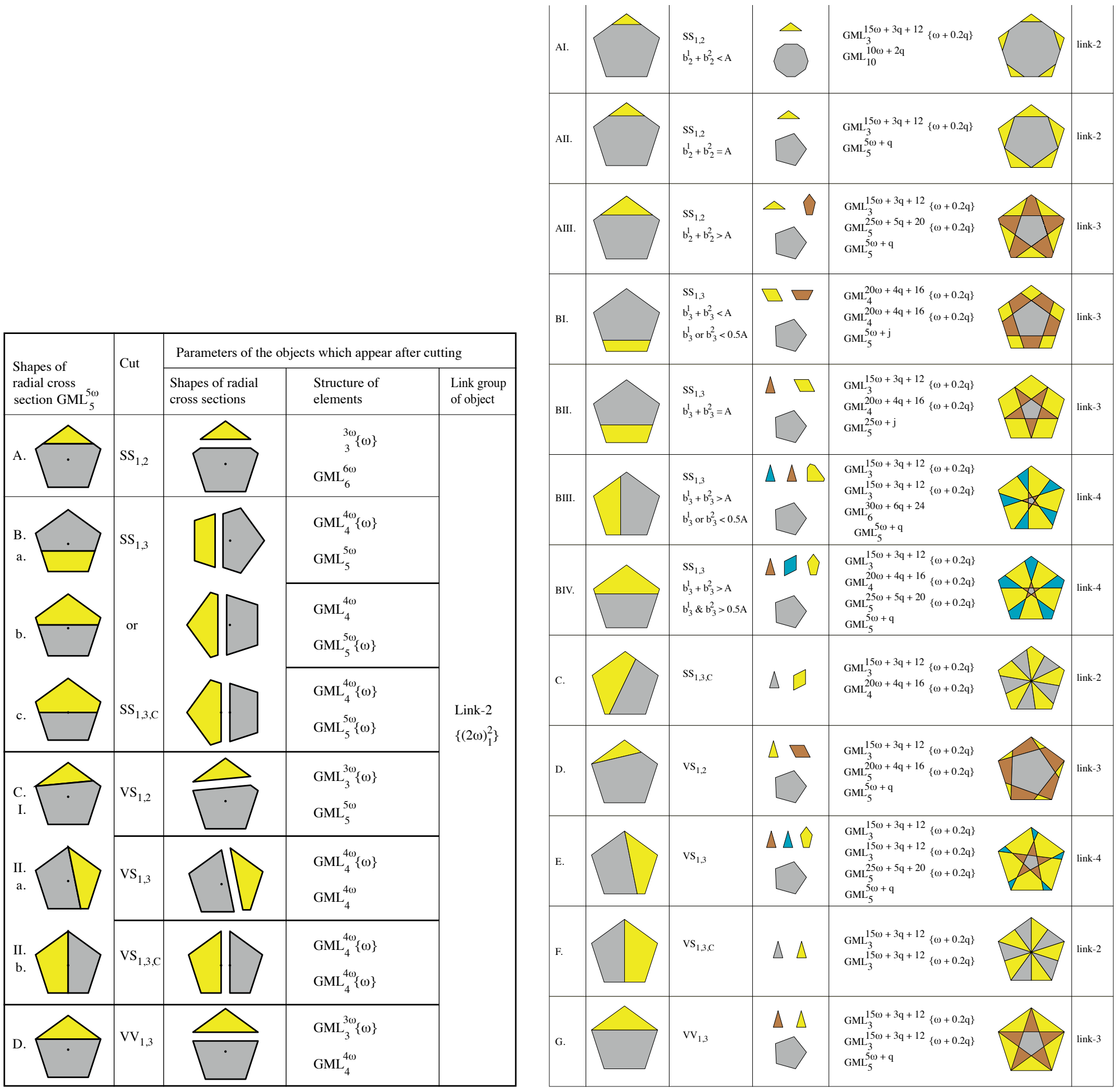

*Figure 16. Cutting* ***a.*** $GML_5^{5\omega}$, ***b.*** $GML_5^{n}$. Number, shape and twisting of resulting bodies and their link number [5]. Different objects are given in different colours.

In Figure 16**a**, in the left column 8 ways of cutting are used, instead of the 5 ways in Figure 10. These 8 ways however, can be reduced to 5 since in this way of cutting the centre of the pentagon is not important:

- All $SS_{1,3}$ cuts (*Ba,b,c*) yield exactly the same results, namely four angular bodies and five angular bodies with the same link-2 number.
- Both *C II* cuts (*CIIa* and *CIIb*), also yield the same result

Hence the number of ways reduces to 5 if we consider topological solution, taking into account only the number of vertices and sides.

In Figure 16**b** in the left column we have 12 ways of cutting with one knife, instead of 5. But here in all cuts except $VV_{1,3}$ the position of the cut is important, depending on whether it is made above, through or below the centre of gravity of the pentagon:

- $SS_{1,2}$ has 3 possible ways of cutting (AI, AII, AIII)
- $SS_{1,3}$ has 5 possibilities (BI, BII, BIII, BIV and C)
- $VS_{1,3}$ has 2 possibilities (E and F)

This explains the difference between 5 ways of cutting for $d_1$ and 12 ways for $d_m$ giving a total of 17.

## 3.2 Projecting and rotating shapes

The methodology used above refers to cutting or drawing lines, from side to side, vertex to vertex, vertex to side, either containing or not, the centre of the shape. The shapes can be considered as the result of *projection*. The transformations *rotation* and *dilation* or *scaling*, show how all these shapes fit together. Focus is on two basic shapes, the invariant one, and the projected one. In section 2.8 the analytic presentations describing the knives (Eqs. 12 and 13), have parameters $\alpha$ and $\delta$, dealing with rotation and dilatation, respectively.

In Figure 12 onto a basic square, other shapes are projected (square, line or rectangle), related to divisors. For $d_1$ a line is projected (upper row; imagine a line by a laser cut); for $d_2 = 2$ a rectangle is projected (lower row; in case of $VV_{1,3}$ or an $SS$ cut through the centre the rectangle reduces to a line); for $d_m = 4$ a square is projected (middle row). The line can be projected in any way, in any position. In Figure 15**c** right column, a grey rectangle is projected onto the yellow square. This rectangle can be narrowed and rotated. In 15**c** (right column) a grey square is projected onto a yellow basic square (AI) and scaled (AI → AII → AIII), then a rotation AIII → BI, a scaling BI → BII, a rotation BI → C or a rotation BII → D. In Figure 17 this is shown as a projection along the ribs of a pyramidal cone (see also Figure 12 central row).

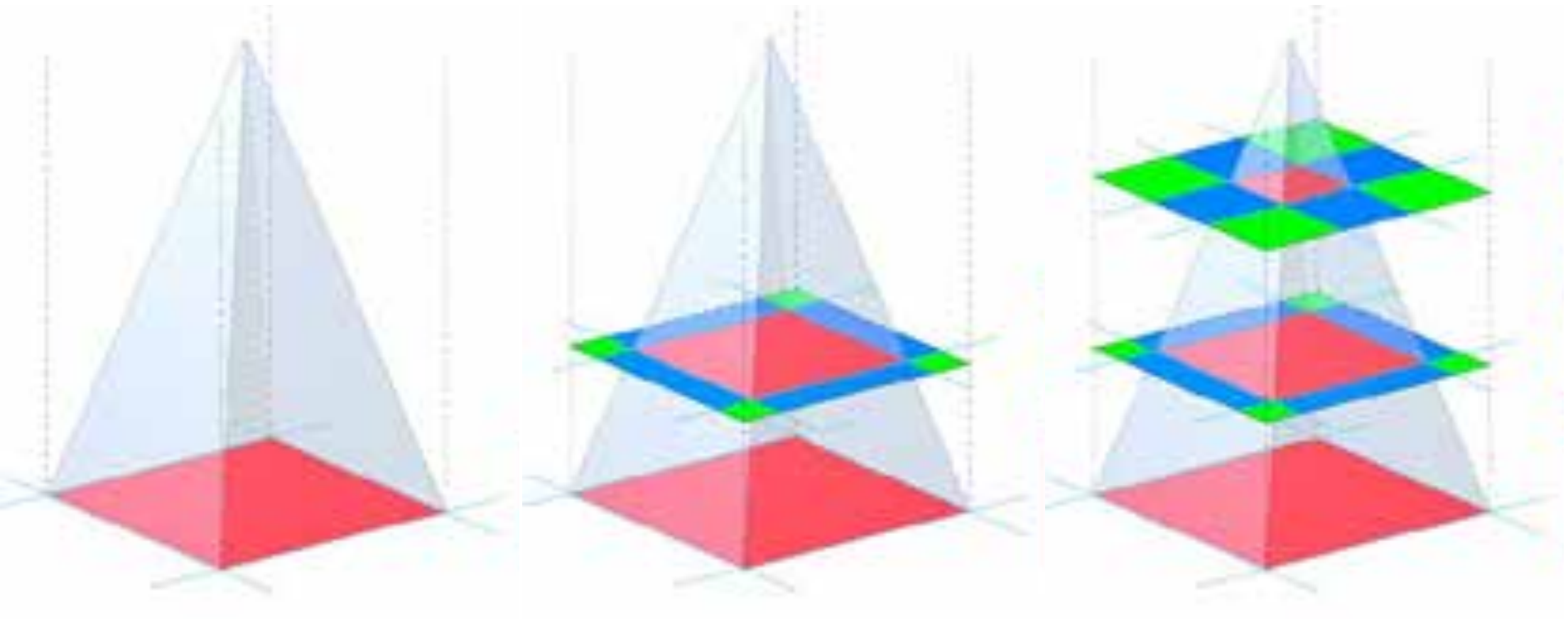

*Figure 17: Projections along a pyramid. The solid blue lines are cutting lines.*

For the pentagon (Figures 10 and 16) the same arguments can be made. In Figure 16**b**, onto a yellow base pentagon, a grey pentagon is projected (AI) and scaled (AI → AII → AIII → G → BIV). The grey pentagon is rotated 36° compared to the yellow pentagon. It can also have the same orientation as in BI and BII. Other instances are related to orientation.

This connects to the classical theory of conic sections, albeit that the cones can have a square cross section (a pyramid) or cross sections of any regular $m$-polygon, and this cone can then be moved through a fixed $m$-polygon. In the broader setting of $GML_m^n$ and *Gielis curves* there is no restriction to prisms. The base cross section, can be swept along a central line, or diminish in size, as in a cone (Figure 17), and need not to be constant.

## 3.3 Equidistant points on the circle

Since all Gielis curves, including self-intersecting curves, can be mapped onto the circle by a continuous Gielis transform, all vertices of an $m$-polygon can be mapped onto the circle. The vertices of a polygon are equispaced points on the circumscribing circle. The inner and outer vertices from a self-intersecting regular polygram, will become equi-spaced points on a circle, but they will be rotated relative to each other (Figure 18 for $m$=23). They can also be coinciding points. In this way, the same phenomena can be studied as points on a circle, using either algebraic approaches (groups and subgroups of the circle) or using geometric-algebraic approaches (Fourier, Chebyshev, complex numbers,….). All circles can be mapped onto one circle, or smaller and larger circles can be used (Figure 18). The self-intersecting Gielis curves, giving rise to multivalued functions, are directly related to Riemann surfaces. Note that a $VV_{1,2}$ cut is possible in the circle, whereas in regular $m$-gons such cut coincides with a side (it does not cut of a separate sector). $VV_{i,i+1}$ cuts (e.g. $V_1 \rightarrow V_2$) do divide the circle into two distinct parts.

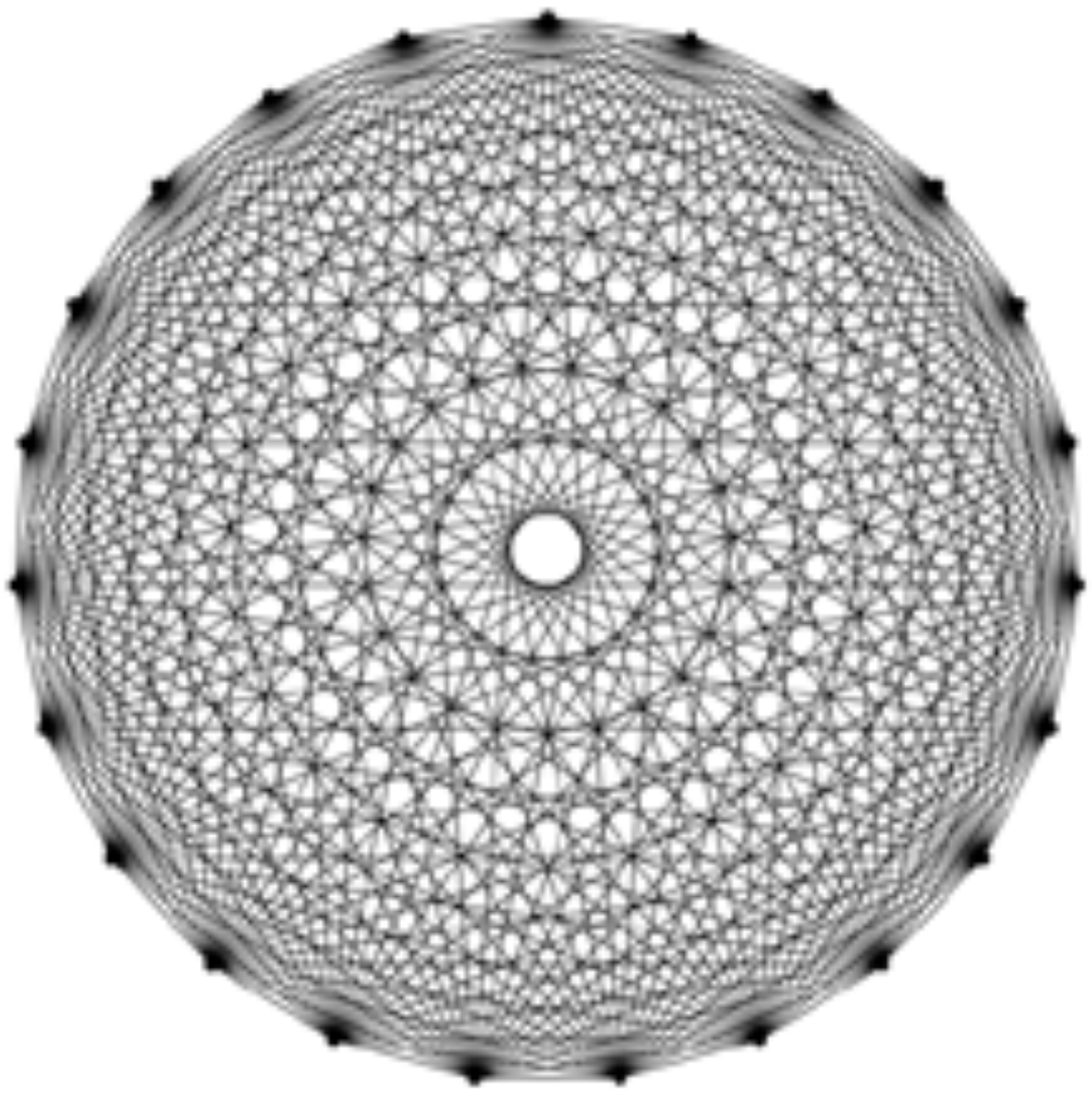

*Figure 18: Connecting equally spaced points on a circle (the circle itself is not shown, only connections; when m is even, connections pass through the center).*

## 3.4 Shortest path and curvature

To go from one point to another point on the curve, if the curve is the only possible way to go, one best follows the curve itself, and this then is the shortest path. Assume a circle, a Lamé curve with exponent $n = 2$, and an inscribed square with $n = 1$. Then the inscribed square, similar to the diagonals connecting four equally spaced points, or using all possible $VV_{i;\, i+1}$ cuts, can be considered as the knives cutting the circle. The straight cuts are the shortest possible ways, using the isotropy of the Euclidean plane.

However, if one considers Lamé curves with $1 < n < 2$ or $n < 1$, curves connecting the points result. They are curved lines, but using a broader definition of curvature, they are the shortest points to connect the two adjacent points on the circle.

From another viewpoint, the outer curve can have one geometry/curvature, while inside it can have another metric/measure/geometry/curvature, using another Lamé or Gielis curve. Now consider the three Lamé curves in Figure 19**a** as separate curves, and suppose that the red lines define a height function, e.g. they are elevated 1000m above the plane.

Then the only way to travel along the circle or any of two other curves safely is to follow precisely that curve. If the circle and the inside square are considered, the circle as the outer curve and the inscribed square as the four knives, both at elevation of 1000 meters, with deep down active volcanoes, tigers, lions and snakes, then it is possible to take different routes to travel from one point to the other, either via the square or along the circle. Connecting two of the four points is equivalent to drawing or cutting with a knife. In this case it is a straight knife, but if the curve is the one with $n = 2/3$, then both the route and the knife are curved, and optimal; following the curved path 1000 m above the plane will be the shortest and only possible path to stay alive anyway.

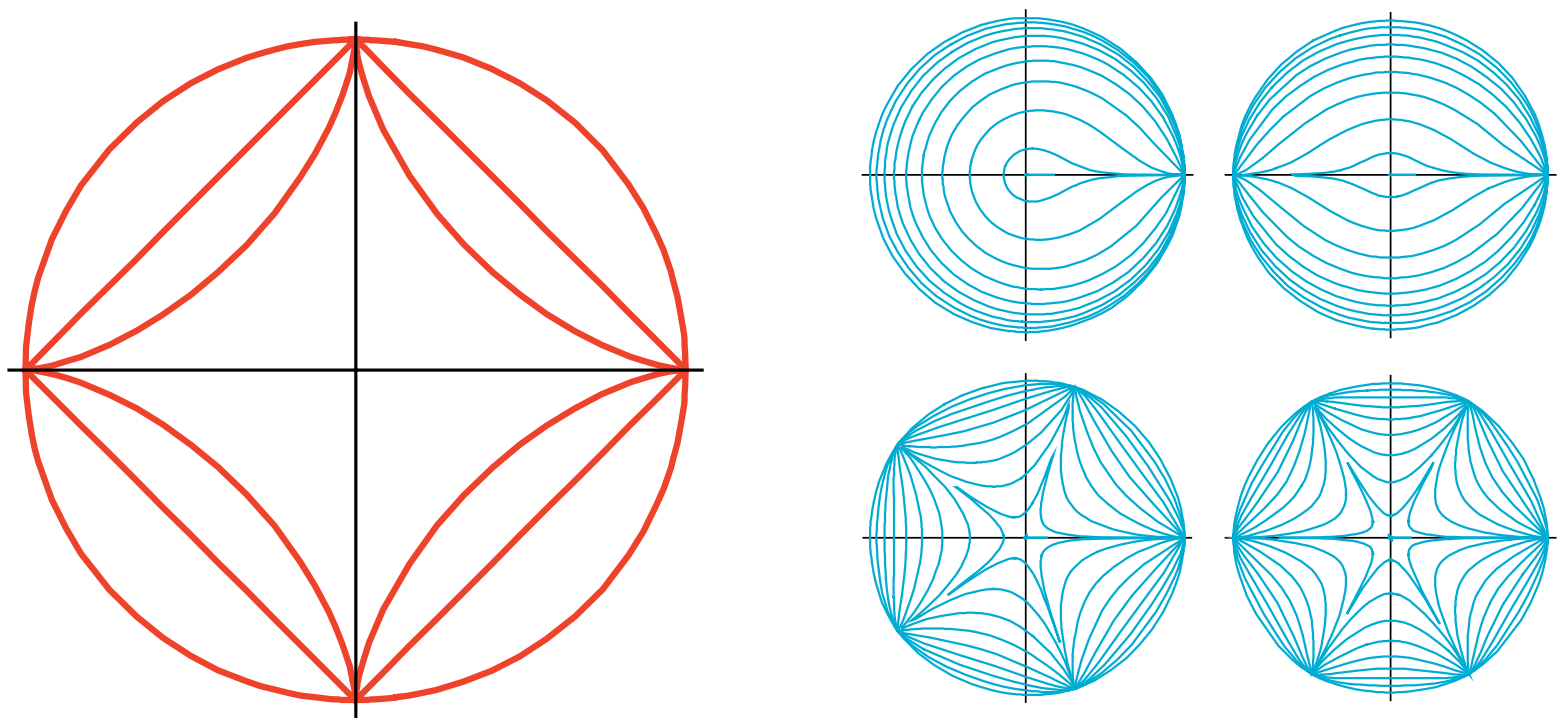

*Figure 19* ***a.*** *Lamé curves for n=2/3, 1 and 2.* ***b.*** *Gielis curves for* $m = 1,\ 2,\ 5,\ 6.$ *Each figure shows various curves for different exponents* [18].

Whether we speak of drawing lines to connect points, using knives to cut bodies or *m*-polygons, considering paths at 1000 m of height as the only way to avoid falling into volcano craters, they are all the same. To travel along the curves for *m*= 2 in figure 19**b**, from one antipodal point to the other, the paths will curved and cannot pass through the centre. Yet these are the shortest possible ways, in the same way as a ray of light follows a path around a gravitational body, which seems curved to us, but is the shortest path or geodesic anyway for the photons.

This inner curve is contained wholly inside the outer curve, but in this way one can define paths connecting two maximal or optimal distances in Gielis curves (Figure 19**b**). Consider a pentagon then one can find infinitely many curves inside the pentagon, which provide paths or curved knives. The same can be done for non-convex curves. In the classical view, a starfish is not a convex curve, but in our view this is less important, since the starfish can be mapped to a regular pentagon or a circle, two convex curves, and inner curves for optimal paths can always be found or defined. All in all, the precise shape of the knives may also be the result of some optimization problem.

Using a curved knife, the cutting can also be applied to concave Lamé-Gielis curves, in the following way: In Figure 20 supercurves were plotted with a HP7475A from the 1980's [19], for given $m$ (e.g. for $m = 4.3$ in Figure 20). The results correspond to self-intersecting Rational Gielis curves, with $m = p/q$, with $p$ a prime number; e.g. $p = 43$ providing curves with $m = 43/10 = 4.3$ with a central hole. The curved lines are now considered as curved knives, following the curvature of the space, thus generating the shortest possible path to connect two points on a circle with equally space points. It is a small step to see how cutting with curved knives corresponds to curved chords and chord diagrams[4].

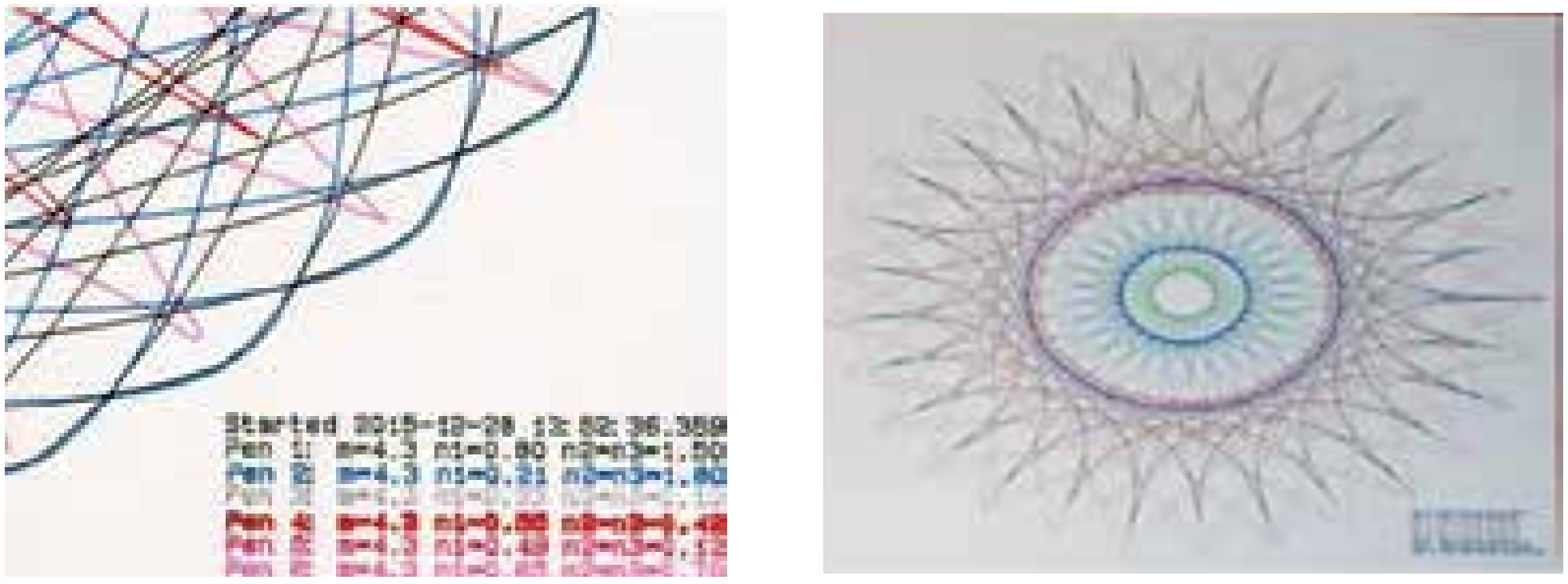

*Figure 20: Gielis curves created with an HP7475A plotter* [19]. ***a.*** *details for a polygon with m = 43/10 with different colours for different exponents.* ***b.*** *the complete curves.*

Having a Pythagorean compact analytic description for both shape and curvature for a wide range of shapes allows for a generalization of the notion of curvature. In [11] a new measure of curvature was introduced, directly related to the shape itself, using Gielis curves as osculating curves, as a generalization of studying curvature with circles.

## 3.5 Rational Gielis curves, R-functions and Flat tori

### 3.5.1 Rational Gielis Curves

The self-intersecting curves for any rational $m$ lead to various sectors in the polygons or cross sections of the *GML* body (Figure 21). Figures 21**a,c** have three zones, while Figure 21**b** has four different zones of different shapes indicated with 4 different colours. In *GML* bodies, when cut and separated, these zones represent different bodies.

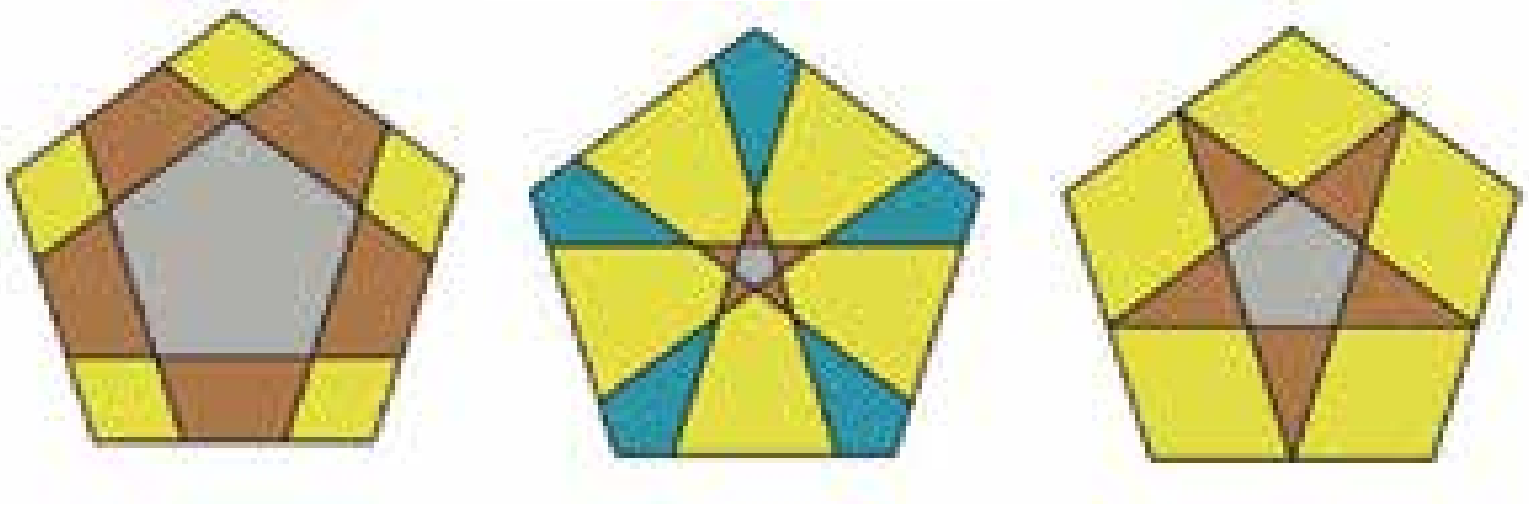

***a*** ***b*** ***c***

*Figure 21: Cutting pentagons from side to side*

[4] Or straight lines in a Poincaré disk.

Self-intersecting Gielis curves can represent the same for planar graphs [20]. For rational $m = p/q$, the number of zones created in Gielis curves is determined by $q$, and the symmetry of the polygons/polygrams is determined by *p*. In Figure 22**a**, a *7/5* RGC or polygram is shown, having 7-symmetry, closing after 5 rotations if drawn in one line.

In Figure 22**a** five different layers or zones can be defined in different shades of blue. Layers $L_0$ to $L_4$ are defined as a combination of layers from inside to outside and all layers have 7 maxima and 7 minima. A ray drawn from the centre *0* in any direction has multiple values indicated by $I_0$ to $I_4$ (red dots). When rotating the ray around the centre, the values of $I_0$ defines the boundaries of $L_0$ and the ray then sweeps the full area of $L_0$. Values of $I_0$ and $I_1$ define the boundaries of $L_1$, and here $I_0$ and $I_1$ coincide at maxima for $L_0$ and at minima for $L_1$. In the same way, values of $I_i$ and $I_{i+1}$ define layer $L_{i+1}$.

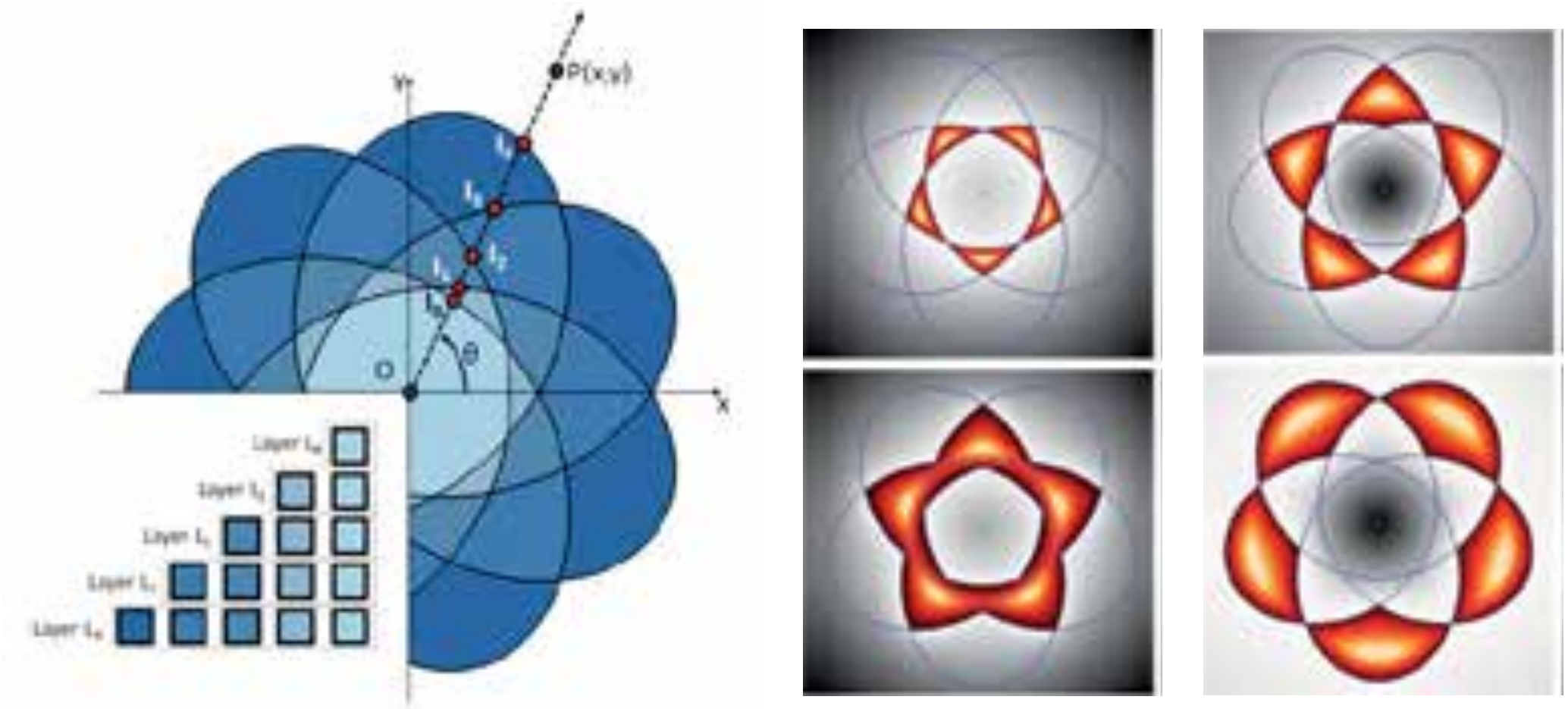

**a** **b**

*Figure 22* ***a.*** *Different layers in Rational Gielis curves.* ***b.*** *RGC for p=5 and q=4 with different zones defined* [20]

The regular polygons and polygrams in the regular polygon and *GML* cutting can be considered as Gielis curves and can be transformed into the shapes in Figure 22**b**. Now, *zones* can be defined, not only as stacked layers *L* in Figure 22**a**, but as separate layers or combinations of layers. We define $l_i$ as separate zones based on the different tones of blue zones in $L_{0,\dots,4}$ in Figure 22**a**. Examples of separate zones or combinations are given in Figure 22**b**; clockwise, from upper left (with $L_0 = l_0$):

- $l_1$ is $L_1 - L_0$
- $l_2$ is $L_2 - L_1$
- $l_3$ is $L_3 - L_2$
- $l_1 + l_2$ is $L_2 - L_0$

The zones and the independent domains separated by lines correspond to self-intersecting Gielis curves. In Table 6 the results of cutting GML's is compared with self-intersecting Gielis curves for symmetry $m = 5$ (Figure 14).

The order of cutting of *GML* or polygons is rearranged, showing the correspondence with self-intersecting curves. The number of independent figures and link number aligns with $q$ and the number independent shapes that can be counted (separated by lines) is a function $f$ of $p$.

In Figure 22**b** Layer $L_1$ has 5 independent shapes, plus the central pentagon. Layer $L_2$ has 11 $(= 5 + 5 + 1)$ independent shapes or *2p+1.* In general $f(p) = ap + 1$, with $a = 0, 1, 2, 3$ and $a + 1 = q$. For $a = 0, q = 1$, and the shapes close in one rotation, for all regular polygons. For $q \neq 1$ the number of rotations to close the shapes is $q$ and the shapes are then $2q\pi$-periodic. The zones $l$ shown in Figure 22**b**, align with the self-intersecting curves in Table 6 with the pentagram $m = 5/2$ and the polygrams $m = 5/3$ and $m = 5/4$.

*Table 6: GML versus Rational Gielis Curves RGC with $m = p/q$*

| Figure 16b | Cut | Number of independent figures in GML | Number of shapes in figure | Link | Number of shapes in figure | *p=5* | *q* & number of layers in RGC | *p/q* |
|---|---|---|---|---|---|---|---|---|
| A I | $SS_{1,2}$ | 2 | **6** | L2 | **6** | p+1 | 2 | 5/2 |
| A II | $SS_{1,2}$ | 2 | **6** | L2 | **6** | p+1 | 2 | 5/2 |
| A III | $SS_{1,2}$ | 3 | **11** | L3 | **11** | 2p+1 | 3 | 5/3 |
| B I | $SS_{1,3}$ | 3 | **11** | L3 | **11** | 2p+1 | 3 | 5/3 |
| B II | $SS_{1,3}$ | 3 | **11** | L3 | **11** | 2p+1 | 3 | 5/3 |
| D | $VS_{1,2}$ | 3 | **11** | L3 | **11** | 2p+1 | 3 | 5/3 |
| G | $VV_{1,3}$ | 3 | **11** | L3 | **11** | 2p+1 | 3 | 5/3 |
| B III | $SS_{1,3}$ | 4 | **16** | L4 | **16** | 3p+1 | 4 | 5/4 |
| B IV | $SS_{1,3}$ | 4 | **16** | L4 | **16** | 3p+1 | 4 | 5/4 |
| E | $VS_{1,3}$ | 4 | **16** | L4 | **16** | 3p+1 | 4 | 5/4 |

### 3.5.2 Boolean operations and R-functions

The number of independent shapes in Table 6 results from counting the isolated zones after cutting, delineated by lines. From a *GML* point of view, they are not independent, but are connected as sectors of a cut *GML*, according to a certain cutting process (Figures 15 and 16 with link numbers). In Gielis curves, these independent figures can be connected into layers or zones (Figure **22a,b**) using *R*-functions [21-23]. This is a method to translate Boolean operations into geometrical figures, and allows for, amongst others, for blending or ensuring differentiability up to some order [20-23]. Basic Boolean operations can be used to define the different layers and separate sectors as in Figure 22 and these can be translated into geometric language [11; 20]..

In Table 7 the most commonly used Boolean and associated *R*-functions are shown. The relationship to Lamé curves (from which the Gielis curves were derived) is clear via $x_1^p + x_2^p$ (and special cases of Lamé curves; the superellipses $|x_1|^p + |x_2|^p$ and circle for *p=2*): The use of exponents ${}^{m}/_{2}$ and ${}^{1}/_{p}$ converts Lamé curves into distances. In table 7 we find $\langle (x_1^2 + x_2^2)^{\frac{m}{2}} \rangle$ or, a special Lamé curve $((x_1^2 + x_2^2)$, raised to exponent $\frac{m}{2}$. This is the Pythagorean part distance $(x_1^2 + x_2^2)^{\frac{1}{2}}$ raised to power *m.* The function $\mathfrak{R}_p$has $\langle [|x_1|^p + |x_2|^p]^{\frac{1}{p}} \rangle$. Originally $p$ is a positive integer but can now be any positive number $p > 1$ [20].

This leads to the remarkable fact that separated regions in *RGC* (and in polygrams in general) can be defined as coherent structures. Using *R*-functions, these operations are translated into geometry, whereby the different layers (Figure 22**a** in blue) or different sectors (the red zones in Figure 22**b**) are defined as single geometrical domains or combinations of single domains [20].

*Table 7: Examples of R-functions for the most common Boolean functions.*

| $\mathfrak{R}_0^m$ | $\mathfrak{R}_p$ |
|---|---|
| $y_1 \equiv -1$ | $y_1 \equiv -1$ |
| $y_2 \equiv \bar{x} \equiv -x$ | $y_2 \equiv \bar{x} \equiv -x$ |
| $y_3 \equiv x_1 \wedge_0^m x_2 \equiv \left(x_1 + x_2 - \sqrt{x_1^2 + x_2^2}\right)(x_1^2 + x_2^2)^{\frac{m}{2}}$ | $y_3 \equiv x_1 \wedge_p x_2 \equiv x_1 + x_2 - [\lvert x_1\rvert^p + \lvert x_2\rvert^p]^{\frac{1}{p}}$ |
| $y_4 \equiv x_1 \vee_0^m x_2 \equiv \left(x_1 + x_2 + \sqrt{x_1^2 + x_2^2}\right)(x_1^2 + x_2^2)^{\frac{m}{2}}$ | $y_4 \equiv x_1 \vee_p x_2 \equiv x_1 + x_2 + [\lvert x_1\rvert^p + \lvert x_2\rvert^p]^{\frac{1}{p}}$ |
| $y_5 \equiv x_1 x_2$ | $y_5 \equiv x_1 \sim_p x_2 \equiv \frac{x_1 x_2}{[\lvert x_1\rvert^p + \lvert x_2\rvert^p]^{\frac{1}{p}}}$ |
| $y_6 \equiv x_1 \rightarrow_0^m x_2 \equiv \left(x_2 - x_1 + \sqrt{x_1^2 + x_2^2}\right)(x_1^2 + x_2^2)^{\frac{m}{2}}$ | $y_6 \equiv x_1 \rightarrow_p x_2 \equiv x_1 - x_2 + [\lvert x_1\rvert^p + \lvert x_2\rvert^p]^{\frac{1}{p}}$ |
| $y_7 \equiv x_1 /_0^m x_2 \equiv \left(\sqrt{x_1^2 + x_2^2} - x_1 - x_2\right)(x_1^2 + x_2^2)^{\frac{m}{2}}$ | $y_7 \equiv x_1 /_p x_2 \equiv -x_1 - x_2 + [\lvert x_1\rvert^p + \lvert x_2\rvert^p]^{\frac{1}{p}}$ |
| | $p > 1$ |

### 3.5.3 Knots and Graphs

In this way we can connect seemingly unconnected zones in self-intersecting graphs in the very same way as the zones are connected when executing the cutting in 3D. Using *R*-functions the zones can be defined in planar graphs, which are connected as in 3D *GML* bodies. In *GML* bodies the resulting figures form separate bodies and in cross sections they show up as differently coloured, separated zones. In planar polygons, they can be connected via *R*-functions, but also by drawing graphs occurring in the cut polygons. Consider for example Figure 23. The boundaries of the yellow zones are defined by the circumscribing pentagon and by the inscribed pentagram (enclosing the brown and grey zone). Connecting the yellow zones, by using the boundaries of the innermost pentagram in grey, and considering only the boundaries, a graph is obtained.

The boundaries of zones 1, 2, 3, 4 and 5 (individual zones separated by small perpendicular blue lines), and connected by green highlighted lines, form a graph. This graph can be also obtained starting from the uppermost vertex of the inner pentagon with green highlighted sides, going around sector $S_1$, then continuing along the green line to enclose $S_5$, then to $S_4$, then $S_3$, and $S_2$, ending this graph exactly at the upper vertex where the trajectory started. The same can be done for the brown sectors, and for the grey pentagon.

The boundary of graph enclosing the inner pentagon closes in one rotation, while the other two graphs self-intersect. *These graphs are homeomorphic to the planar projection of a $5_1$ knot (or RGC with m=5/2) figure 23**b**, or to Figure 23**c**.*

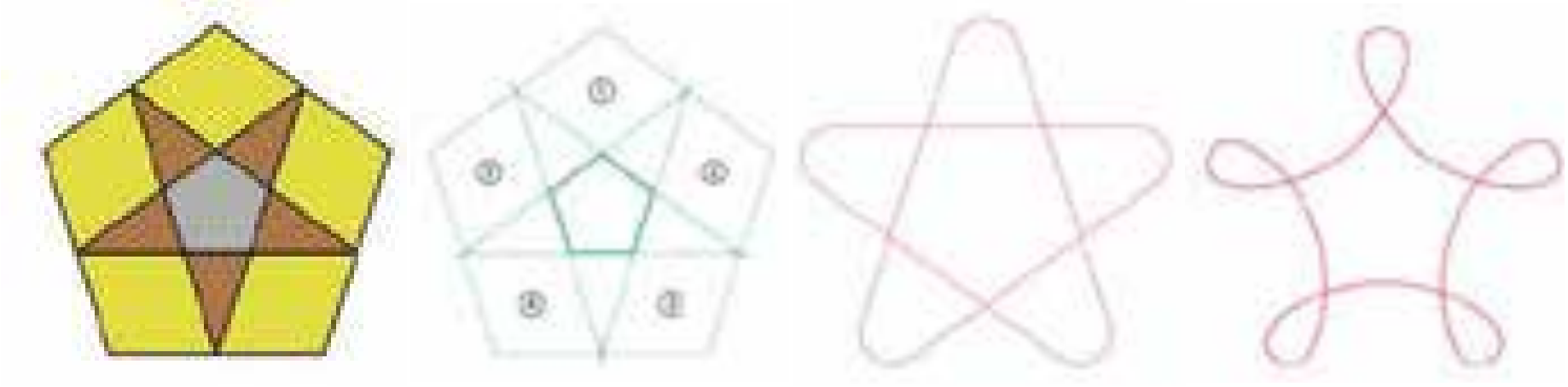

*Figure 23 **a.** cross section of $GML_5^n$ for $SS_{1,3}$ cut, case $BII$ **b.** generation of graphs for case $BII$. **.c.** 5/2 knot, **d.** graph with loops corresponding to sectors 1-5 in 23**b**.*

| | | | | | |
|---|---|---|---|---|---|
| BII. | 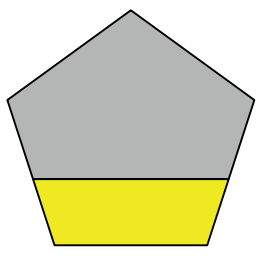 | $S_{1,3}$<br>$b_3^1 + b_3^2 = A$ | 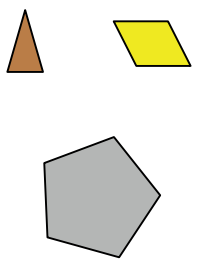 | $GML_3^{15\omega + 3q + 12}$ $\{\omega + 0.2q\}$<br>$GML_4^{20\omega + 4q + 16}$ $\{\omega + 0.2q\}$<br>$GML_5^{25\omega + j}$ 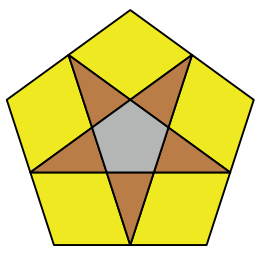 | link-3 |

*Figure 24. Resulting structures and bodies from BII cut of $GML_5$* [5]

These graphs also relate to the process of cutting $GML$ in the following way. When cutting $GML$ with a pentagonal cross section (as example), according to $S_{1,3}$ cut, case $BII$, the result is a Link-3 of three separate bodies. The inner body is pentagonal, the second body has a triangular cross section and the outer one has a quadrangular cross section. In Figure 24 it is indicated how these bodies $GML_3, GML_4$ and $GML_5$ are twisted. For $\omega = 1$ and $j = 0 = q$ the result is that the $GML_5$ body is twisted 25 times (five full rotations), the $GML_3$ body is twisted 27 times (9 full rotations), and the $GML_4$ body is twisted 36 times (also 9 full rotations). The overall structure is a link of 3 bodies, one of which ($GML_5$) is the unknot with circle as basic line, while $GML_3$, is a $5_1$knot, with the pentagram (or $RGC$ with $m = 5/2$) as basic line.

The study of $GML$ bodies and surfaces revealed complicated structures on interlinked bodies with knot structures, in particular the simplest knots $k/2$ (with $k$ integer) or torus knots. These results also appear in the methodology of cutting planar regular $m$-polygons and in the methodology of self-intersecting Gielis curves and graphs (at least for $\omega = 1$ and $j = 0 = q$ as in the $B\ II$ case in Figure 24).

As a consequence, when cutting $m$-regular polygons, the link number and knot characteristics of the original $GML$ from which they were cut can be found, counting the number of independent sectors, their shape and coherence as layers. Their coherence is established via $R$-functions, named after the Ukranian mathematician V.L. Ravchev (1926-2005).

## 3.6 Galilei, Euler, Catalan and Cayley

### 3.6.1 The problem of dividing *m*-polygons with non-crossing diagonals

Euler considered the following problem: in how many ways can one divide a polygon into triangles, in such a way that diagonals do not cross, i.e. similar to cutting with non-crossing knives [24]. For a pentagon using two diagonals this gives a total of 5 possibilities (Figure 25). For a square, only two diagonals are possible.

In the case of a pentagon, the cuts that are allowed are symmetric ones, $VV_{1,2}$ and $VV_{1,3}$. Using diagonals means vertex-to-vertex cuts and using the same notation as before these cuts are $VV_{1,3}$ and $VV_{2,4}$ in the square, giving 2 possibilities (Figure 26**a**). The upper square gives $VV_{1,3}$ (or $VV_{2,4}$) while the lower one gives both possibilities in superposition. In the case of a pentagon, $VV_{1,3}$ and $VV_{1,4}$ diagonal knives (Section 2.10) can be used. For the *GML* case they are symmetrical, but for the Euler problem, they are different (Figure 25). In Figure 26**b**, in the lower pentagon, all possible cuts are given with $VV_{1,3}$ and $VV_{1,4}$, starting from all five vertices. So the number of possible ways of cutting is 5. The lower pentagonal figure in Figure 26**b** can be considered as a composition of all $VV_{i,i+2}$ and $VV_{i,i+3}$ cuts of Figure 25.

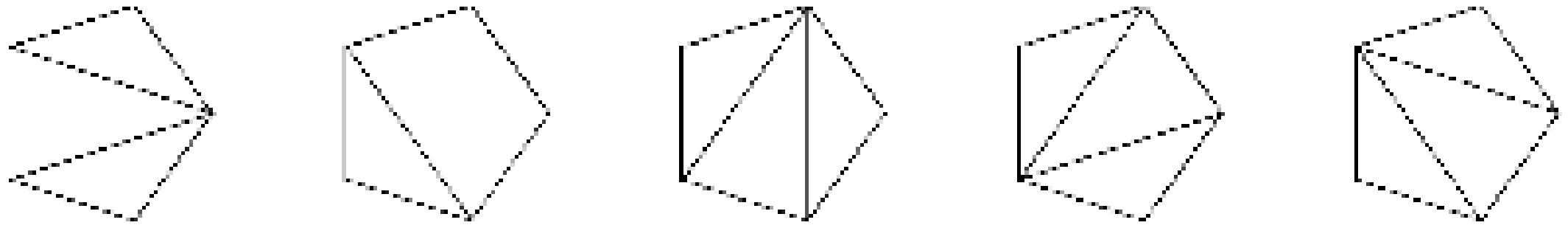

*Figure 25: Five ways of cutting a pentagon, with diagonal knives*

In the hexagon (Figure 26**c**) the possible cuts are also related to $VV_{i,i+2}$ and $VV_{i,i+3}$ cuts, and the total number of possibilities, is obtained when superposing the two lower hexagons. Possible cuts with $VV_{i,i+3}$ and $VV_{i,i+2}$ give the lower left and lower right hexagon, respectively. The number of possible cuts is 14 and the superposition of the lower two hexagons gives all possibilities. In other words, by superposition of the two hexagons, all possible individual paths for the Euler problem can be traced.

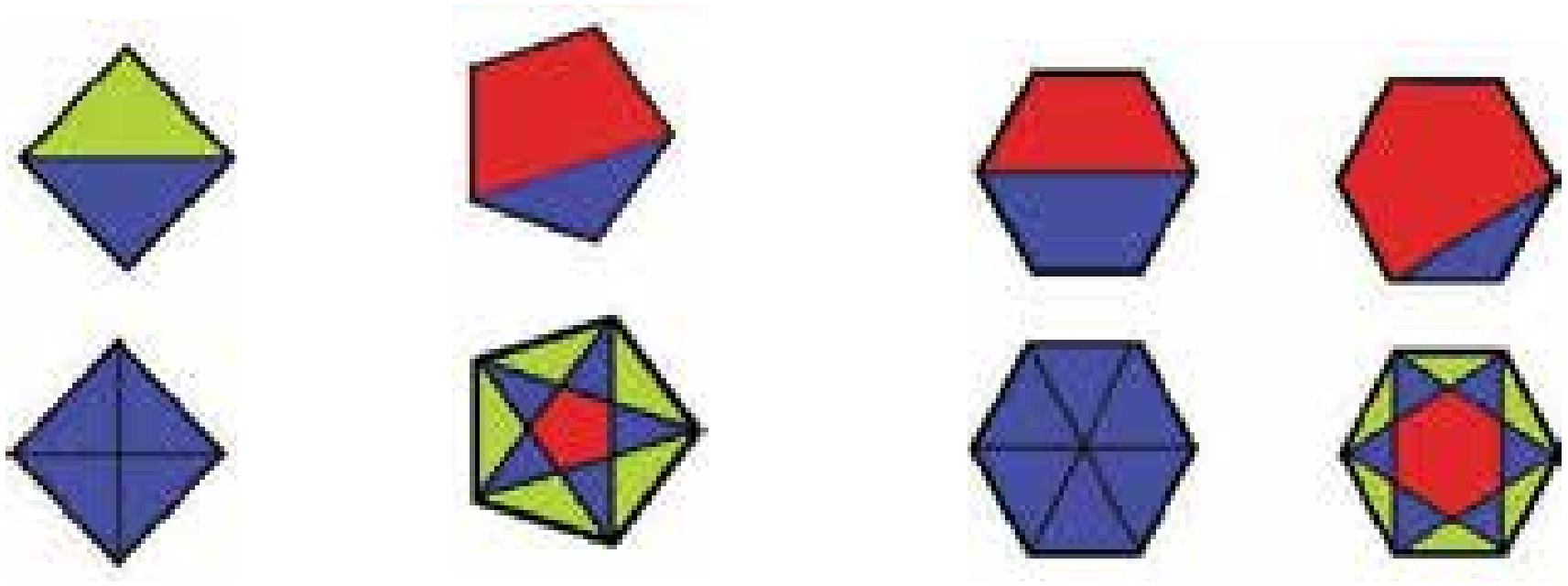

**a** **b** **c**

*Figure 26: Cutting square, pentagon and hexagon with diagonal knives*

If one observes the upper row Fig26**b**, a division by a $VV_{1,3}$ cut (divisor 1) gives a triangle and a pentagon. The latter can be divided then into triangles, by $VV_{1,4}$ and $VV_{1,5}$ cuts. In Figure 26**c**, a $VV_{1,4}$ divides the hexagon into 2 quadrilaterals, which can further be cut along their diagonals. Or one can start from $VV_{1,4}$ cuts, giving two triangles and one rectangle. This shows that this can be considered as a recursive procedure.

In heptagons $VV_{i,i+2}$ and $VV_{i,i+3}$ cuts are possible, and again, the superposition of the two lower figures in Figure 28 left, would give all possibilities. For octagons, $VV_{i,i+2}$, $VV_{i,i+3}$ and $VV_{i,i+4}$ are possible diagonals. Overlaying all three lower octagons would give all possible solutions. In figure 28 this is shown for the octagon, leading to 132 different possibilities.

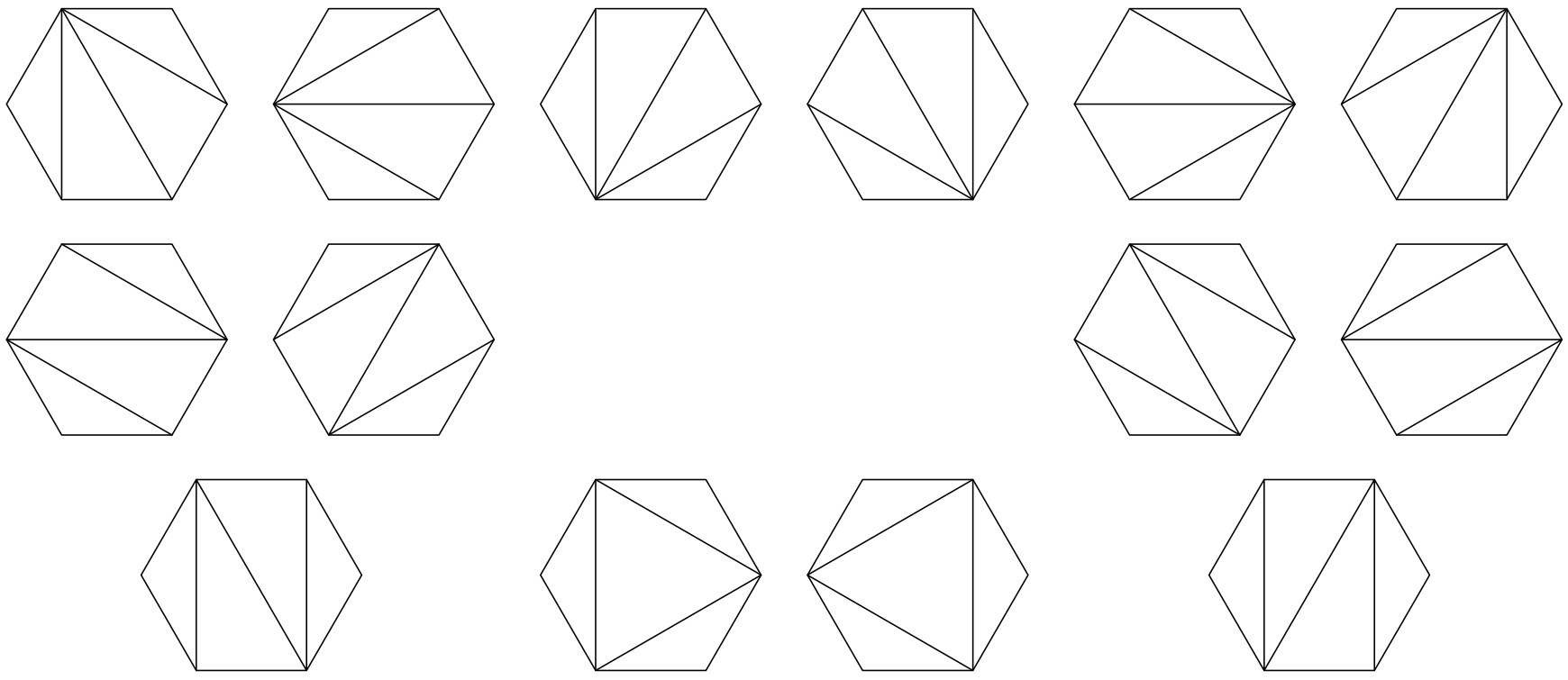

*Figure 27: 14 ways of cutting a hexagon, with non-crossing diagonals*

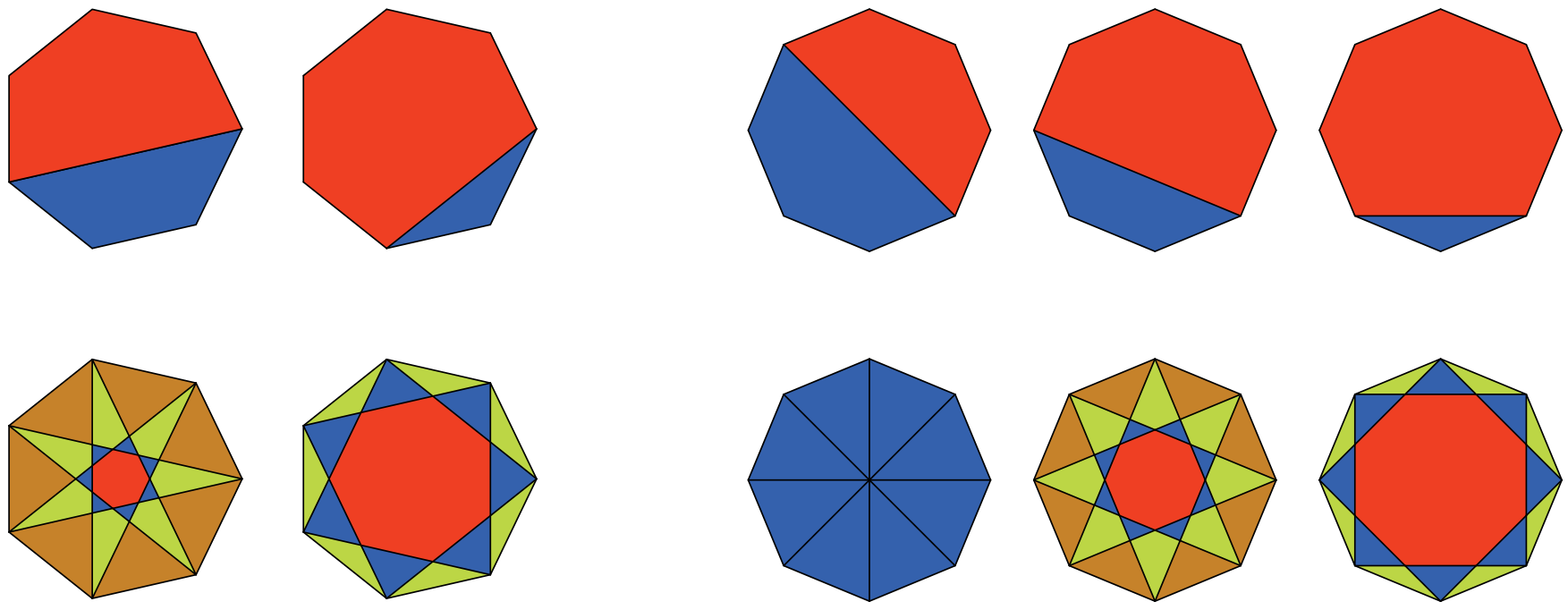

*Figure 28: Cutting heptagons and octagons with diagonal knives*

With superposition of the figures, diagonal knives seem to cross, but this is only superficially so, since the operations are carried out independently in *GML* rotating through one fixed knife. The question or problem that could be asked is the following:

Describe all possible ways of cutting *m*-polygons into triangles by using the allowed *VV* cuts. This can be carried out in 3D by using one knife each to cut the bodies. Each of the cuts corresponds to *d*-knives. The result will be the lower set of hexagons. Superposition of figures will give all possible paths.

The solution will be obtained in one operation if the process is executed in 3D by using three fixed knives simultaneously (example for the hexagon) if one moves a *GML* body through the knives.

### 3.6.2 Euler, Segner, Lamé and Catalan

The total number of possible cuttings with non-crossing knives of an $m$-gon, the problem originally studied by Euler [24; 25] is given in Table 8 for $m = 5, \dots 10$. The number $r$ in the second row gives with the number of vertices $r + 2$, starting from $r = 3$, for a pentagon. This sequence of numbers, starting *1, 2, 5, 14, 42,....* are Catalan numbers, which play a very important role in many fields of mathematics. The number $r$ in Table 8 is the position in the row of Catalan numbers that start with 1 and 2.

*Table 8: Possible cuts for m-gons*

| | | | | | | |
|---|---|---|---|---|---|---|
| *m-gon* | 5 | 6 | 7 | 8 | 9 | 10 |
| $r = m - 2$ | 3 | 4 | 5 | 6 | 7 | 8 |
| *Total N°* | 5 | 14 | 42 | 132 | 429 | 1430 |

Euler found a solution to the problem in 1751. Euler's number of possibilities $E_m$:

$$E_m = \frac{2.6.10\ldots(4m-10)}{(m-1)!} \qquad (14)$$

In 1758 the Hungarian mathematician Segner found a recursive formula. Let the vertices of any convex $m$-gon be *1, 2, 3,..., m.* The triangle $\Delta 1rm$ divides the polygons into an $r$-gon and an $m + 1 - r$ gon (Figure 29). From this Segner's recursive formula follows for $(E_2 = 1)$:

$$E_m = E_2 E_{m-1} + E_3 E_{m-2} + \cdots + E_{m-1} E_2 \quad (15)$$

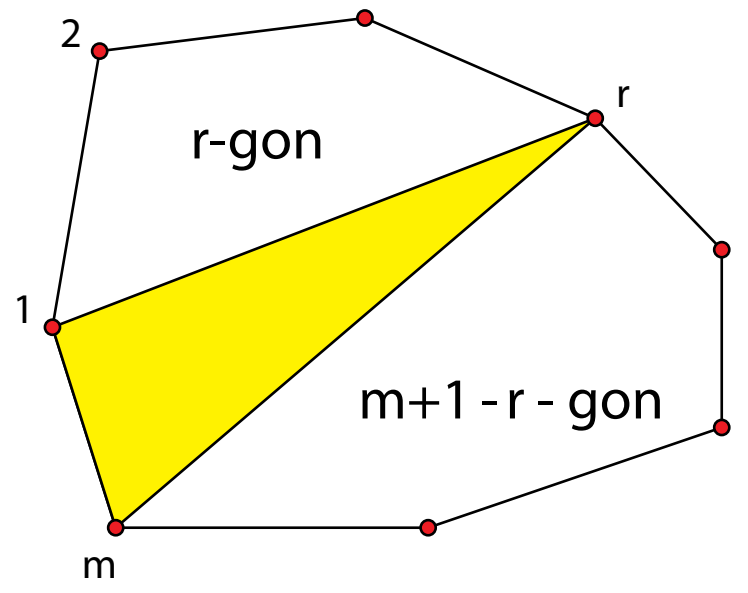

*Figure 29: Dividing an m-gon into a triangle, r-gon and $(m + 1 - r)$-gon*

For his newly established journal Joseph Liouville in 1838 asked proofs to combine the solution of Euler and the recursive one of Segner. Liouville published the geometric proof of Gabriel Lamé as the most elegant [26]. Later the same journal published solutions of Olinde Rodrigues and Eugene Catalan. They linked the sequence of numbers more to combinatorial problems, for example the way in which a sequence of operations can be carried out using brackets. A recent book collects over 200 problems and questions in which Catalan numbers appear [27].

### 3.6.3 Euler, Liouville and Cayley

This is related to many other problems, which has involved many important mathematicians, such as Euler, Liouville and Cayley. One example: *For $m$ equally spaced points and all possible connections / diagonals drawn how many separate parts result?* Examples are given in Figure 30 and the answer is OEIS *A006533.* The drawings once more correspond to our overlayed figures. It is possible to draw figures for any $m$ via [*https://oeis.org/A006561/a006561.html*](https://oeis.org/A006561/a006561.html).

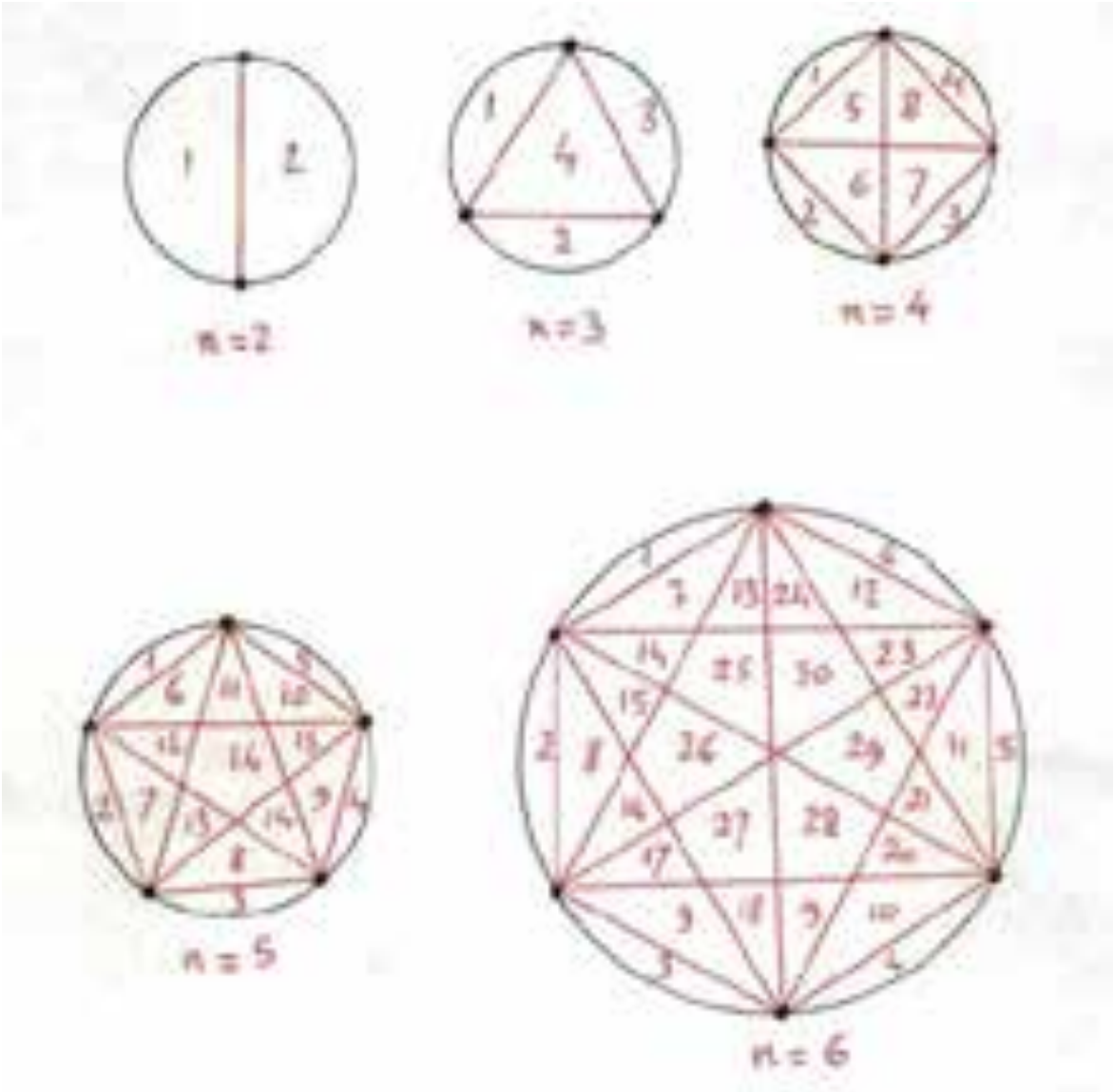

*Figure 30: Original drawings of Jos Meeus (*[*https://oeis.org/A006561/a006561.html*](https://oeis.org/A006561/a006561.html)*). For $m = n = 5$ the original Pythagorean pentagon and pentagram result.*

For even $m$, diagonals of opposite points cross in the centre, while for odd $m$ a central $m$-polygon is created, with the same shape as the original polygon, but rotated relative to the original polygon (Figure 31). In *GML*/planar polygon strategy, diagonals through the centre can be created for odd $m$, by using a vertex to side cut, in particular a vertex to the middle of the opposite side, as for example in a pentagon. This is directly related to Table 6 and the sectors and layers in Rational Gielis Curves *RGC* (Section 3.5).

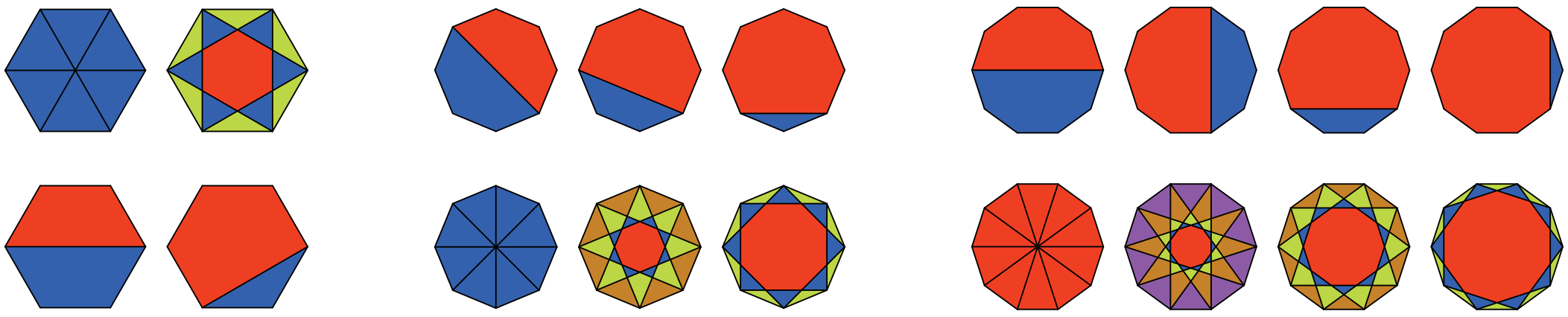

a. *m=6*  b. *m=8*  c. *m=10*

*Figure 31: Diagonals through centre for ($VV_{1;(m+2)/2}$ cuts) and central m-polygons created for all other VV cuts. For $m = 6, 8, 10$*

### 3.6.4 Cayley's coloured partitions

A more general question on coloured partitions of a convex polygon by non-crossing diagonals was initiated by Cayley. For any positive integers *a* and *b,* all coloured partitions made by non-crossing diagonals of a convex polygons into polygons, whose number of sides is congruent to *a modulo b* can be obtained by recurrence relations and an explicit representation in terms of partial Bell polynomials [28]. Other relations can be found with Motzkin numbers and Dyke paths. This is very much related to the Euler problem in 3.6.1 and Figure 27, whereby the inscribed convex polygons reduce the number of vertices on the original polygons. In Figure 32 all paths of such inscribed polygons can be traced.

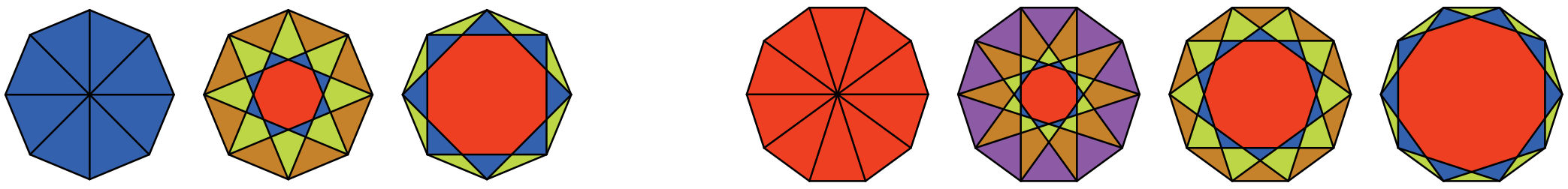

Figure 32: Rational Gielis Curves with *m*= 8/2, 8/3 and *m* = 10/2, 10/3, 10/4, 10/1

### 3.6.5 Galilei, cycloids and clams

The non-crossing diagonals also appear in a problem of rolling polygons. If a polygon is rolled on a line the vertices form a polygonal arch. If one wants to compute the area under the polygonal arch this can be done by forming sets of blue, green and pink triangles in Figure 33. These can be arranged in the original polygon and this proves that the area under the polygonal arch is three times the area of the generating polygon [29]. The arrangement of the triangles in the polygons is one of the solutions of the Euler problem of cutting a polygon with non-crossing diagonals.

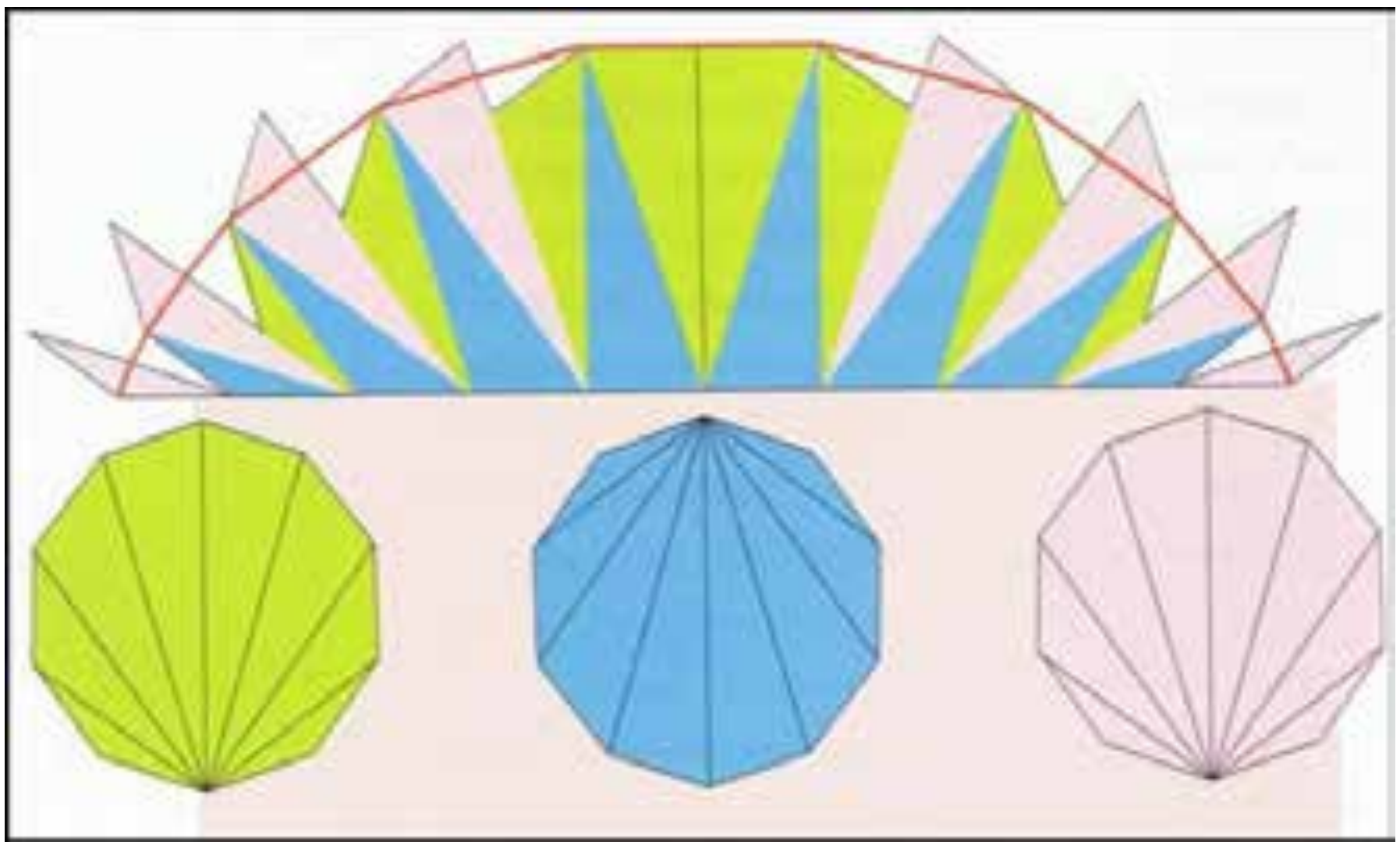

*Figure 33: Determining the area under a polygonal arch*

As a corollary, the area under a cycloid, generated by a fixed point on the circle as the rolling curve is three times the area of the circle. This corollary implicitly assumes that the number of sides of the polygon increases to infinity, with the circle as the limiting case.

With Gielis transformations however, we can keep the symmetry of the polygons fixed and morph them into circles by changing the exponents in order to have the transformation equal 1, i.e. exponents $n$ are 2, the Pythagorean case. Piecewise linear and continuous are seamlessly combined.

If one considers the diagonals starting from the top in the blue polygon in Figure 33, the vertices represent the points reached by a ball, either dropped (central diagonal) or rolling off an inclined plane, in a fixed amount of time $t_1$. Figure 34**a** is a figure from Galileo Galilei's famous book, in which he discusses the phenomenon of a ball under friction and the force of gravitation [30]. Points E, G and I are points reached after a time $t_1$. The greater circle, with points B, H and F is the locus for all points that are reached by a ball at a time $t_2 > t_1$. In nature such phenomena are observed in growth patterns in bivalve mollusks (Figures 34**b,c**), and plant leaves as already observed in the original book *Inventing the Circle* [10]. So in some sense, we have come full circle after almost two decades.

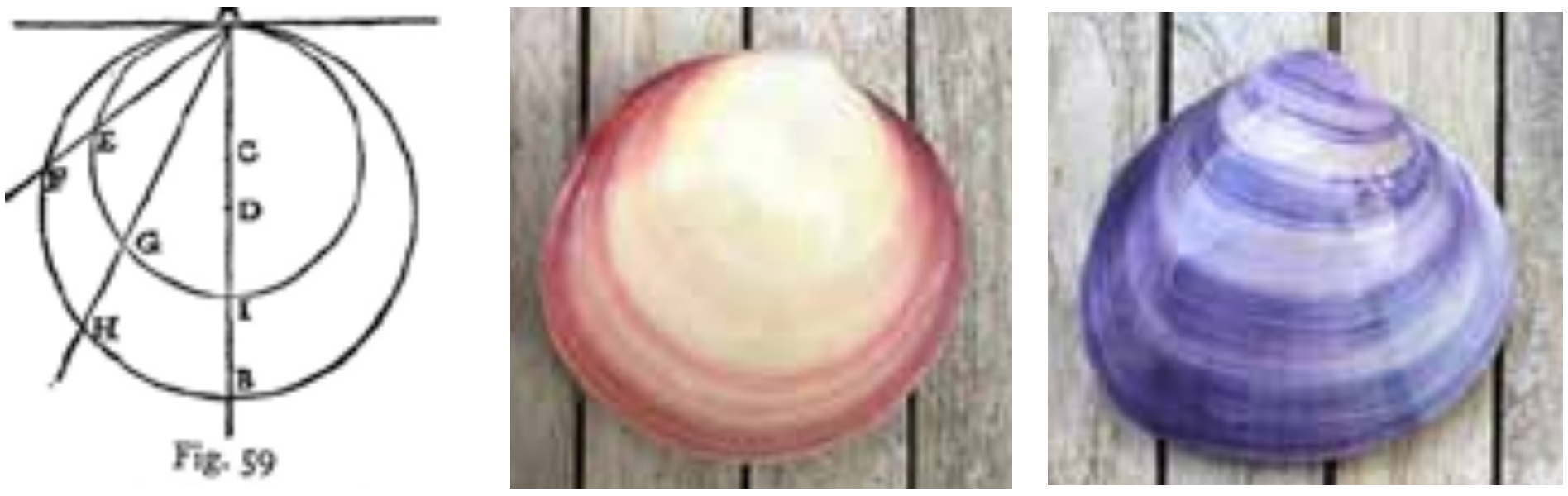

*Figure 34: **a**. Galilei's original figure, **b,c**. Growth rings in shells*

## 3.7 *GML* and Gielis curves as connectors

The theory of Generalized Möbius-Listing surfaces and bodies had been merged earlier with Gielis curves, either as cross section or as basic line. In the sections *Expositio* and *Determinatio,* it is shown that the connection between *GML* and planar Gielis curves runs much deeper, unveiling the relation between planar geometry, conic sections, *GML* surfaces and bodies, torus knots and planar knot projection, Flatland and branes, and more, and a wide range of methods and problems based on groups and equidistant points on a circle, revealing deep connections and equivalency of geometrical, arithmetical, logical, geometric-logical, combinatorial and algebraic methods.

In particular, a 3D problem can be solved by planar geometry but the problem of cutting *GML* bodies and surfaces is more general than related problems of non-crossing diagonals, or the number of sectors when connecting equally spaced points on a circle. Rotation allows for starting with one cut only and repeating this procedure, to find all possibilities of cutting. With rotation and scaling, *VV* cuts (classical diagonals) can become *VS* and *SS* cuts, or any possible cut dividing the $m$-polygon in exactly two parts with one knife. Solving the general case will allow to consider many other mathematical problems as parts of the general solution.

# 4 Constructio

## 4.1 Cutting a regular polygon

Consider a regular $m$-polygon. Vertices and sides are numbered from 1 to $m$. In Figure 35 the example of a pentagon is shown.

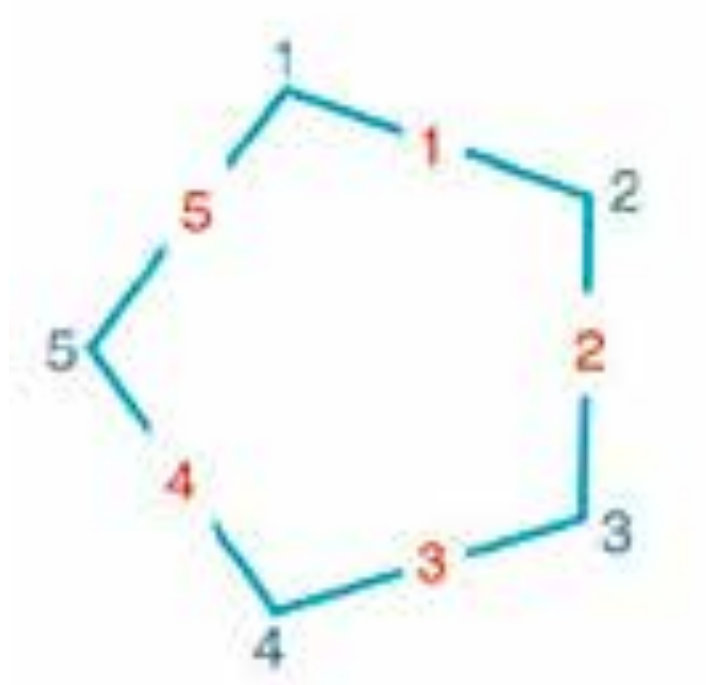

*Figure 35: Numbering of pentagon, vertices in blue, sides in red*

A regular $m$-polygon can be cut in various ways, from vertex to vertex (notation $VV_{ij}$, e.g. $VV_{1,3}$ from vertex 1 to vertex 3), from vertex to side (notation $VS_{ik}$, e.g. $VS_{1,3}$ from vertex 1 to side 3), or from side to side (notation $SS_{kl}$, e.g. $SS_{1,3}$ from side 1 to side 3). $VV_{max}$, $VS_{max}$ and $SS_{max}$ are the cuts from and to vertices or sides with maximal separation (*max* is not necessarily the longest length). For example in Figure 35 $VV_{max} = VV_{1,3} = VV_{1,4}$; $VS_{max} = VS_{1,3}$ through the centre, and $SS_{max} = SS_{1,3}$. If the line in $SS_{max}$ cuts side 3 in the middle, and the cut of side 1 is moved to vertex 1, then $SS_{max}$ can be made arbitrarily close to $VS_{max} = VS_{1,3}$.

The following general facts can be observed in Figure 36:

- In a polygon for even $m$, $VV_{max}$ goes through the centre of the polygon ($m = 6$), and in a polygon for odd m, $VS_{max}$ goes through the centre, crossing the side opposite the vertex in the middle (*m=7*).
- Going from a regular $m$-polygon to a *(m+1)*-polygon introduces one extra vertex and one extra side. From $m = 4$ to $m = 5$ a line converts into a wedge of $VV_{1,3} = VV_{1,4}$, giving one extra vertex and one extra side. Considering $VS_{max}$ a wedge is created in $m$ = even polygons ($m = 8$), and one line (through the centre) in $m$ = odd polygons ($m = 7, m = 9$)
- In $m = 10$, all cuts or diagonals are drawn from one vertex. Besides the $VV_{max}$ the other *VV* cuts are two by two symmetrical (solid and dashed lines). In the case of $m$ = odd, the same can be said for *VS* cuts ($m = 7$).
- A single cut divides an $m$-polygon into two parts, which are defined by their shape and number of vertices and sides. In $m = 4$ the square is divided by the red diagonal into two triangles. In $m = 6$ the hexagon is divided into two equal quadrilaterals or trapezoids. In $m = 11$ the polygons is divided into an octagonal shape with 8 vertices, and a pentagonal shape with 5 vertices.

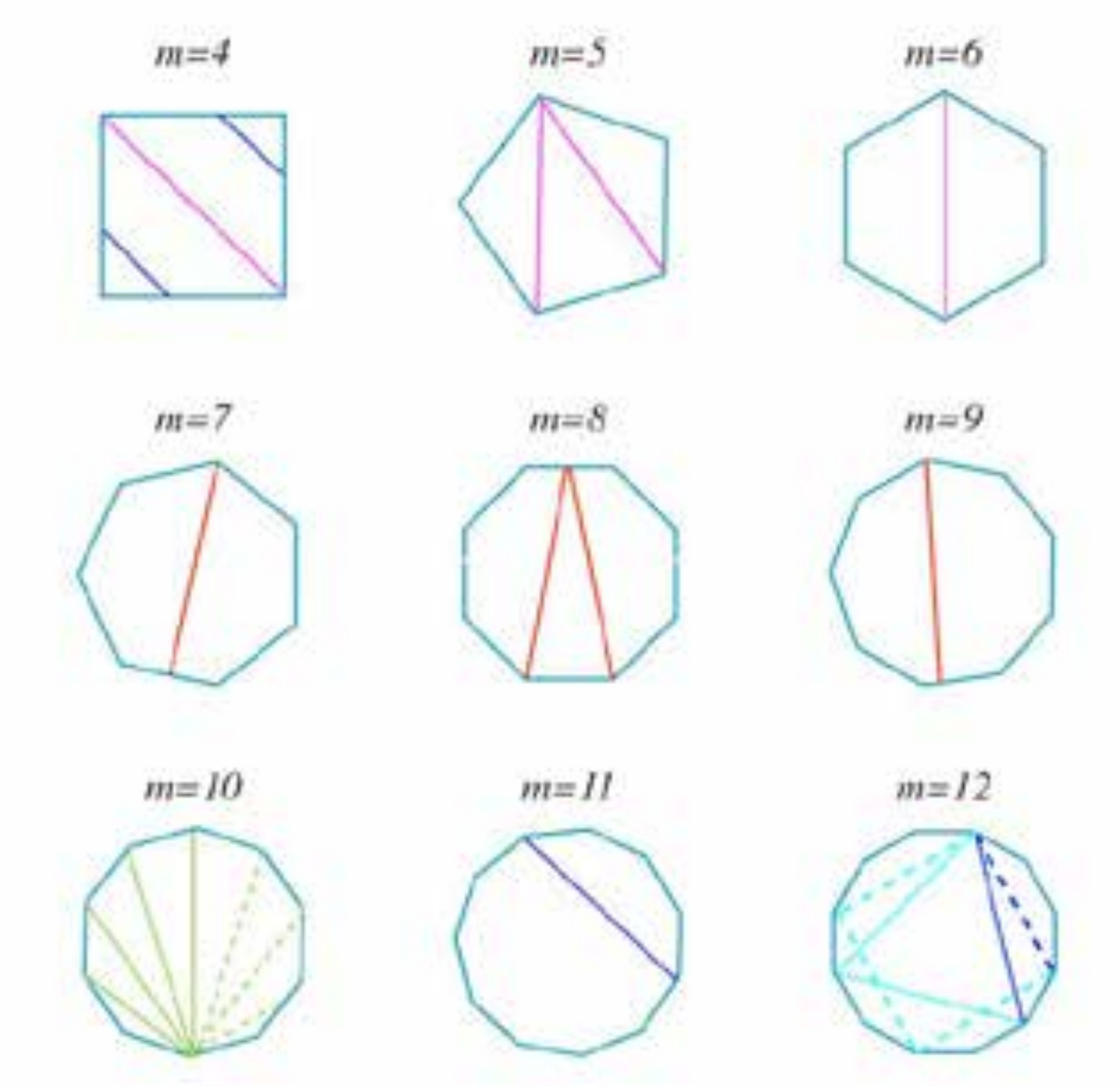

Figure 36: Understanding modes of cutting

## 4.2 Cutting and divisors of m

There is a definite relation between the ways of cutting and the number of divisors of $m$. In Figure 36, $m = 11$, a $VV_{1,5}$ cut is related to the smallest divisor of $m$, $d_{min}$=1. For $m = 12$:

- $VV_{1,4}$ cut (dark blue dashed line), is repeated as $VV_{4,7}$ and as $VV_{7,10}$ and $VV_{10,1}$. In fact, the latter cuts are obtained from the original $V_{1,4}$ cut (blue dashed line) by a rotation by $\frac{2\pi}{4}$. This gives a square, related to divisor 4.
- $VV_{1,5}$ cut (blue solid line), is repeated as $VV_{5,9}$ and as $VV_{9,1}$. In fact, the latter cuts are obtained from the original $V_{1,5}$ cut (blue solid line) by a rotation by $\frac{2\pi}{3}$. This gives a triangle, related to divisor 3.
- $VV_{1,3}$ cut will give a hexagon. If this is cut is repeated 6 times over an angle of $\frac{2\pi}{6}$, related to divisor 6.
- $VV_{1,2}$ cut does not divide the polygon, but coincides with the original side 1. Rotation this by $\frac{2\pi}{12}$ will give the original dodecagon. An inscribed dodecagon is obtained via a $SS_{1,2}$ cut, rotated twelve times by $\frac{2\pi}{12}$. Both cases are related to $d_{max} = 12$.

In $m = 4$ the upper blue line is a $SS_{1,2}$ cut from side 1 to side 2. When rotated by 180° the lower blue line is obtained, for the second smallest divisor $d_2 = 2$.

The inscribed figures (triangle and square in $m = 12$) close in one rotation but this can be generalized, for $m = \frac{p}{q}$ with *p, q* rational numbers and relative prime. .

- In $m = 5$ the sequence $VV_{1,3}, VV_{3;5}, VV_{5,2}, VV_{2,4}, VV_{4,1}$ will generate a pentagram in the pentagon, i.e. a figure that closes in 2 rotations, having 5 angles that are spaced $\frac{4\pi}{5} = 144°$ apart. This generates the classic Pythagorean pentagram that led to the discovery of irrational numbers and the golden ratio, corresponding to a *RGC* with $m = 5/2$.
- In $m = 7$ the sequence $VV_{1,3}, VV_{3,5}, VV_{5,7}, VV_{7,2}$ , $VV_{2,4}, VV_{4,6}, VV_{6,1}$ generates a heptagram, closing in 2 rotations, corresponding to a *RGC* with $m = \frac{7}{2}$
- Also in $m = 7$ the sequence $VV_{1,4}, VV_{4,7}, VV_{7,3}, VV_{3,6}, VV_{6,2}, VV_{2,5}, VV_{5,1}$ generates a heptagram, closing in 3 rotations, corresponding to a *RGC* with $m = \frac{7}{3}$.

## 4.3 Rotations and scaling

For $m = 16$ results of cutting with $d_m$ knives for divisors 2, 4, 8 and 16 are shown in Figure 37. They are rotations of the $d_1$ knife.

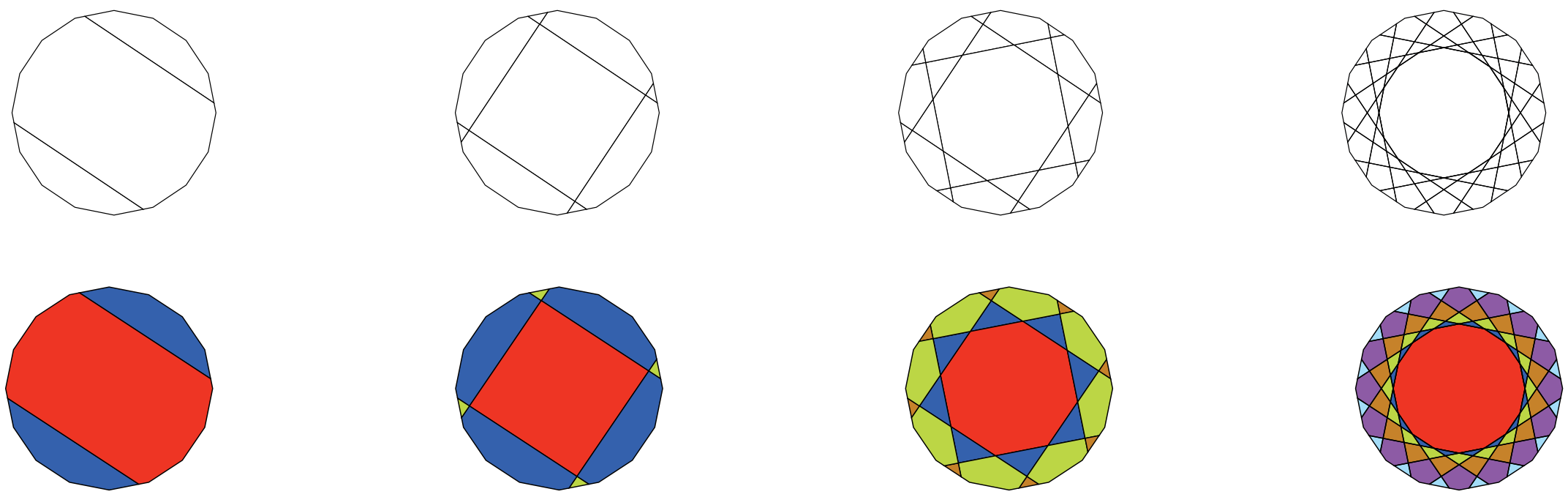

*Figure 37: Cutting with $d_i$ knives ($d_2$ =2, $d_3$ =4, $d_4$ =8, $d_5$ = 16) and the resulting sectors*

The relation of divisors and rotations show that *VV* and *SS* cuts can be transformed into each other by rotations:

- In $m = 6$, the $VV_{max} = VV_{1,4}$ cut or diagonal can be rotated every 60° and all diagonals meet in the centre (Result is 6 diagonals that coincide 2 by 2)
- This shape can then be rotated by 30°, resulting in $SS_{max}$. The rotation can in fact be done for any angle.

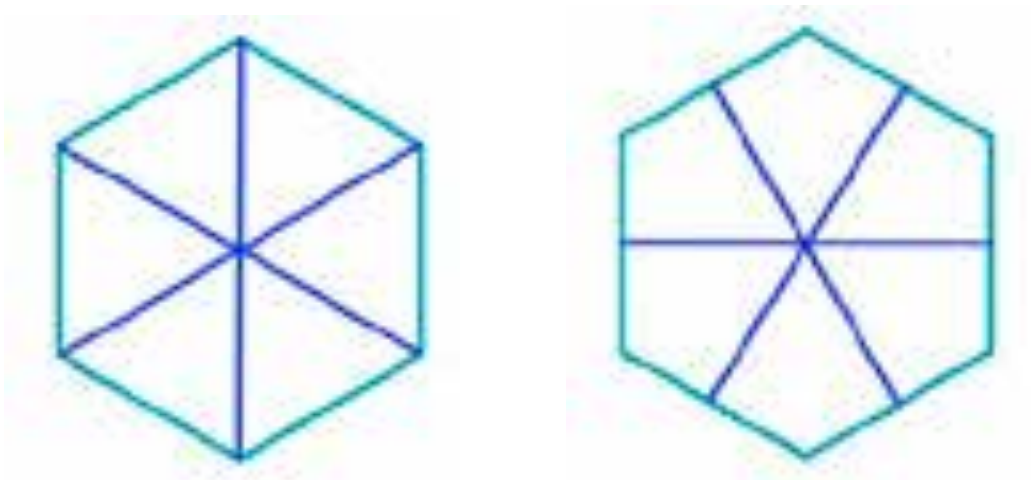

*Figure 38: Converting $VV_{max}$ into SS cut through centre of the sides by rotation*

There are many ways in which the figures can be transformed into the other figures, using rotation and scaling. A first example sequence of rotations and scaling could be:

- If the red square in Figure 39**f** is rotated by 45° the result is **c,**
- If the red square in **c** is scaled to size zero, d results.
- If the cross in **d** is rotated by 45° the result is **a.**

A second example in Figure 39:

- If the red square in **f** is scaled to a larger size, **e** results; when it is scaled to **a** smaller size **g** results.
- When the inscribed figure in **g** is rotated so that one of the sides of the small yellow triangles ends in a vertex, we obtain **b**.

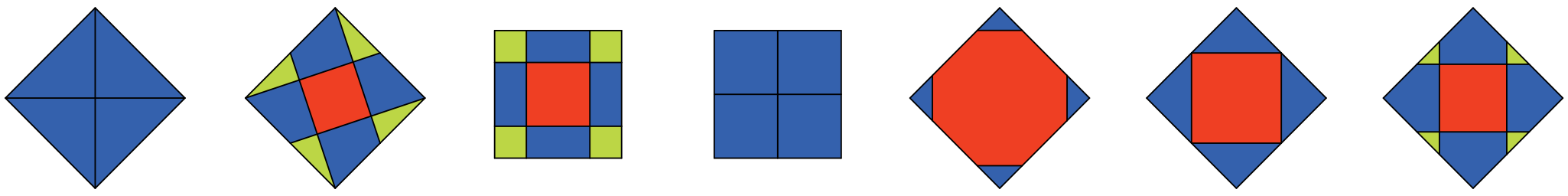

**a** **b** **c** **d** **e** **f** **g**

*Figure 39: Ways of cutting square*

This shows that all figures can be considered as transformations of an inscribed square (and in general the inscribed $m$-polygon inside an $m$-polygon). The analytic representation of the $d_i$ -knife was given in section 2.7 (Equations 12 and 13) [8].

## 4.4 The way of cutting determines the final result

In Figure 39 from left to right, all possible cuts are shown for a square

- **a**. One case of *VV* cut, with two $VV_{max}$ diagonals
- **b.** One case of *VS* cut, with a $VS_{1,2}$ cut and its rotations over $\frac{2\pi}{4}$. It total 3 different shapes are created. Four green triangles, four blue quadrilaterals and one red square. The smaller red square can be considered as the inscribed square rotated and scaled to smaller size.
- **c** and **d**. Two cases of $SS_{max} = SS_{1,3}$ cuts and rotations. One *SS* cut does not pass through the centre (**c.**) and the other one passes through the centre (**d**). The former creates 3 different sets of quadrilaterals indicated by different colours. The latter creates four different squares. (In *GML* in **d.** these four different squares form one body, and in **c.** each of the coloured zones creates 3 separated bodies)
- Three cases of $SS_{1,2}$ cuts (**e**, **f**, **g**). It is clear that the result depends on where the cut is made. The middle figure **f.** is the inscribed square, while **e.** and **g.** are scaled version (larger and smaller, without rotation).
- This also generates different shapes. In **e.** four triangles and one octagon; in **f.** four triangles and one square and in **g.** one set of 4 triangles, one set of 4 pentagons and a central square. Again in *GML* and in rational Gielis curves RGC they will form different bodies or layers.

The same logic is applied to hexagons in Figure 40, with more possibilities of cutting, namely two *VV* cuts, two *VS* cuts and eight *SS* cuts.

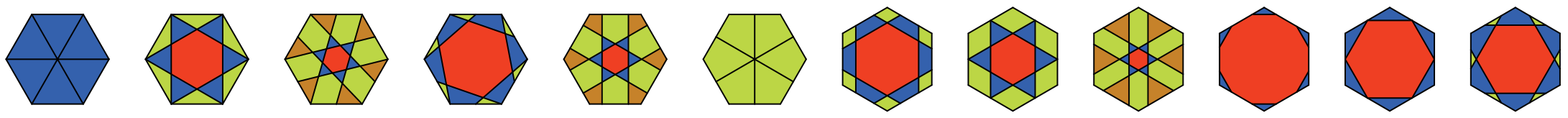

*Figure 40: Hexagon cutting*

## 4.5 Inheritance of possible cuts for different divisors

In figure 41 *VV* and *VS* cuts are shown for $m = 4, 6$. A square has three divisors, so 3 rows; for the hexagon there are four rows corresponding to 4 divisors. For $m$ prime only 2 rows result, for divisors $m$ and 1.

- For even $m$ the vertical columns of square and hexagon in Figure 41 show that the possible cuts are the same for all divisors, as the result of rotating cut 1 over the relevant angle, related to divisor. The identification of vertices and of the knife (e.g. for $d_2$ =2 lower row Figure 41) links the cutting of m-polygons to cutting of $GML_m^n$ bodies
- For odd $m$ (in case of pentagon in Figure 42), the *VV* cuts are inherited from cut 1 via rotations. This is not the case for the *VS* cuts of the pentagon in Figure 42. However, one figure is missing in the upper row, namely the *VS* cut not through the centre. If this figure is also considered, then also for the pentagon the number of cuts is fully inherited for the two divisors, as in the case of the even *m.*

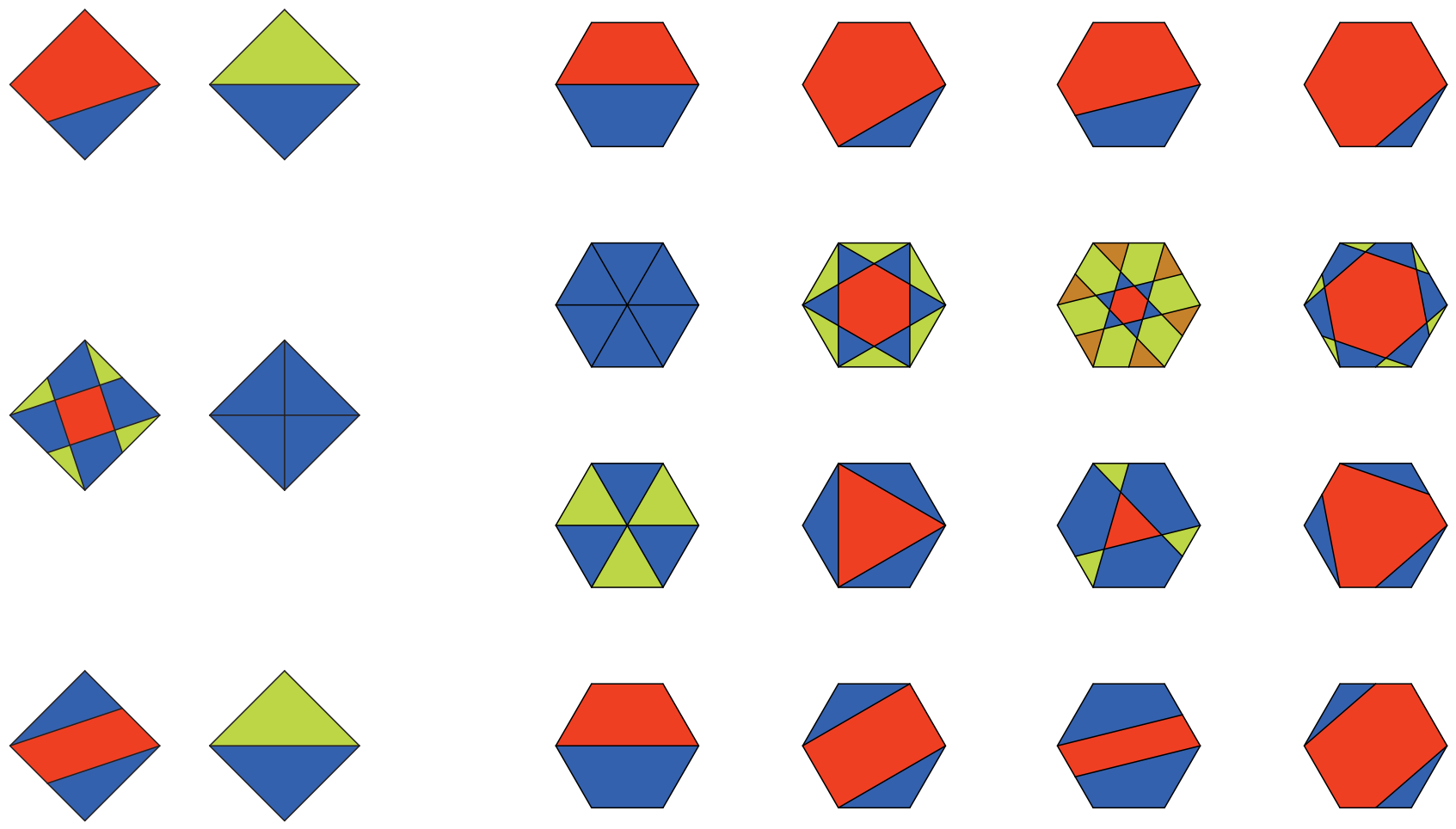

*Figure 41: VV and VS cuts of square and hexagon*

- The reason why in Figure 42 for the pentagon (and in general for any odd $m$) only one figure is shown for *VS* (the cut through the centre) is that the number of sectors and the number of vertices and sides of the resulting polygons remains the same. In the case of the pentagon and $VS_{1,3} = VS_{max}$, two quadrilateral figures are created, whether or not the cut goes through the centre. The quadrilateral shapes share the topological characteristic of four vertices and four sides.
- In case of the $VS_{1,3} = VS_{max}$ cut going through the centre, the two quadrilateral shapes have also the exactly same shape, a geometrical characteristic. If $VS_{1,3} \neq VS_{max}$, i.e. not going through the centre, then the geometrical shapes of the two quadrilaterals is different.

So, if this **geometrical** characteristic is considered, also for divisor 1 or one cut, two different shapes need to be considered, whereby the two shapes are different or the same. In the **topological** case, these two shapes reduce to one shape as in Figure 42. The same line of reasoning can be considered for *SS* cuts, where in Figure 42 there are 2 line cuts in upper row versus 8 in the lower row.

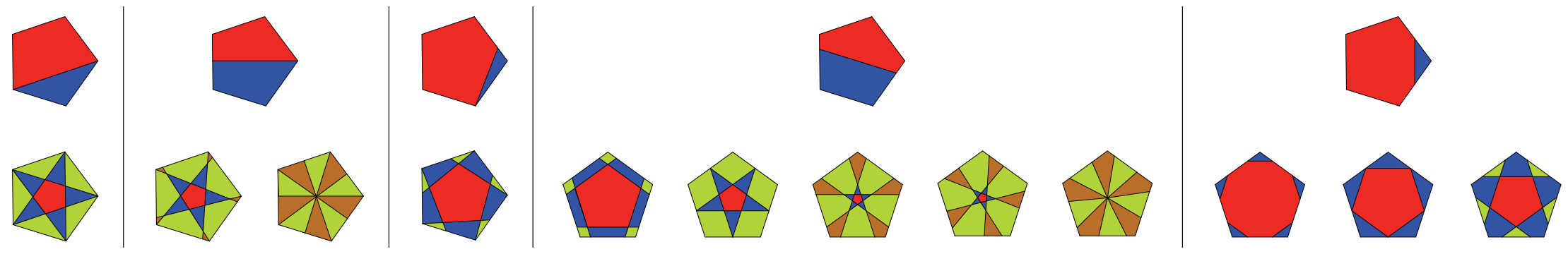

*Figure 42: All possible cuttings of a pentagon via VV, VS and SS cuts*

In general, starting from the maximum divisor, or max cuts equal to $m$, for other divisors the number of possibilities is inherited precisely, in the **geometrical** sense.

## 4.6 $d$-knives and projective geometry

The analytic definition of $d$-knife is a construction, with $m$ straight lines ([8], section 2.8)

$$\sin\left(\alpha + \frac{2\pi}{m} i\right) x_i + \cos\left(\alpha + \frac{2\pi}{m} i\right) y_i + \delta = 0, \quad i = 0, 1, .., m-1; \quad -\frac{\pi}{m} \leq \alpha \leq \frac{\pi}{m},$$

The precise cutting with a $d_i$ −knife can also be considered as a projection of a shape with symmetry $m$ or one of the divisors of $m$, onto the basic shape of the cross section of the *GML* with symmetry $m$, along a cone with the same cross section (Figure 43). Cutting with a $d_m$ −knife is equivalent to projecting an $m$-regular polygon onto an $m$-regular polygon with zooming parameter $\delta$. In Figure 43 the projection of a red square onto a green square (the basic one) is shown. At the top of the pyramid, the red square reduces to one point and the blue zone narrows to lines, as in Figures 39**c** and 39**d**.

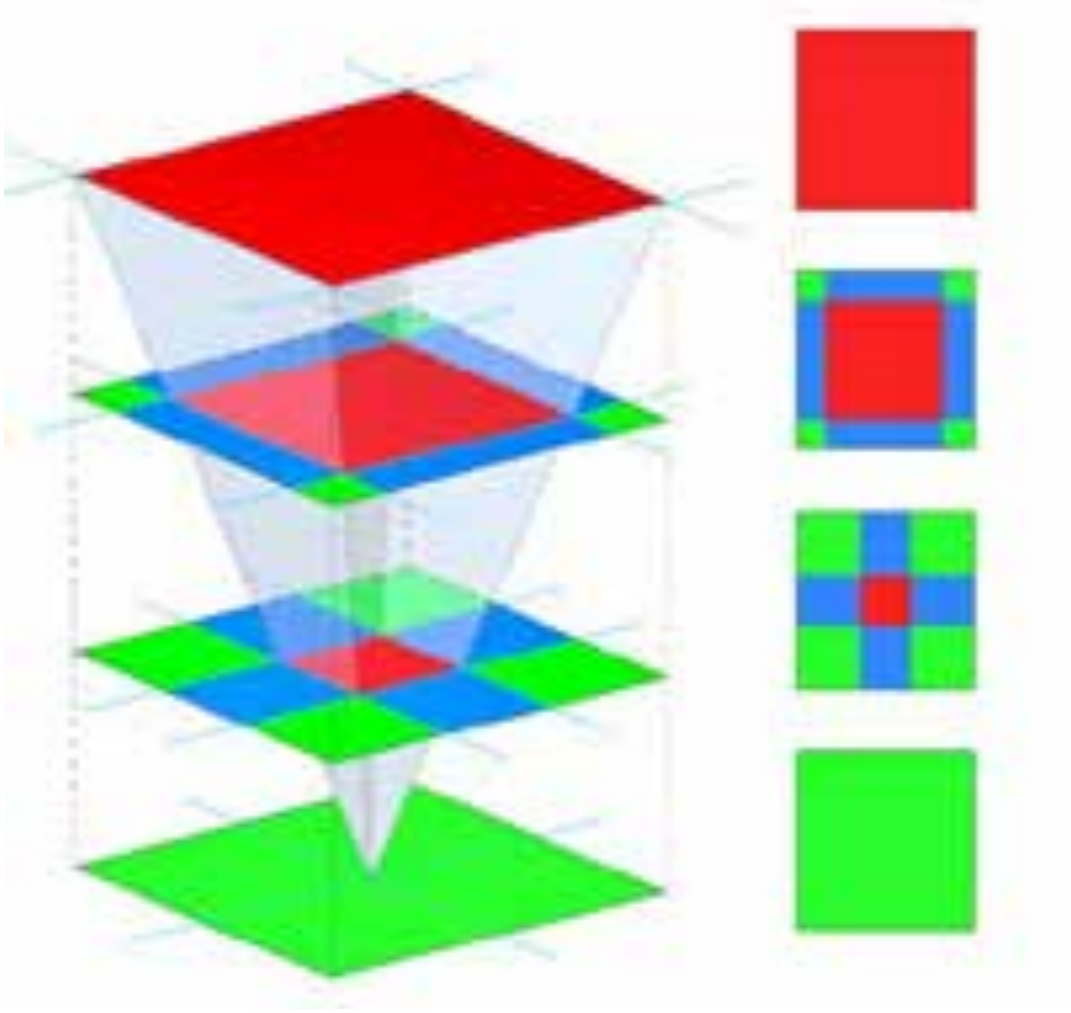

*Figure 43: Projective equivalence with cutting with $d_i$knife for squares $d_m = 4$*

For the smallest divisor 1 the cut is with one knife, and this cut can be repeated by rotation as in Figure 37 for other divisors. For even numbers the second smallest divisor $d_2 = 2$ (a rectangle is projected which at the top of the pyramid also reduces to a line (compare Figure 41 lower row), and so on. This projection corresponds to the dilation parameter $\delta$. Rotation parameter $\alpha$ corresponds to a rotation of the prismatic pyramid relative to the base polygon.

## 4.7 The number of independent objects after cutting

When cutting $m$-regular polygons, $GML$ bodies (or more generally Generalized Rotating and Twisting $GRT_m^n$ surfaces and bodies, whereby the basic line is a Gielis curve with $m$ rational), this leads to a number of independent sectors in the plane cross-sections. Tables 9 and 10 give for even and odd values of $m$ all possibilities:

- $d_i\,, i = 1,2, \ldots, m$ are different divisors of the number $m$. $N$ is the number of different nontrivial ( $d_i \neq 1, m$) divisors of the number $m$.
- if $m = 2k + 1$ then all divisors $d_i\,, i = 1,2, \ldots, N$ of $m$ are odd numbers.
- if $m = 2k$ some of its divisors $d_i\,, i = 1,2, \ldots, N$ may be odd.

*Table 9: Example cuts of $GRT_m^n$ bodies for odd $m = 3, 5, 9, 27$*

| m= 2·1+1 | d$_i$-knives | SS -cuts | SS$_C$ -cuts | VS -cuts | VS$_C$ -cuts | VV- cuts | VV$_C$ -cuts | |
|---|---|---|---|---|---|---|---|---|
| $d_m$ = 3 | 3-knives | 4 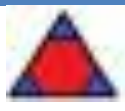 7 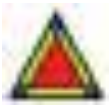 | 6 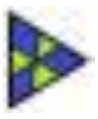 | 7 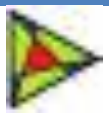 | 6 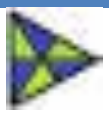 | 0 | 0 | k=1 |
| $d_1$ = 1 | 1-knife | 2 | 2 | 2 | 2 | 0 | 0 | |

| m= 2·2+1 | d$_i$-knives | SS | SS$_C$ | VS | VS$_C$ | VV | VV$_C$ | |
|---|---|---|---|---|---|---|---|---|
| $d_m$= 5 | 5-knives | 6 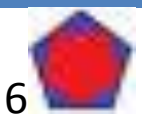 11 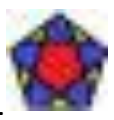 16 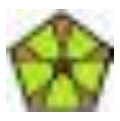 | 10 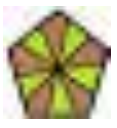 | 11 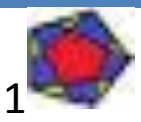 16 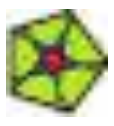 | 10 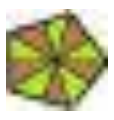 | 11 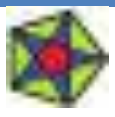 | 0 | k=2 |
| $d_1$ = 1 | 1-knife | 2 | 2 | 2 | 2 | 2 | 0 | |

| m= 2·4+1 | d$_i$-knives | SS | SS$_C$ | VS | VS$_C$ | VV | VV$_C$ | |
|---|---|---|---|---|---|---|---|---|
| $d_m$= 9 | 9-knives | 10 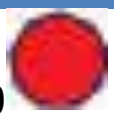 19 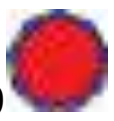 28 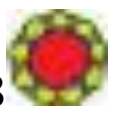 37 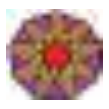 46 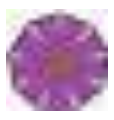 | 18 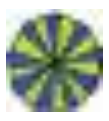 | 19 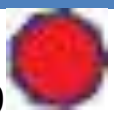, 28 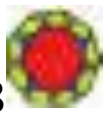 37 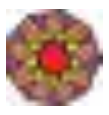 46 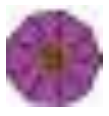 | 10 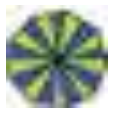 | 19 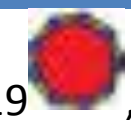, 28 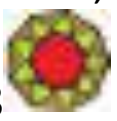 37 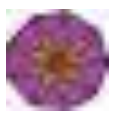 | 0 | k=4 |
| $\boldsymbol{d_2}$ = **2·1+1** | 3-knives | 4 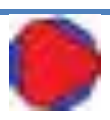 7 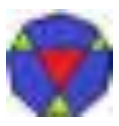 | 6 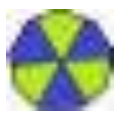 | 4 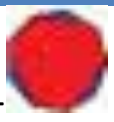 7 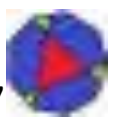 | 6 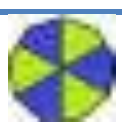 | 4 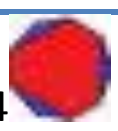 7 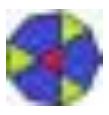 | 0 | $k_1$=1 |
| $\boldsymbol{d_1}$ = **1** | 1-knife | 2 | 2 | 2 | 2 | 2 | 0 | |

| m= 2·13+1 | $d_i$-knives | SS | $SS_C$ | VS | $VS_C$ | VV | $VV_C$ | |
|---|---|---|---|---|---|---|---|---|
| $\boldsymbol{d_m}$= **27** | 27-knives | 28, 55, 82, 109, 136,163,190,217 , 244, 271, 298 , 325, 351, 379 | 54 | 55,82,109,13 6,163,190, 217, 244,271,298 , 325,351,379 | 54 | 55,82,109, 136,163,190,2 17, 244,271,298 , 325, 351 | 0 | k=13 |
| $\boldsymbol{d_3}$ **= 9 2·4+1** | 9-knives | 10,19,28,37,46 | 18 | 10,19,28,37,4 6 | 18 | 10,19,28,37,46 | 0 | $k_1$=4 |
| $\boldsymbol{d_2}$ **= 3 2·1+1** | 3-knives | 4,7 | 6 | 4,7 | 6 | 4,7 | 0 | $k_2$=1 |
| $\boldsymbol{d_1}$ **= 1** | 1-knife | 2 | 2 | 2 | 2 | 2 except m=3 | 0 | |

**Remark**: For odd $m$ a complete correspondence is found with Rational Gielis Curves for $VV$-cuts. For example for $m = 5$ (Section 3.5 and Table 6): if $m = \frac{p}{q} = \frac{5}{2}, \frac{5}{3}, \frac{5}{4}$ the number of independent shapes are 6, 11 and 16 (for $p = 5$, resp. $p + 1, 2p + 1$ and $3p + 1$).

*Table 10: Example cuts of $GRT_m^n$ bodies for even $m = 4, 6, 8, 10, 24$*

| m= 2·2 | $d_i$-knives | SS | $SS_C$ | VS | $VS_C$ | VV | $VV_C$ | |
|---|---|---|---|---|---|---|---|---|
| $\boldsymbol{d_m}$= **2·2** | 4-knives | 5 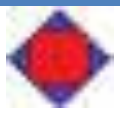 9 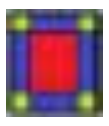 | 0 | 9 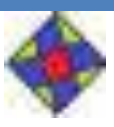 | 0 | 0 | 0 | m=2k |
| $\boldsymbol{d_2}$ = **2·1** | 2-knives | 3 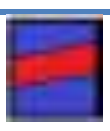 | 4* 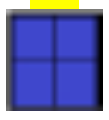 | 3 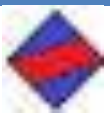 | 0 | 0 | 4* 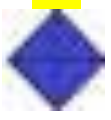 | $d_1$ is even and $k_1$=1 |
| $\boldsymbol{d_1}$ = **1** | 1-knives | 2 | 2 | 2 | 2 | 2 | 2 | |

| m= 2·3 | $d_i$-knives | SS | $SS_C$ | VS | $VS_C$ | VV | $VV_C$ | |
|---|---|---|---|---|---|---|---|---|
| $\boldsymbol{d_m}$= **2·3** | 6-knives | 7 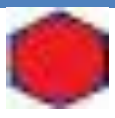 13 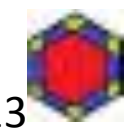 19 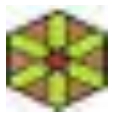 | 0 | 13 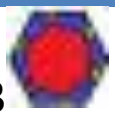 19 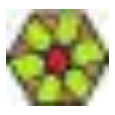 | 0 | 13 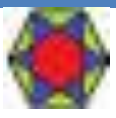 | 0 | k=3 |
| $\boldsymbol{d_3}$= **m/2 = 2·1+ 1** | 3-knives | 4 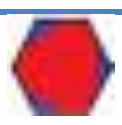 7 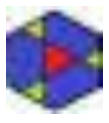 | 6 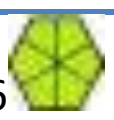 | 4 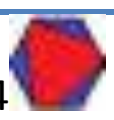 7 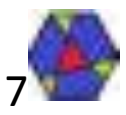 | 0 | 4 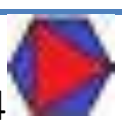 | 6 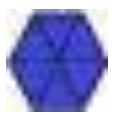 | $d_1$ is odd $k_1$=1 |
| $\boldsymbol{d_2}$ = **2·1** | 2-knives | 3 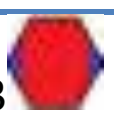 | 2* | 3 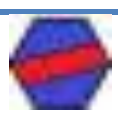 | 0 | 3 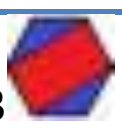 | 2* | $d_2$ is even $k_2$=1 |
| $\boldsymbol{d_1}$ = **1** | 1-knives | 2 | 2 | 2 | 2 | 2 | 2 | |

| m= 4·2 | $d_i$-knives | SS | $SS_C$ | VS | $VS_C$ | VV | $VV_C$ | |
|---|---|---|---|---|---|---|---|---|
| $d_m$= 4·1 | 8-knives | 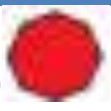 9 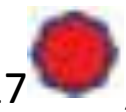 17 , 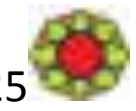 25 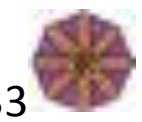 33 | 0 | 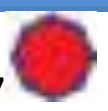 17 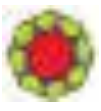 25 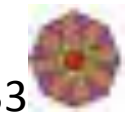 33 | 0 | 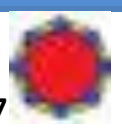 17 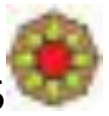 25 | 0 | k=4 |
| $d_3$=2 ·2 | 4-knives | 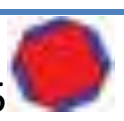 5 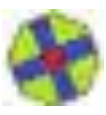 9 | 8 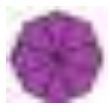 | 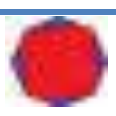 5 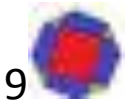 9 | 0 | 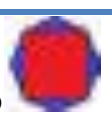 5 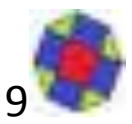 9 | 8 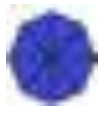 | $d_1$ is even $k_1$=2 |
| $d_2$ = 2·1 | 2-knives | 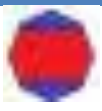 3 | 4* 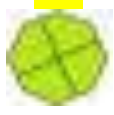 | 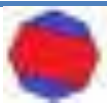 3 | 0 | 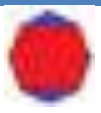 3 | 4* 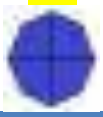 | $d_2$ is even $k_2$=1 |
| $d_1$ = 1 | $1_i$-knives | 2 | 2 | 2 | 2 | 2 | 2 | |

| m= 2·5 | $d_i$-knives | SS | $SS_C$ | VS | $VS_C$ | VV | $VV_C$ | |
|---|---|---|---|---|---|---|---|---|
| $d_m$= 2·5 | 10-knives | j·10+1 j=1,2,3,4,5 (11,21,31,41,51) | 0 | j·10+1 j=2,3,4,5 (21,31,41,51) | 0 | j·10+1 j=2,3,4 (21,31,41) | 0 | k=5 |
| $d_3$=2 ·2+1 | 5-knives | j·5+1 j=1,2,3 (6,11,16) | 10 | j·5+1 j=1,2,3 (6,11,16) | 0 | j·5+1 j=1,2 (6,11) | 10 | $k_1$=2 |
| $d_2$ = 2·1 | $2_i$-knives | j· $2_i$+1 j=1 (3) | 4* | j· $2_i$+1 j=1 (3) | 0 | j· $2_i$+1 j=1 (3) | 4* | $d_2$ is even $k_2$=1 |
| $d_1$ = 1 | $1_i$-knives | 2 | 2 | 2 | 2 | 2 | 2 | |

| m= 2·12 | $d_i$-knives | SS | $SS_C$ | VS | $VS_C$ | VV | $VV_C$ | |
|---|---|---|---|---|---|---|---|---|
| $d_m$= 2·12 | 24-knives | j·24+1 where j=1,2,…,12 | 0 | j·24+1 where j=2,3,…,12 | 0 | j·10+1 where j=2,3,…,11 | 0 | k=12 |
| $d_7$=2 ·6 | 12-knives | j·12+1 where j=1,2,…,6 | 2·12 | j·12+1 where j=1,2,…,6 | 0 | j·12+1 where j=1,2,…,6 | 2·12 | $k_1$=6 |
| $d_6$=2 ·4 | 8-knives | j· 8+1 where j=1,2,3,4 | 2·4* | j· 8+1 where j=1,2,3,4 | 0 | j· 8+1 where j=1,2,3,4 | 2·4* | $k_2$=4 |
| $d_5$=2 ·3 | 6-knives | j· 6+1 where j=1,2,3 | 2·3* | j· 6+1 where j=1,2,3 | 0 | j· 6+1 where j=1,2,3 | 2·3* | $k_3$=3 |
| $d_4$=2 ·2 | 4-knives | j· 4+1 where j=1,2 | 2·4* | j· 4+1 where j=1,2 | 0 | j· 4+1 where j=1,2 | 2·4* | $k_4$=2 |
| $d_3$=2 ·1+1 | 3-knives | j· 3+1 where j=1,2 | 2·3 | j· 3+1 where j=1,2 | 0 | j· 3+1 where j=1,2 | 2·3 | $k_5$=2 |
| $d_2$ = 2·1 | 2-knives | j· 2+1 where j=1 | 4* | j· 2+1 where j=1 | 0 | j· 2+1 where j=1 | 4* | $k_2$=1 |
| $d_1$ = 1 | $1_i$-knives | 2 | 2 | 2 | 2 | 2 | 2 | |

## 4.8 Using radial knives

The Möbius phenomenon only occurs for *m* even, when the knife cuts through the centre, since only then congruent figures are obtained through cutting that can be connected throughout with a twist. Going through the centre is an important condition, but whether the cuts are VV, VS or SS is less important: Möbius phenomena can occur in all cases. This Möbius phenomenon never occurs for *m*=odd; while then congruent figures may be obtained, they only have mirror symmetry, not rotational.

The main condition, congruent figures and rotational symmetry, is achieved when using radial knives (see section 2.10), emanating from the centre *C* of the figure, to the perimeter, so *CV* or *CS* cuts result with the length of the knife half compared to the previous knives. In this case the Möbius phenomenon can be achieved in all cases, irrespective whether *m* is odd or even. The word knife can be substituted by other words, like secants, diagonals, rays and chords, but has the advantage that the separated zones can actually be separated.

# 5 Demonstratio

## 5.1 The geometrical solution

The total number of ways of cutting an $m$-polygon according to the rules described above with $d$-knives for the geometrical case, for $m$ = even and $m$= odd respectively is as follows: *VV, VS* and *SS* cuts increase from by 1, 1 and 3 respectively from a given even or odd number to the next even or odd number. As a result, the total number of ways of cutting using a $d_m$ knife increases by 5 to each subsequent even or odd number (Table 11, Subtotal), which gives the sequence $(2), 2,\ 7, 7, 12, 12, 22, 22, 27, 27, 32, 32, 37, 37$ … for $m\ =\ 2, 3, 4, 5 \ldots 15$. Taking sums the sequence 4, 14, 24, 34, 44, 54, 64, 74…. results, which is monotonically increasing.

*Table 11: Number of possible cuts for even and odd m*

| m = even | Cut type | 2 | 4 | 6 | 8 | 10 | 12 | 14 |
|---|---|---|---|---|---|---|---|---|
| **(m-2)/2** | VV | 0 | 1 | 2 | 3 | 4 | 5 | 6 |
| **(m-2)/2** | VS | 0 | 1 | 2 | 3 | 4 | 5 | 6 |
| **Step +3** | SS | 2 | 5 | 8 | 11 | 14 | 17 | 20 |
| | Subtotal | 2 | 7 | 12 | 17 | 22 | 27 | 32 |
| | Divisors | 2 | 3 | 4 | 4 | 4 | 6 | 4 |
| **TOTAL** | | **4** | **21** | **48** | **68** | **88** | **162** | **128** |

| m = odd | Cut type | 3 | 5 | 7 | 9 | 11 | 13 | 15 |
|---|---|---|---|---|---|---|---|---|
| **(m-3)/2** | VV | 0 | 1 | 2 | 3 | 4 | 5 | 6 |
| **(m+1)/2** | VS | 2 | 3 | 4 | 5 | 6 | 7 | 8 |
| **Step + 3** | SS | 5 | 8 | 11 | 14 | 17 | 20 | 23 |
| | Subtotal | 7 | 12 | 17 | 22 | 27 | 32 | 37 |
| | Divisors | 2 | 2 | 2 | 3 | 2 | 2 | 4 |
| **TOTAL** | | **14** | **24** | **34** | **66** | **54** | **64** | **148** |

Since the number of possibilities is determined by $d_1$ and $d_m$ and is inherited by the other divisors given identification of vertices and knives (Figure 44 for divisors 1, 2 and 6 in a hexagon), the total number is then the subtotal times the number of divisors. The identification links planar geometry to 3D $GML_m^n$ bodies. If one follows the $d_1$-knife along the basic line of the $GML_6^6$ body, the different positions of the knives indicated by arrows in Figure 44 for $d_1$ function as a clock, relative to the torus circumscribing the $GML_6^6$ body. For other knives the clock arithmetic is the same, albeit with more hands.

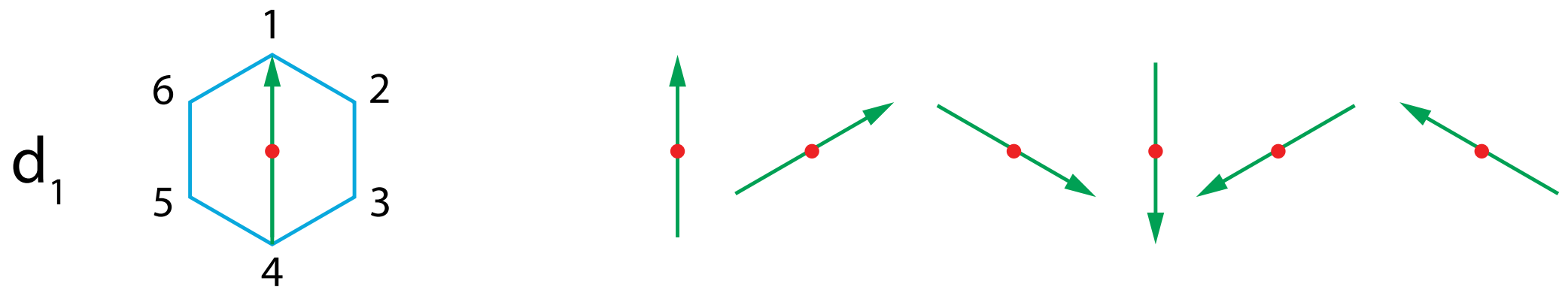

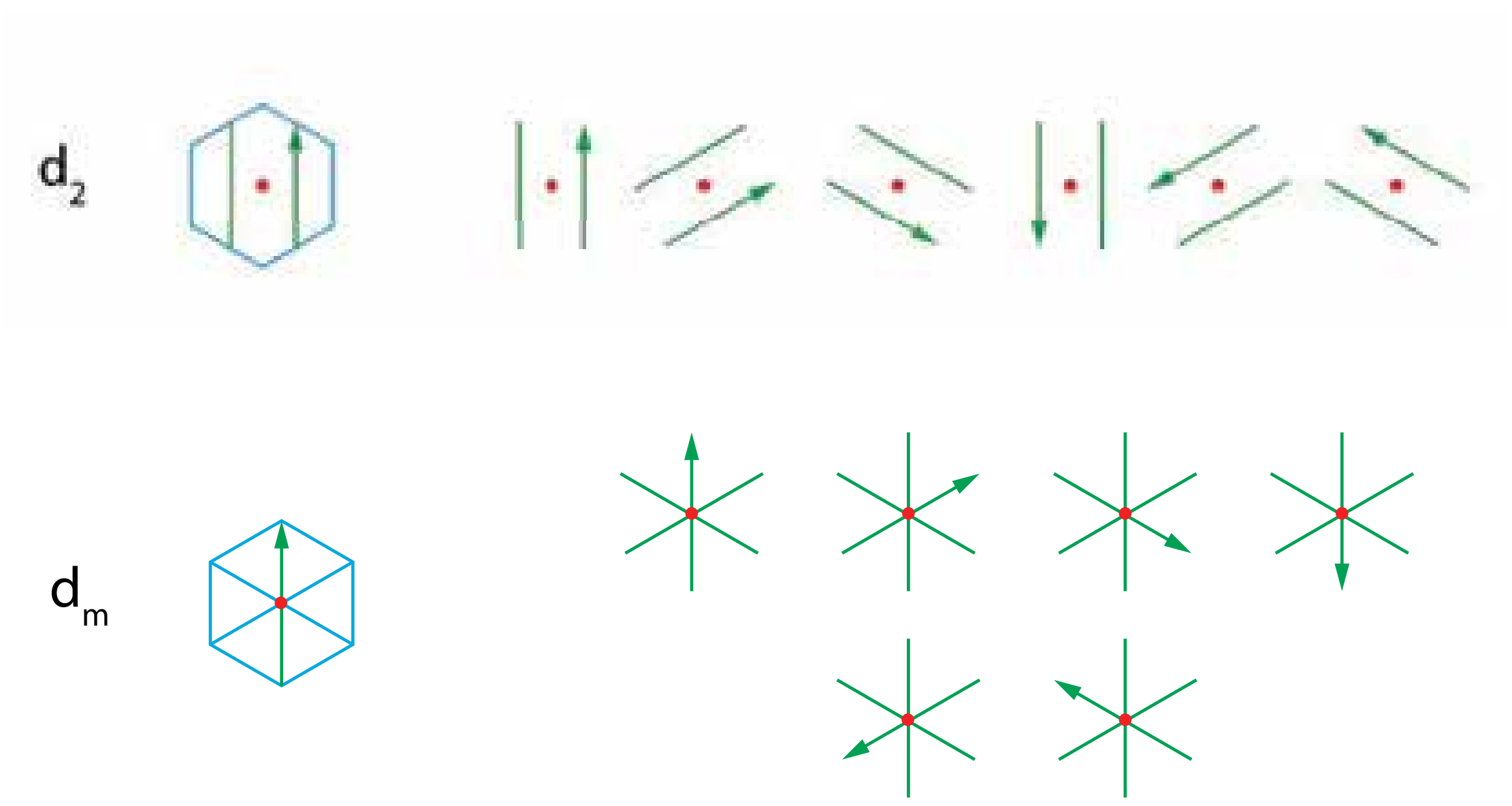

*Figure 44: Inheritance from $d_1$ by all d-knives through identification.*

Table 12 gives the product of the subtotals with positive integers. Entries in rows can be computed as $u_n = u_{n-1} + u_{n-2} - u_{n-3}$, with $u_0 = 2$ for even and $u_0 = 3$ for odd numbers. In bold red are the totals of Table 11 for $m$ odd and in green for $m$ even.

*Table 12: Product of subtotals of Table 11 and number of divisors.*

| N°Divisors | 2 | 7 | 12 | 17 | 22 | 27 | 32 | 37 |
|---|---|---|---|---|---|---|---|---|
| **2** | **4** | **14** | **24** | **34** | 44 | **54** | **64** | 74 |
| **3** | 6 | **21** | 36 | 51 | **66** | 81 | 96 | 111 |
| **4** | 8 | 28 | **48** | **68** | **88** | 108 | **128** | **148** |
| **5** | 10 | 35 | 60 | 85 | 110 | 135 | 160 | 185 |
| **6** | 12 | 42 | 72 | 102 | 132 | **162** | 192 | 222 |

The number of possible cuts can be given by a recurrence formula. For $N_m$, the number of ways of cutting an $m$-gon for one divisor, whereby $N_m^{SS}$ stands for the number of $SS$ cuts for $m$ and $N_{m-2}^{SS}$ for the number of $SS$ cuts for the polygon with $(m-2)$ (i.e. the previous odd or even number) and with $k$ a natural number.

- For even $m$ $(= 2k)$: $N_{m=2k} = m + 1 + N_{m-2}^{SS}$ (1a)
- For odd $m$ $(= 2k + 1)$: $N_{m=2k+1} = m + 2 + N_{m-2}^{SS}$ (1b)

If the number of SS cuts is kept separate, taking into account the step +3, this part of the general formula is recursive. Since $N_m^S = \left(N_{m-2}^S + 3\right)$, it follows that $N_{m-4}^{SS} + 6 = N_{m-6}^{SS} + 9 = \ldots$

Because of the exact inheritance for the geometrical case the total number of ways of cutting for all divisors is then the above formula times the number of divisors $N_m^{div}$ of a number *m.*

**Theorem 1** *(the geometrical solution):*

*The total number of different ways of cutting an m-polygon* $\Xi_m^{geo}$ *is the number of 1 or m cuts, times the number of divisors of m.*

- For even $m$ $(= 2k)$: $\Xi_m^{geo} = N_m^{div}(m + 1 + N_{m-2}^S)$
- For odd $m$ $(= 2k + 1)$: $\Xi_m^{geo} = N_m^{div}(m + 2 + N_{m-2}^S)$

When considering polygons with convex sides, $VV_{i,i+1}$ are possible, so the number of cuts increases by $m$. When knives are used to cut a circle from equally spaced points, the $VV_{i,i+1}$ cuts need to be added, in particular $m$ cuts.

- For even $m$ $(= 2k)$: $N_{m=2k} = 2m + 1 + N_{m-2}^{SS}$ (2a)
- For odd $m$ $(= 2k + 1)$: $N_{m=2k+1} = 2m + 2 + N_{m-2}^{SS}$ (2b)

## 5.2 The topological solution

The number of ways of cutting an $m$-polygon according to the rules for the topological solutions is given in Table 13 for *for* $m = 3,4, \dots ,16$. Annex 1 graphically shows all possibilities for *for* $m = 6,7,8,9,10$. From this follows Theorem 2 (the topological solution).

*Table 13: All possible cuttings for m-regular polygons for* $m = 3,4, \dots ,16$

| *m* | *SS* – cuts | *VS* cuts | *VV*-cuts | Total |
|---|---|---|---|---|
| 3 | 1+5 = 6 | 3 | 0 | **9** |
| 4 | 1+2+1+3 +2+1= 3 +4 +3= 10 | 3 | 3 | **16** |
| 5 | 1+5+1+3 = 6 +4 = 10 | 5 | 2 | **17** |
| 6 | 1+2+1+3 +1+3 +2+3+1+2+1+1= 11+6 +4= 21 | 8 | 8 | **37** |
| 7 | 1+5+1+3 +1+3 = 6 +4+4 = 14 | 7 | 4 | **25** |
| 8 | 1+2+1+3 +1+3 +1+3 +2+ 1+ 3+1+2+1+1+1= 15+7 +5= 27 | 12 | 12 | **51** |
| 9 | 1+5+1+3 +1+3 +1+3 +3+3+1+1 = 6 +4+4 +4 +8= 26 | 14 | 9 | **49** |
| 10 | 1+2+1+3 +1+3 +1+3 +1+3 +2+3+1+3+1+2+1+1+1+1=19 +10+6=35 | 16 | 16 | **67** |
| 11 | 1+5+1+3 +1+3+1+3 +1+3 = 6 +4+4 +4+4 = 22 | 11 | 8 | **41** |
| 12 | 59 | 30 | 30 | **119** |
| 13 | 1+5+1+3 +1+3+1+3 +1+3+1+3 = 6 +4+4 +4+4 +4= 26 | 13 | 10 | **49** |
| 14 | 49 | 24 | 24 | **97** |
| 15 | 54 | 31 | 24 | **109** |

**Theorem 2 (the topological solution)**

If $m = 2k + 1$ and has $N$ nontrivial divisors $d_2, d_3 \ldots d_{N+1}$ and $d_1 \equiv 1, d_{N+2} \equiv d_m \equiv m$, then the number of all possible variants of cutting of $GML_m^n$ bodies is

$$\Xi_m^{top} = 8k + 1 + 3Nk + \sum_{i=2}^{N+1} \left[\frac{k}{d_i}\right] + 2N$$

If $m = 2k$ and has $N$ nontrivial divisors $d_2, d_3 \ldots d_{N+1}$ and $d_1 \equiv 1, d_{N+2} \equiv d_m \equiv m$, then the number of all possible variants of cutting of $GML_m^n$ bodies is

$$\Xi_m^{top} = 8k - 5 + 3Nk + \sum_{i=2}^{N+1} \left[\frac{k-1}{d_i}\right] - N$$

This is in a slightly different form compared to Theorem 1: $\sum_{all\ div}(N_{m=2k}) = \Xi_m^{top}$ depending on whether the total number of variants is expressed in terms of total number of divisors or total number of non-trivial divisors $N$ (excluding $d_1$ and $d_m$).

The complete proof of Theorem 2 will be given in a separate paper. It is based on the fundamental facts from the theory of cyclic groups with a finite number of elements ($m$);

1. The number of cyclic subgroups is the number ($N$) of nontrivial divisors of $m$.
2. The number of elements in each subgroup is the number of transactions and equal to the $gcd(m, i)$
3. The number of cuts is either 3 or 1, (3 or 1 mod8) and this is determined by the property of the subgroup and the property of the cut line – i.e. when the number of cuts is three, then the ends of the survey line lie on the same strings of the initial polygon, except for the case when $k = \left[{}^{m}/_{2}\right] + 1$
4. If $k = \left[{}^{m}/_{2}\right] + 1$ and for an odd number $m$ the number of cuts is 5, and for even $m$ the number of cuts is 2. This is determined by the property of the subgroup and the property of the cut line. In the latter case $m = 2k$ an important role is played by the rotational symmetry.

## 5.3 Number of independent bodies or zones

**Theorem 3:** The number of independent objects depends on 1) the type of cut, 2) the divisor or $d_i$-knife, and 3) the parity of $m$ and its divisors. For odd or even $m$ the number of possible independent objects appearing after different cuts of $m$-polygons and $GRT_m^n$ bodies with $d_i$ –knives is given in Table 14 for odd and even *m.*

*Table 14*

| m= $2k+1$ | $d_i$-knives | SS | $SS_C$ | VS | $VS_C$ | VV | $VV_C$ | |
|---|---|---|---|---|---|---|---|---|
| $d_m$= $2k+1$ | m-knives | j·m+1 where j=1,2,…,k+1 | 2·m | j·m+1 where j=2,…,k+1 | 2·m | j·m+1 where j=2,…,k | 0 | |
| $d_i$ = $2ki+1$ | $d_i$ - knives | j· $d_i$+1 where j=1,2,…,$k_i$+1 | 2·$d_i$ | j· $d_i$+1 where j=1,2,…,$k_i$+1 | 2·$d_i$ | j· $d_i$+1 where j=1,2,…,$k_i$+1 , | 0 | i=1,2,…,N |
| $d_1$ =1 | 1-knife | 2 | 2 | 2 | 2 | 2, except m=3 | 2 | |

| m= 2k | $d_i$-knives | SS | $SS_C$ | VS | $VS_C$ | VV | $VV_C$ | remarks |
|---|---|---|---|---|---|---|---|---|
| $d_m$= $2k$ | m-knives | j·m+1 where j=1,2,…,k | 0 | j·m+1 where j=2,…,k | 0 | j·m+1 where j=2,…,k-1 | 0 | |
| $d_i$ =m/2 $2k_1+1$ | m/2-knives | j· $d_1$+1 where j=1,2,…,$k_1$+1 | 2·$d_1$= m | j· $d_1$+1 where j=1,2,…,$k_1$+1 | 2·$d_1$= m | j· $d_1$+1 where j=1,2,…,$k_1$ | 0 | m/2 is odd |
| $d_1$=m/2=2 $k_1$ | | j· $d_i$+1 where j=1,2,…,$k_1$ | | j· $d_i$+1 where j=1,2,…,$k_1$ | | j· $d_i$+1 where j=1,2,…,$k_1$ | | m/2 is even |
| $d_i$ = $2ki+1$ | $d_i$ - knives | j· $d_i$+1 where j=1,2,…,$k_i$+1 | 2·$d_i$ | j· $d_i$+1 where j=1,2,…,$k_i$+1 | | j· $d_i$+1 where j=1,2,…,$k_i$+1 | 2·$d_i$ | $d_i$ is odd |
| $d_i$ = $2k_1$ | | j· $d_i$+1 where j=1,2,…,$k_i$ | 2·$d_i$* | j· $d_i$+1 where j=1,2,…,$k_i$ | | j· $d_i$+1 where j=1,2,…,$k_i$ | 2·$d_i$ * | $d_i$ is even |
| $d_1$ =1 | 1-knives | 2 | 2 | 2 | 2 | 2 | 2 | |

2·$d_i$* - When $d_i$* is an even number, the $d_i$ -knives turn into $\frac{d_i}{2}$ - knives and therefore the result is as if it were acting similar to the $\frac{d_i}{2}$ knives.

## 5.4 Möbius phenomena

**Theorem 4**: In the cutting of *GML* bodies with chordal knives, the Möbius phenomenon with one resulting body and link number 1 can appear only for *m* even and when the knife cuts through the centre. In the cutting of *GML* bodies the Möbius phenomenon can appear for both *m* odd and *m,* when the knife is a radial knife, a ray starting at the centre of the polygon.

The radial knife is actually the position vector. This acts as the hand of a clock with discrete ticks (Figure 45) but with continuous movement within a $GML_m^n$ body, keeping a fixed direction from centre to the vertex. During this continuous movement, the various vertices of the polygon coincide with the upper circle in a torus in discrete steps. The tip of the position vector traces out a toroidal line, wound around the torus, now in a continuous way.

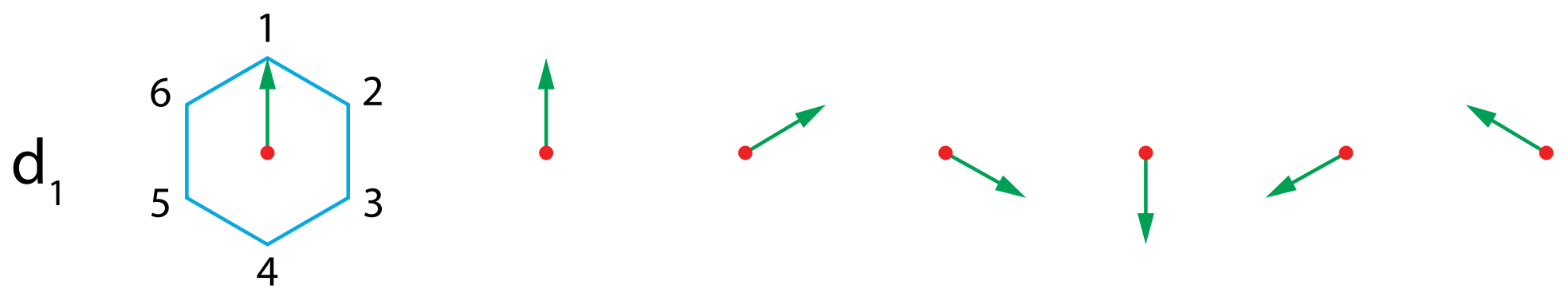

*Figure 45: Radial knife or position vector*

# 6 Conclusio

## 6.1 Möbius before and after Möbius

The Möbius strip is an icon of mathematics. It was named after the German mathematician Möbius, but known already 2 millenia ago (Figure 46). This "twisted cylinder" is obtained by giving one end of a rectangular strip of paper a twist before joining both ends. Simple as this may seem Möbius strips are one-sided, non-orientable surfaces with some counterintuitive properties. Tracing out a path on a Möbius strip will show that it has only one side and the whole strip can be painted in one color without lifting the paint brush. If some form of identification is done along the central line before joining, one will find that orientations are reversed after rejoining. Johann-Benedict Listing (1808–1882) and August Ferdinand Möbius (1790–1868) discovered these remarkable objects in the mid-nineteenth century. The name Topology - the study of shape characteristics independent of measuring - was used first by Listing. While topology is a branch of geometry, Listing chose the name topology to avoid reference to measurements and referred to "Modalität und Quantität" (modality and quantity) to distinguish between shape and magnitude.

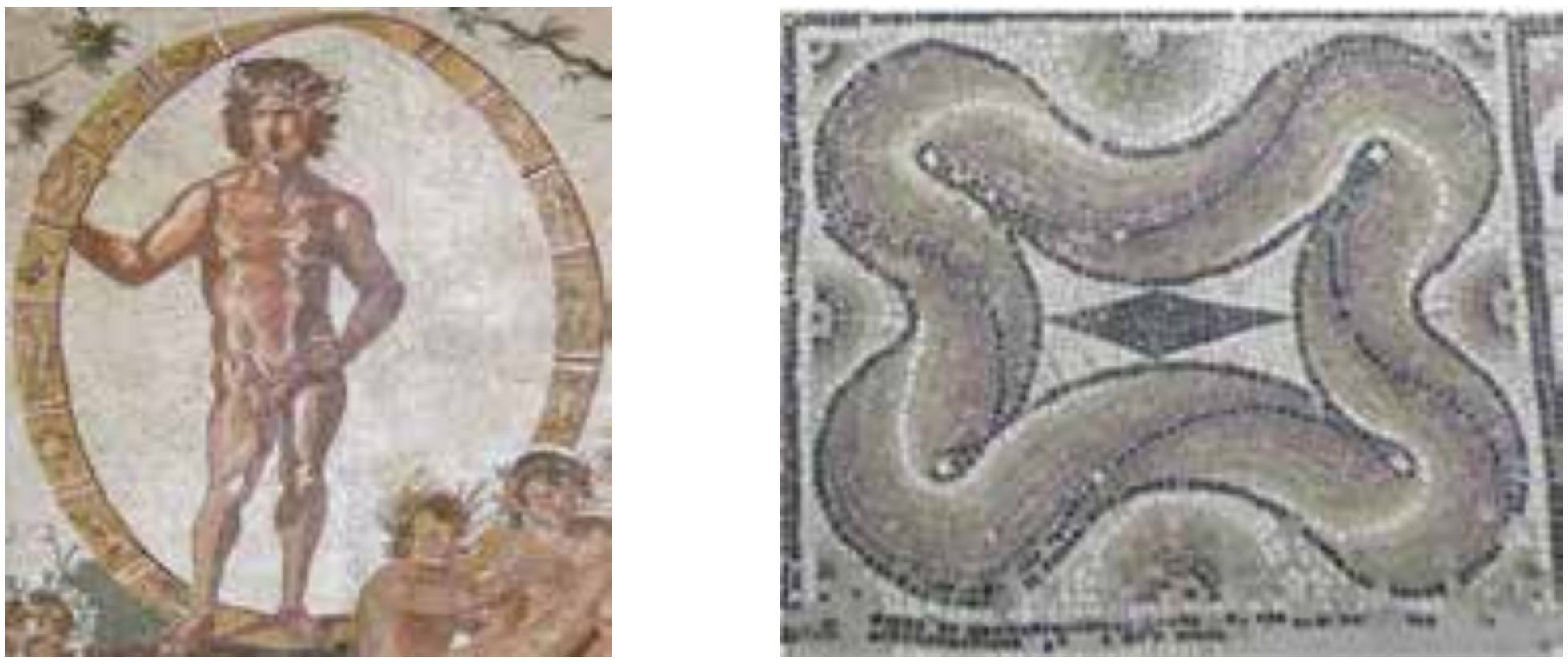

**a** **b**

*Figure 46: **a**. Möbius before Möbius [31] and **b.** GML's before the present*

In *On Growth and Form*, D'Arcy Thompson [32] summarized the field and its initial historical development: "*...in this study of a segmenting egg we are on the verge of a subject adumbrated by Leibniz, studied more deeply by Euler, and greatly developed of recent years. It is the* Geometria situs *of Gauss, the* Analysis situs *of Riemann, the* Theory of Partitions *of Cayley, of* Spatial Complexes *or* Topology *of Johann Benedict Listing... Leibniz has pointed out there was room for an analysis of mere position, apart from magnitude:* "Je croy qu'il nous faut encor une autre analyse, qui nous exprime directement *situm*, comme l'Algèbre exprime *magnitudinem*". .. *Leibniz used it to explain the game of solitaire, Euler to explain the knight's move on the chessboard or the routes over the bridges of a town. Vandermonde created a* géometrie de tissage, *which Leibniz himself had foreseen, to describe the intricate complexity of interwoven threads in a satin or brocade. Listing, in a famous paper, admired by Maxwell, Cayley and Tait, gave a new name to this new "algorithm", and shewed its application to the curvature of a twining stem or tendril, the aestivation of a flower, the spiral of a snail-shell, the scales on a fir-* *cone, and many common things. The theory of "spatial complexes", as illustrated especially by knots, is a large part of the subject*". There are many things to which the new Geometria Situs could be applied.

## 6.2 Geometrical topology

*Analysis situs* or topology developed over many centuries, and has developed into highly abstract mathematics at one extreme and into various educational tools for displaying the magic of mathematics at the other. Möbius bands are icons of mathematics [33], and so are knots and links. It is also linked to geometry, in particular, the inner geometry of surfaces, starting with Gauss' *Theorema Egregium*. Indeed, a torus can be constructed topologically by folding a piece of paper into a cylinder and connecting the ends. For this reason the intrinsic geometry of the torus and the piece of paper (exemplifying part of the plane) is the same, and hence the torus is considered as a flat surface (in contrast to a sphere). In similar ways Möbius bands, Klein Bottles and Projective planes can be constructed by taking a piece of paper, identifying opposite sides and connecting those in the appropriate way. The projective plane is a disk to which a Möbius band is glued.

Nevertheless, these are recipes or instructions to generate surfaces, and one of the main lacunae in the geometry / topology distinction is the absence of general but simple analytic representations. Even the whole field of knots still only leads to tabulations, not classifications. Our analytic representations or equations substitute for "recipes" or "(computer) algorithms" to "generate" Möbius strips, tori, helices or more complex shapes. We no longer need all types of different parametrizations for different shapes, but we can combine these shapes now into single equations, allowing for continuous transformations from circle or sphere into any other shape, including knots and certain polyhedrons, irrespective of topological constraints.

In this sense transformations from sphere to torus e.g. with intermediates as ring, horn and spindle tori, or to knots are very natural, and continuous. The definition of a torus as $S^1 \times S^1$, or fibre bundles with unit intervals as fibre over manifolds [34], or processes such as making holes and using surgery or pinching, can be approached from another, continuous geometric perspective. In this perspective shapes can be transformed into each other and cut into predictable substructures. The reduction of "recipes" to equations, with a concomitant reduction of shape complexity, has always led to progress in science. This "old kind of science" has been extremely valuable for mankind, from Greek mathematics to the present day. By focusing on some very elementary ideas, such as the Pythagorean structure of Gielis curves and Gaspard Monge's idea of describing complex movements as superpositions of elementary ones for *GML* [35], highly complex structures and very many of the well-known shapes in mathematics can be described in very simple analytic representations. There are one-to-one mappings from one shape to the other, and these mappings are geometrical, at times violating the rules of topology, but at the same time bridging topology and geometry.

In our work we focused so far on the 'forward' approach, but we believe that uniqueness results can be proven, so that from observing and counting complex structures it should be possible to reconstruct the original *GML* body or surfaces from which it was derived.

## 6.3 The general case of cutting GML

The generalization of Möbius ribbons to Generalized Möbius-Listing bodies and surfaces was initiated in the late nineties of last century [1]. The motivation was the fact that the solution of partial differential equations strongly depends on the topological properties of the domain on which the problem is considered. The first results were obtained on Möbius strips, giving a general solution for cutting of ribbons for any number of twisting and for any number of cuttings. Consequently the description was broadened to include surfaces and bodies and in the recent decade, classifications were also obtained for $GML$ bodies and surfaces with convex cross-sections for symmetries $m = 2, 3, 4, 5, 6$.

These developments became possible by integrating $GML$ and Gielis curves and bodies, a geometric transformation describing a very wide range of abstract and natural shapes. This integration made it possible to develop strategies to achieve a general classification for cutting of $GML$ bodies of any symmetry, starting from cross sections that are $m$-regular polygons. The original problem could be reduced to the problem of cutting $m$-regular polygons with crossing diagonals or knives. With the proposed methods even the topological characteristics of the resulting bodies and surfaces and their link number can be found with the 2D strategy of cutting polygons.

Determining the number of variants after cutting has both geometrical and topological solutions, resulting in Theorems 1 and 2. These theorems not only give the precise number of variants after cutting, but also allow for determining the precise shape of the resulting bodies and their topology. Annex 1 gives all relevant graphics from $m = 6$ to $m = 10$. Theorem 3 gives the number of different objects resulting from specific ways of cutting, and from this one can count the connected zones in polygons or the link number after cutting of $GTR_m^n$ bodies. It generalizes the cutting of polygons and partitions in the sense of Euler and Cayley, with $VV$ cuts and diagonals. Because of these connections, combinatorial or other approaches could have been used, but the strategy of using planar geometry (as a reduction of a 3D problem), has the advantage of revealing many connections and a broad scope of applications, and the additional advantage of having a strong educational component.

Theorem 4 gives the precise conditions under which the Möbius phenomenon will occur for regular polygons with both odd and even symmetry $m$, when using radial knives. The radial knife is nothing but the position vector, which is static if defined as $m$-polygons, but dynamic using Gielis transformations, allowing to morph the polygon into circles, starfish and so on. The tip of the knife traces out a toroidal line on the circumscribing torus, whereas the knife itself describes a surface. When the $GML_m^n$ body is not closed, a generalized cylinder results. The generalized cylinder can also be a generalized cone. The tip of the arrow traces out a helix, but with $m$-symmetry. The arrow itself than traces out a surface, in particular a Riemann surface. If the arrow has a given depth (thickness in the direction of the cylinder), it will trace out a body or helical shell. If the arrow has a given depth in this case, it will trace out a shell body as observed in mollusc shells.

## 6.4 Future work

The results are obtained for regular polygons, but can easily be generalized. First, using Gielis curves and curved knives greatly generalize these results, for convex and concave domains. Second, knives can be the result of some optimization problem in convex geometry. Third, only one knife in one direction was used, but $k$-knives [3] with multiple parallel blades can be used to make multiple parallel cuts at once; in the extreme the knife can be a Cantor-knife, such that the resulting cut gives in cross section the Cantor continuum. Fourth, the definition of $d_i$-knives starts from $d_1$-knife by $i$ rotations with the centre of the polygon as centre of rotation, but different movements of $d_1$-knives will give other results. One example is the Cremona construction of the cardioid, as envelope of a pencil of lines. This is realized with $VV$ cuts of the type $VV_{i,2i}$ connecting $n$ equally spaced points on the circle.

These results open up many new possible areas of research, in particular for a unified and general geometric framework for the study of natural phenomena in biology, chemistry and physics, both in the small and in the large. It also opens up many challenges.

In physical realisations Möbius bands will take on shapes, which are physically dependent, among others, on width, length and material [36]. Within the more general framework of $GML_m^n$ or cylinders with non-circular cross sections, the situation is very different. In the case of cylinders subjected to torsion, the cross sections are known to remain planar. The research of Adhémar Jean Claude Barré de Saint-Venant (1797-1886) on the torsion of prisms and non-circular rods [37] showed how the internal body is subjected to stress resulting in distorted, anticlastic cross sections (Figure 47**a,b**). The lines that remain in-plane are $VV$ and $SS_C$ cuts perpendicular to the sides for the square or $VS_C$ cuts for the triangle. This phenomenon was observed in the horns of rams (Figure 47**c**) [32]. The stresses and strains and the displacements (Figure 47) connect shape and differential equations. Elasticity theory then enters the field, a wide research domain of elastic and prismatic cusped shells, plates and beams [38], (Figure 1).

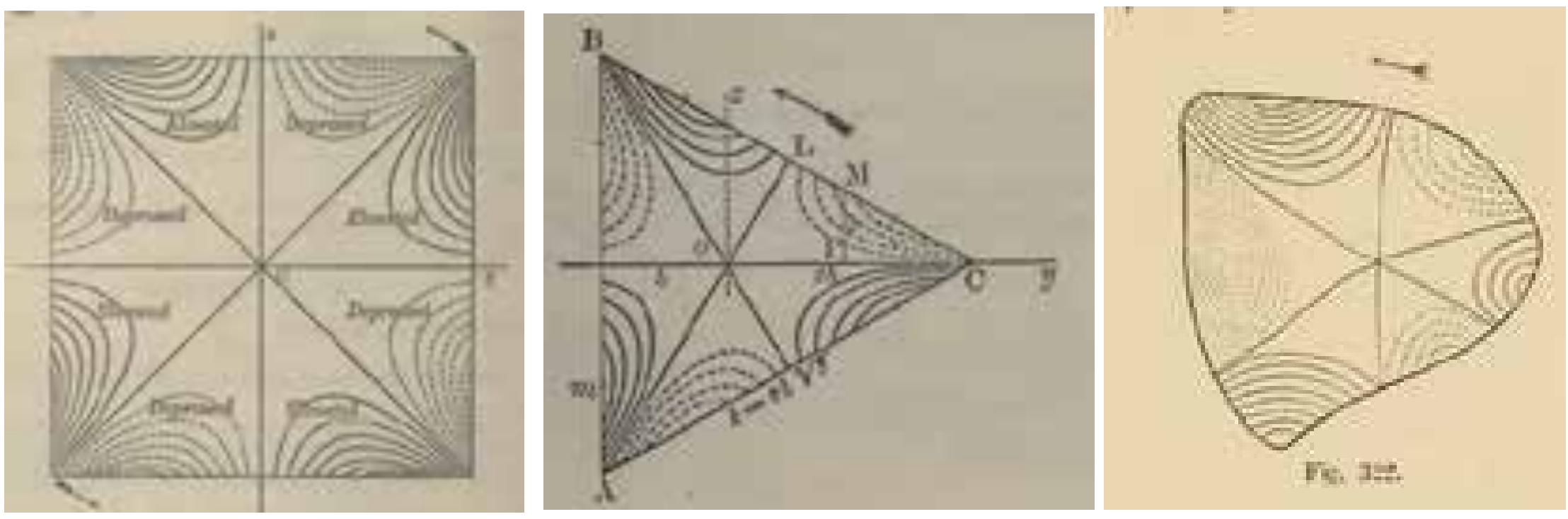

**a** **b** **c**

*Figure 47: Torsion of prisms with square (**a**) and triangular (**b**) cross-sections. Dotted lines are depressed regions and full lines are elevated parts of the cross sections.* **c.** Anticlastic surfaces in the horns of rams [32]

Precise knowledge of the domains is essential for solving boundary value problems and this was precisely the motivation for the study of $GML_m^n$ surfaces and bodies (in particular, the starting point was the study of Saint-Venant's Principle on complex domains). The same motivation underlies the development of $R$-functions [11; 22].

Our methods allow for the precise description of a wide variety of domains in a uniform geometric way, based on $GML/GTR$ and plane geometry. With Gielis transformations continuous transformations among domains are added. This leads to naturally adapted coordinate systems on which boundary value problems can be solved, which is the essence of mathematical physics, according to Gabriel Lamé, who was also one of the main contributors to elasticity theory [37].

# 8 Annex 1

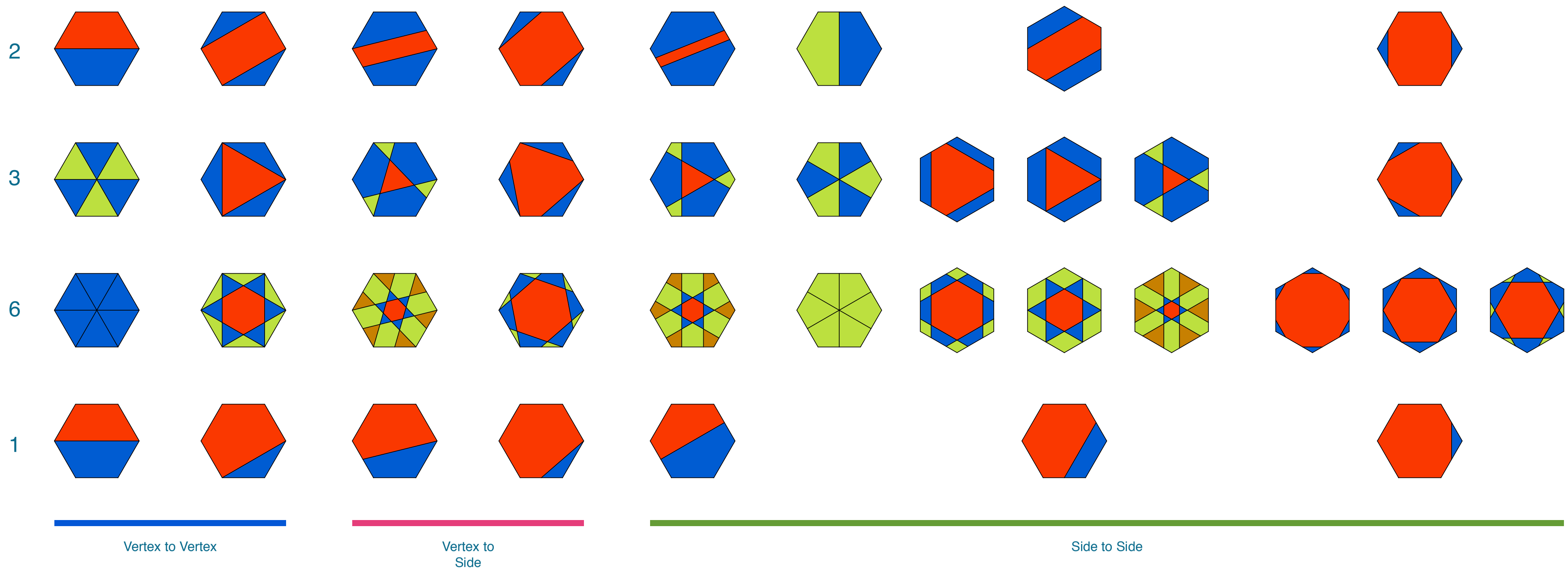

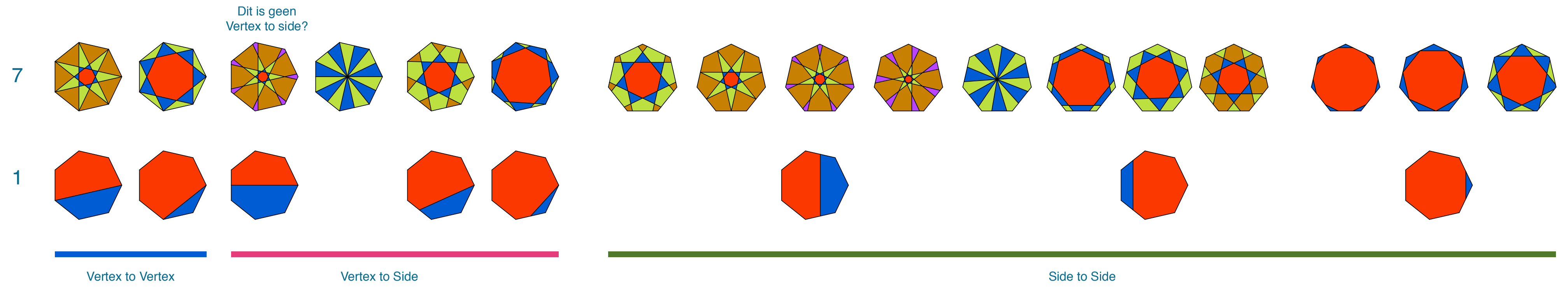

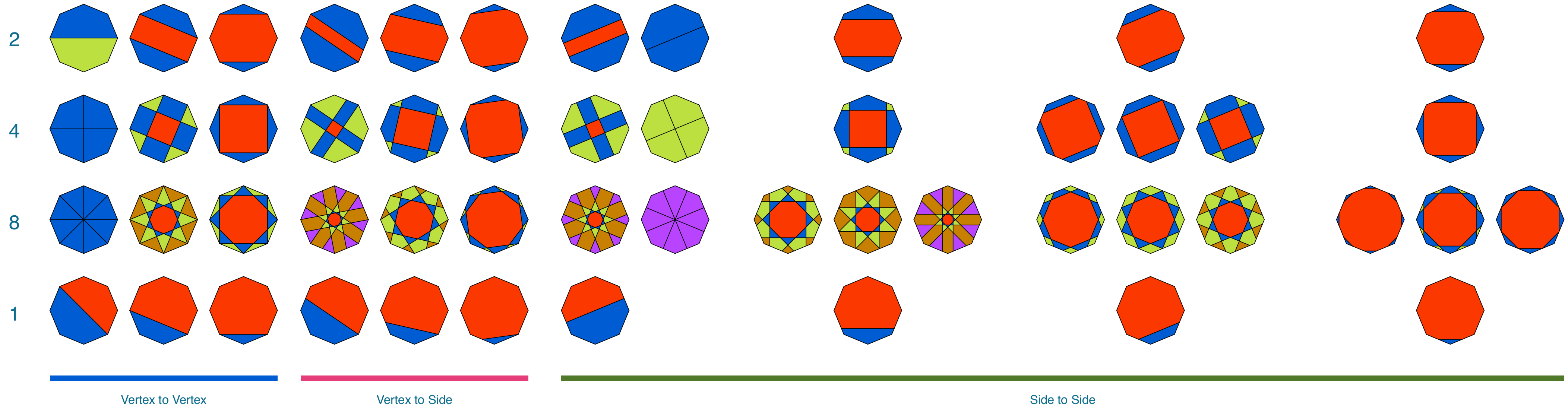

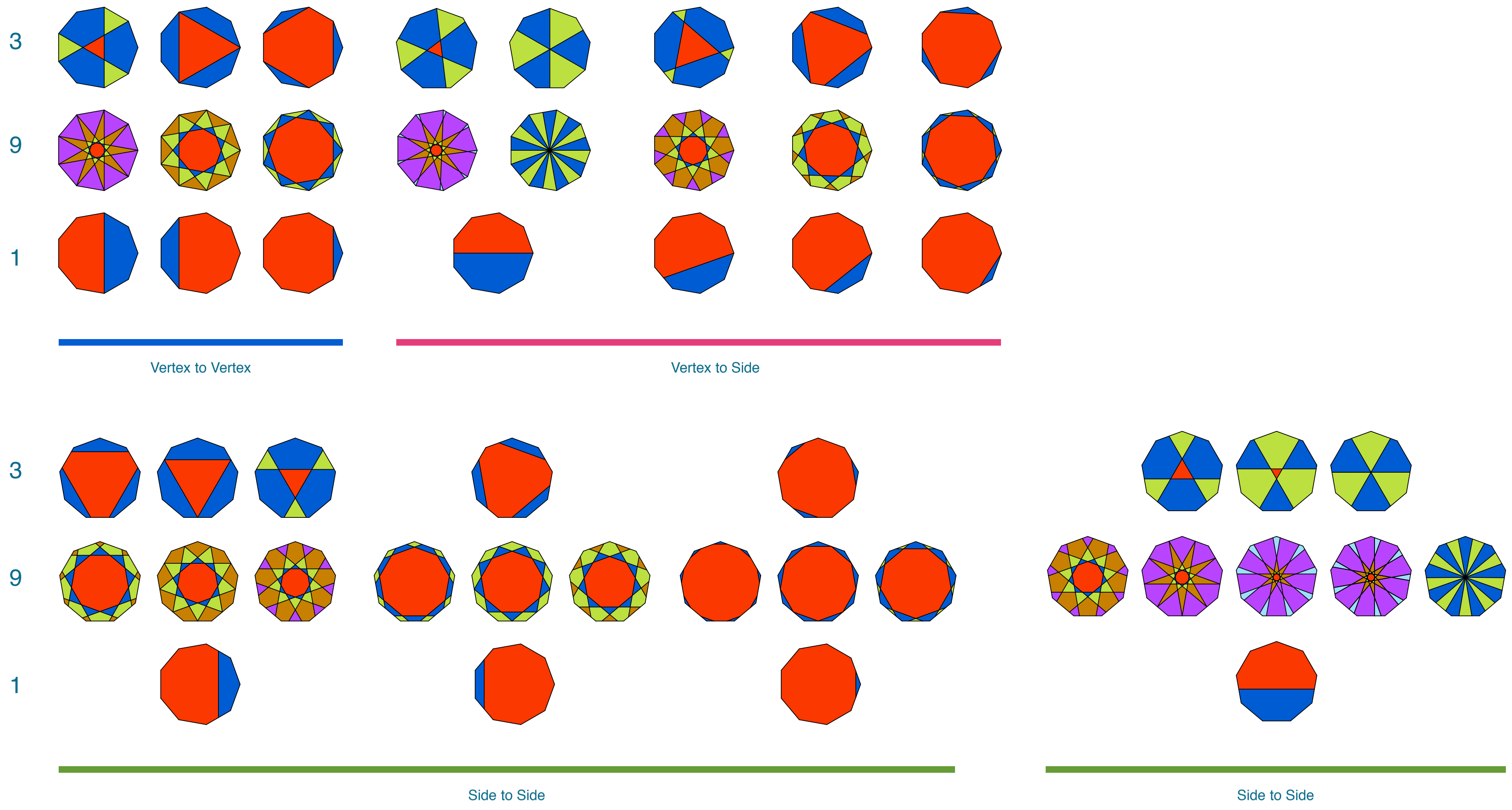

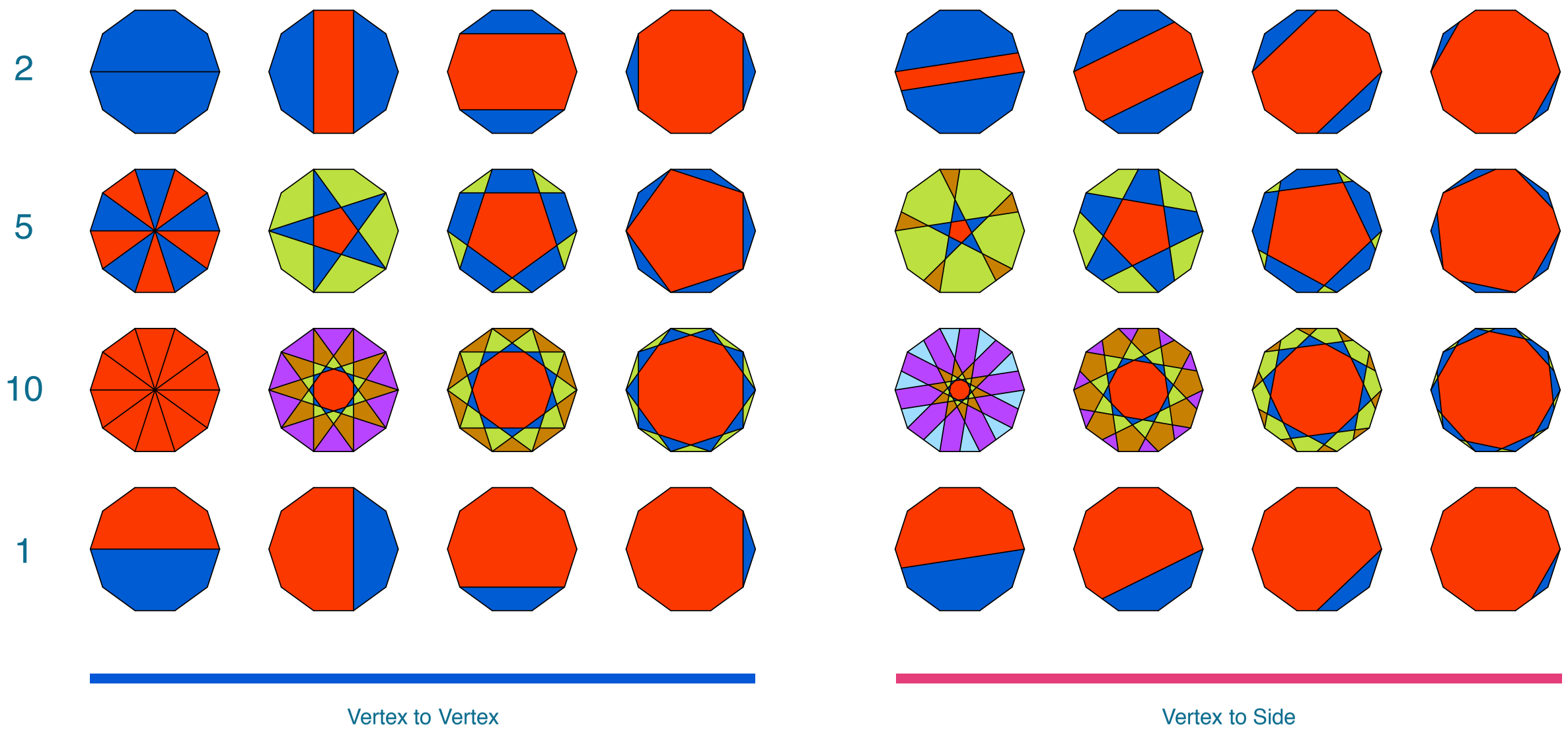

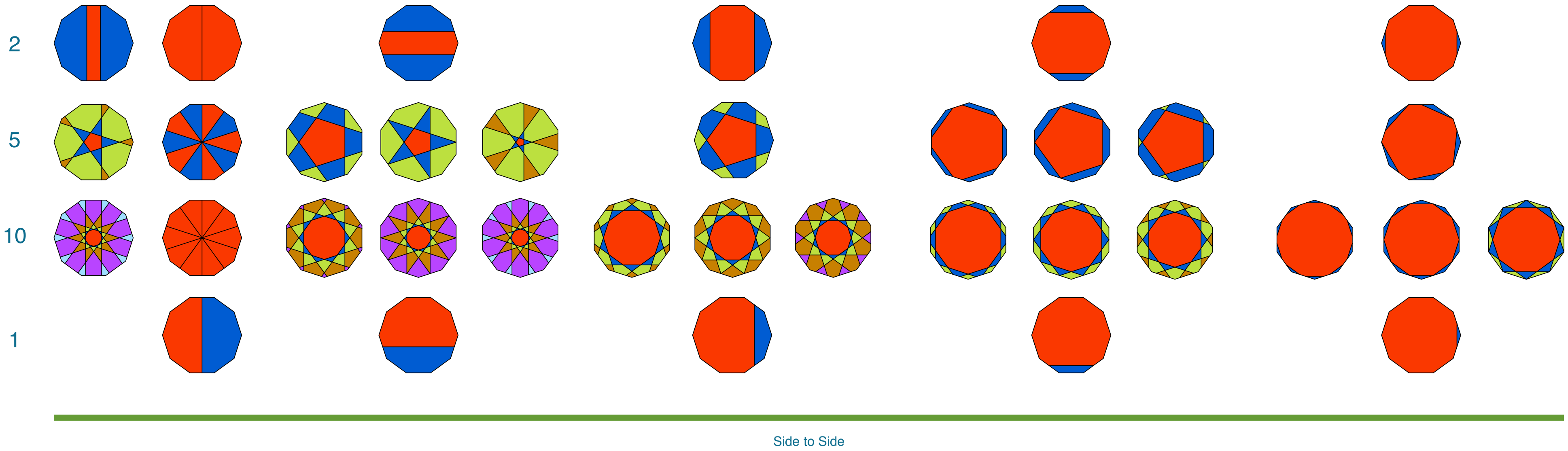
2
5
10
1
Side to Side